\documentclass[letterpaper]{article}
\usepackage{authblk}
\usepackage{url}
\usepackage[margin=1in,footskip=0.25in]{geometry}
\usepackage[bbgreekl]{mathbbol}
\usepackage{amssymb}             % AMS Math
\DeclareSymbolFontAlphabet{\amsmathbb}{AMSb}%
\usepackage{eucal}
\usepackage{graphicx}
\usepackage{bm,upgreek}
\usepackage[table]{xcolor}% http://ctan.org/pkg/xcolor
\usepackage{slashbox}
\usepackage[symbol,bottom]{footmisc}

\usepackage{url}
\usepackage{amsmath}
\usepackage[linesnumbered,ruled]{algorithm2e}
\usepackage{algpseudocode}
\usepackage{float}
\usepackage{lipsum}
\usepackage{subfigure}
\usepackage{wrapfig}
\usepackage{color}
\usepackage{easylist}
\usepackage{lineno}
\usepackage{stmaryrd}
\usepackage{booktabs}
\usepackage{cleveref}
\usepackage{multirow}
\usepackage{hhline}
\usepackage{filecontents}
\usepackage{enumitem}
\usepackage{cite}
\newlist{steps}{enumerate}{1}
\setlist[steps, 1]{label = Step \arabic*:}
\usepackage{tikz}
\usepackage{pifont}

\newcommand{\xmark}{\text{\ding{55}}}
\def \grad{\nabla}
\def \half{\frac{1}{2}}

\newcommand{\volume}{{\ooalign{\hfil$V$\hfil\cr\kern0.08em--\hfil\cr}}}

\usepackage{blindtext}

\def \CC{\bbsigma_\text{C}}

\def \CC{\mathbb{C}}
\def \FF{\mathbb{F}}
\def \II{\mathbb{I}}
\def \PP{\mathbb{P}^\text{s}}

\def \sigmaf{\bbsigma^\text{f}}

\def \Ncycle{N_\text{cycle}} 
\def \F{\vec{F}}

\def \N{\vec{N}}

\def \U{\vec{U}}

\def \X{\vec{\chi}}
\def \Y{\vec{\xi}}

\def \s{\vec{X}}

\def \u{\vec{u}}

\def \x{\vec{x}}

\def \P{\varphi}
\def \grad{\nabla}
\def \half{\frac{1}{2}}

\def \fs{\text{fs}}

\def \Omegas{\Omega^\text{s}}

\def \DA{{\mathrm d} A}

\def \rhos{\rho^{\text{s}}}

\def \Ks{\kappa^{\text{s}}}

\def \rhof{\rho^{\text{f}}}
\def \muf{\mu^{\text{f}}}
\def \cL{\mathcal L}
\def \cE{\mathcal E}

\def \cJ{\vec{\mathcal J}}

\def \cS{\vec{\mathcal S}}
\def \cT{\mathcal T}

\def \dtf{\Delta t^\text{f}}
\def \dt{\Delta t}
\def \dts{\Delta t^\text{s}}

\def \cS{\mathcal S}
\def \cT{\mathcal T}

\def \mfac{M_\text{fac}}

\def \Gs{G^\text{s}}
\def \Ks{K^\text{s}}

\def \St{\text{St}}

\def \L2{L^2}
\def \Linf{L^{\infty}}

\def \grad{\nabla}
\def \half{\frac{1}{2}}

\def \Mfac {M_{\text{fac}}}

\renewcommand{\vec}[1]{\ensuremath\boldsymbol{#1}}
\newcommand{\norm}[1]{\left|\left|#1\right|\right|}

\title{A sharp interface Lagrangian-Eulerian method for flexible-body fluid-structure interaction}

\author[1]{Ebrahim M.~Kolahdouz} %\corref{cor1}}
\author[2]{David R. Wells}
\author[2]{Simone Rossi}
\author[3]{Kenneth I. Aycock}
\author[3]{Brent A.~Craven}
\author[4,5,6,7]{Boyce E.~Griffith} %\corref{cor2}}
\affil[1]{Center for Computational Biology, Flatiron Institute, Simons Foundation, New York, NY, USA}
\affil[2]{Department of Mathematics, University of North Carolina, Chapel Hill, NC, USA}
\affil[3]{Division of Applied Mechanics, Office of Science and Engineering Laboratories, Center for Devices and Radiological Health, United States Food and Drug Administration, Silver Spring, MD, USA}
\affil[4]{Departments of Mathematics and Biomedical Engineering, University of North Carolina, Chapel Hill, NC, USA}
\affil[5]{Carolina Center for Interdisciplinary Applied Mathematics, University of North Carolina, Chapel Hill, NC, USA}
\affil[6]{Computational Medicine Program, University of North Carolina, Chapel Hill, NC, USA}
\affil[7]{McAllister Heart Institute, University of North Carolina, Chapel Hill, NC, USA}
\affil[ ]{\texttt{ekolahdouz@flatironinstitute.org} and \texttt{boyceg@email.unc.edu}}
\begin{document}

\maketitle

\begin{abstract}
This paper introduces a sharp-interface approach to simulating 
fluid-structure interaction (FSI) involving 
flexible bodies described by general nonlinear material models and across a broad range of mass density ratios.
This new flexible-body immersed Lagrangian-Eulerian (ILE) scheme 
extends our prior work on integrating partitioned and immersed approaches to rigid-body FSI. Our numerical approach
incorporates the geometrical and domain solution flexibility
of the immersed boundary (IB) method with an accuracy comparable to body-fitted approaches 
that sharply resolve flows and stresses up to the fluid-structure interface.
Unlike many IB methods, our ILE formulation uses distinct momentum equations for the fluid and solid subregions 
with a Dirichlet-Neumann coupling strategy that connects fluid and solid subproblems through simple interface conditions. 
As in earlier work, we use approximate Lagrange multiplier forces to treat the kinematic interface conditions along the fluid-structure interface. 
This penalty approach simplifies the linear solvers needed by our formulation by introducing two representations of the fluid-structure interface, 
one that moves with the fluid and another 
that moves with the structure, that are connected by stiff springs.
This approach also enables the use of multi-rate time stepping, which allows us to use different time step sizes for the fluid and structure subproblems.
Our fluid solver relies on an immersed interface method (IIM) for discrete surfaces 
to impose stress jump
conditions along complex interfaces while enabling the use of fast structured-grid solvers for the incompressible Navier-Stokes equations. 
The dynamics of the volumetric structural mesh are determined using a standard finite element approach to
large-deformation nonlinear elasticity via a nearly incompressible solid mechanics formulation.
This formulation also readily accommodates 
compressible structures with a constant total volume, and it can handle fully
compressible solid structures for cases in which at least part of the solid
boundary does not contact the incompressible fluid.
Selected grid convergence studies demonstrate
  second-order convergence in volume conservation and in the pointwise discrepancies between corresponding positions of the two interface representations
   as well as between first and second-order convergence in the structural displacements.
The time stepping scheme is also demonstrated to yield second-order convergence.
 To assess and validate the robustness and accuracy of the new algorithm,
comparisons are made with computational and experimental FSI benchmarks.
Test cases include both smooth and sharp geometries in various flow conditions. 
We also demonstrate the capabilities of this methodology by applying it to model the transport and capture of a 
geometrically realistic, deformable blood clot in an
inferior vena cava filter.
\end{abstract}


%%%%%%%%%%%%%%%%%%%%%%%%%%%%%%%%%%%%%%%%%%%%%%%%%%%%%%%%%%%%%%%%%%%%%%%%%%%%%%%%%%%%%%%%%%%%%%%%%%%%%%%%%%%%%%%%%%%%%%%%%%%
\section{Introduction}
\label{sec:intro}

Computer modeling and simulation 
are powerful tools for analyzing nonlinear fluid-structure interaction (FSI)
 that can include large deformations and displacements of flexible bodies immersed in fluid. 
 Well-known approaches to FSI include \emph{partitioned formulations} that use separate descriptions of the fluid and solid subregions.
One common discretization approach for such formulations is to use non-overlapping, body-fitted domain discretizations \cite{nobile1999stability,farhat2001discrete}.
Body-fitted discretizations facilitate the specification of precise boundary conditions, which can resolve flow features up to the fluid-structure interface.
However, body-fitted approaches require complex and potentially expensive mesh manipulations via grid generation, regeneration, and mesh morphing \cite{cottrell2009isogeometric}.
Further, grid regeneration generally will require interpolating the computed solution to a new mesh which  \cite{degand2002three, luke2012fast},
can introduce errors and instabilities.
\emph{Immersed formulations}, such as Peskin's immersed boundary (IB) method \cite{peskin1972flow,peskin1977,peskin2002immersed}, 
are alternatives to body-fitted approaches that avoid the need to use conforming discretizations of the fluid and structure. 
Immersed methods thereby can circumvent the mesh generation difficulties of body-fitted approaches and
facilitate very large structural deformations.
Immersed methods 
also can simplify developing efficient structured-grid solvers. 
The key challenge in these methods, however, is 
related to developing effective and accurate coupling operators and boundary conditions that link the Eulerian and Lagrangian variables; 
see a recent review by Griffith and Patankar \cite{griffith2020immersed}
along with those by Peskin \cite{peskin2002immersed}, Mittal and Iaccarino \cite{mittal2005immersed}, and Sotiropoulos and Yang \cite{sotiropoulos2014immersed}.
This paper extends our earlier work to integrate partitioned and immersed approaches to FSI by generalizing our recently developed 
immersed Lagrangian-Eulerian formulation for fluid-rigid structure interaction  \cite{kolahdouz2021sharp} to problems 
involving flexible structures described by nonlinear elastic material models.
To do so, in this work, the equations of motion for the 
rigid body are replaced 
by the equations of elastodynamics. 
We introduce mathematical formulations for both incompressible and compressible structures. Numerical methods for the incompressible case adopt a 
nearly incompressible formulation 
as a penalty method for exactly incompressible structural models.

The immersed boundary method was created to treat FSI problems involving thin flexible 
interfaces. Subsequent extensions were developed to treat 
structures of finite thickness, 
most of which use 
integral transforms with regularized Dirac delta function kernels to connect the Eulerian and Lagrangian 
variables \cite{mcqueen2000three,mcqueen2001heart,peskin2002immersed,griffith2005order,griffith2007adaptive}. 
In the original IB method, 
fiber-based elasticity models describe the elasticity of the immersed structure, and regularized 
Dirac delta functions connect the Eulerian and Lagrangian variables.
These regularized delta functions prolong structural forces to the Cartesian grid and restrict Cartesian grid velocities onto the Lagrangian mesh \cite{peskin1972flow}.
The immersed finite element (IFE) method is an extension of the IB method that uses finite element approximations for both the Lagrangian and the Eulerian equations
with a class of kernel functions based on reproducing kernel particle methods (RKPM) for the fluid-solid coupling \cite{zhang2004immersed,zhang2007immersed}. 
Boffi et al.~\cite{boffi2008hyper} presented a fully variational IB formulation that allows for a nonlinear mechanics-based description of the 
Lagrangian structure within the framework of large-deformation continuum mechanics.
Regularization in this formulation is achieved implicitly through the finite element (FE) basis functions.
The immersed structural potential method of Gil et al.~\cite{gil2010immersed,gil2013enhanced} 
uses a meshless method for the description of hyperelastic immersed structure
 and a family of spline-based kernel functions for interaction between fluid and solid.
Devendran and Peskin \cite{devendran2012immersed} developed an energy-based IB method 
 that allows for a nodal
approximation of elastic forces directly calculated from energy functional without the use of stress tensors,  
though in the process, analytical calculations of the strain energy functional derivatives
with respect to the coefficients of the discretized displacement field in the FE space are required.  
Griffith and Luo introduced a hybrid finite difference/finite element 
 approach that uses a FE discretization of the solid structure while retaining a finite difference scheme for the Eulerian variables
  on the Cartesian grid \cite{BEGriffith17-ibfe}. As in the conventional IB method,
   a regularized delta function is used in approximations to the integral transforms, 
  but the nodal velocity of the immersed structure is obtained by projecting 
  an intermediate Lagrangian velocity field 
  onto the finite element space.
  Recent work demonstrated that this method is equivalent to the method of Devendran and Peskin for specific choices of 
  FE basis functions and quadrature rules \cite{wells2021nodal}. 
  Mortar methods have also been used to treat interface conditions in immersed formulations \cite{baaijens2001fictitious,hesch2014mortar,nestola2019immersed}.
  In such approaches, information transfer between fluid and solid subdomains is done  
  in a weak sense, by using an $L^2$-projection operator.
   Methods based on distributed Lagrange multipliers (DLM) form another class of immersed formulations,
and since their original development by Glowinski et al. \cite{glowinski1994fictitious,glowinski1999distributed},
 different variations of this approach have been reported for fluid-flexible~structure interaction \cite{van2004combined,yu2005dlm,bhalla2013unified,boffi2015finite,kadapa2016fictitious}.
 To our knowledge, except for our earlier work on rigid-body FSI \cite{kolahdouz2021sharp}, all existing DLM 
 formulations use volumetric coupling schemes to connect the fluid and solid variables.
  Other immersed formulations for FSI include the moving least square direct forcing immersed 
  boundary method \cite{de2016moving,vanella2020direct,spandan2018fast}
  and fully Eulerian approaches based on finite difference \cite{sugiyama2011full, valkov2015eulerian,rycroft2020reference}
  or finite element \cite{richter2013fully,richter2010finite} techniques.
  In all of these approaches, the equations of fluid dynamics are solved throughout the computational domain, including in the solid subregion.
  Formulations that use a Cartesian grid for the fluid domain but do not extend the fluid solution to the interior 
  of the solid domain 
  have been developed based on the
   extended finite element method (XFEM) \cite{mayer20103d,burman2014unfitted,zonca2018unfitted}, 
     overlapping Chimera-like methods \cite{wall2006large,banks2012deforming, miller2014overset, serino2019stable},
   cut-cell methods \cite{mittal2008versatile}, and the cut-finite element method \cite{massing2015nitsche,schott2019monolithic}.
   These methods can be considered immersed in that they use a fixed Eulerian grid,  
   although to increase accuracy
   near the fluid-structure interface, they typically adopt approaches such as
   local modifications to the finite difference stencils 
   that make them similar to body-fitted discretizations. 

As in our earlier work \cite{kolahdouz2021sharp}, our goal is to create a numerical scheme with the geometrical and domain solution flexibility
of the IB method and an accuracy comparable to body-fitted approaches
for which sharp resolution of flows and stresses is achieved up to the fluid-structure interface.
Unlike many IB methods, our ILE approach uses distinct momentum equations for the fluid and solid subregions, and
it uses a Dirichlet-Neumann coupling strategy that connects fluid and solid subproblems through interface conditions. 
Our mathematical formulation can treat exactly incompressible structures, as in conventional IB methods. Unlike standard IB methods, however, 
it also can readily handle compressible structures with a constant total volume. Further, it can model fully compressible structures that are only partially immersed in fluid, 
such as if the elastic structure fully encloses an incompressible fluid.
To avoid solving monolithic systems of equations that involve both fluid and solid variables, however, we use two distinct 
representations of the fluid-structure 
interface: one that moves with the fluid and a second that moves with the structure.
 These two representations are connected by 
 approximate Lagrange multiplier forces that account for the kinematic interface conditions. 
Results demonstrate that this simple approach can control discrepancies between the two interface configurations effectively.
The deformations of the structure are determined 
as a result of interplay between the forces of intrinsic stresses associated with the solid constitutive model 
and the exterior fluid traction forces that are imposed as boundary conditions to the solid domain.
 Consequently, our approach can be viewed as a DLM scheme that uses interfacial, rather than volumetric, FSI coupling.
The dynamics of the volumetric structural mesh are determined using a standard Lagrangian finite element method for
large-deformation nonlinear elasticity.
We adopt a nearly incompressible solid mechanics formulation for most of the test cases considered in this work.
We also demonstrate the method's 
ability to model fully compressible structures
 for cases in which part of the solid
boundary does not contact the fluid, such as when the elastic structure fully encloses an incompressible fluid.
Our fluid solver relies on an immersed interface method (IIM) for discrete surfaces \cite{kolahdouz2020immersed} 
to impose stress jump conditions along complex interfaces and allows for the use of fast structured-grid solvers for the incompressible Navier-Stokes equations. 
Previous studies using the IIM have primarily focused on interfaces with prescribed motion 
\cite{xu2006, xu2008immersed, le2006immersed,kolahdouz2020immersed} or 
on FSI models involving thin elastic membranes \cite{le2006immersed, tan2009immersed, tan2010level, thekkethil2019level,layton2009using}
with simple material descriptions. 
 Our IIM formulation enables us to evaluate exterior fluid traction along fluid-structure interfaces described by $C^0$ surface representations, which facilitates its application to 
 FE models of complex geometries and structural mechanics.
%~ Our IIM formulation enables us to evaluate exterior fluid traction along the fluid-structure interface 
%~ described by the $C^0$ continuous surface representation.

We use both computational and experimental flexible-body FSI benchmarks to assess and validate the robustness and accuracy of the new algorithm.
Test cases include two dimensional cases of 
a soft disk in a lid driven cavity \cite{zhao2008fixed,BEGriffith17-ibfe}, 
a version of the Turek-Hron benchmark problem \cite{turek2006proposal}, and 
the damped structural instability of a fully enclosed fluid reservoir \cite{kuttler2006solution}.
We also present results from three-dimensional cases, including a horizontal flexible plate inside a flow phantom
based on experimental results proposed by Hessenthaler et al. \cite{hessenthaler2017experiment}, 
and a demonstration case applying our methodology to model the transport and capture of a geometrically realistic, deformable blood clot in an
inferior vena cava (IVC) filter \cite{riley2021vitro}.
 IVC filters are cardiovascular devices that are implanted in the IVC, a large vein in the abdomen through which blood flow from the lower 
 extremities returns to the heart. IVC filters capture clots shed from the lower extremities, as can occur in deep vein thrombosis, before they can migrate 
 to the lungs and cause a potentially fatal pulmonary embolism.
 Modeling the transport and fluid-structure interaction of flexible blood clots 
in the venous vasculature is critical to computationally predicting the performance of embolic protection devices like IVC
blood clot filters. 

Numerical tests demonstrate that approximate Lagrange multiplier forces can be constructed to achieve
pointwise second-order convergence for the discrepancy between the positions of the two interface representations in our algorithm.
We also report
at least second-order accuracy in structural volume representing the incompressibility error of the Lagrangian domain.
Our time integration scheme is demonstrated to yield second-order accuracy for the structural displacement, and the overall scheme yields between 
first-~and second-order accuracy in a full spatio-temporal convergence study.
Additionally, our simulations demonstrate that structures with low and near-equal density ratios 
are successfully modeled with reasonable time step sizes 
 with no apparent restrictions caused by artificial added mass instabilities.
 Instabilities due to the artificial added mass effect have been previously
reported for weakly coupled FSI schemes with explicit time stepping schemes
for systems with low or nearly equal structure-to-fluid density ratios \cite{forster2007artificial,liu2014stable,young2012hybrid}. 
Though no theoretical proof is provided in this work, we hypothesize that the evident robustness of 
 the method to artificial added mass instabilities stems from its
consistent treatment of the momentum of the fluid near the fluid-structure interface.
%~ in the immersed formulation
%~ that avoids inducing temporal discontinuities in the fluid acceleration. 
%~ Such instabilities can occur if the mass density of the solid $\rhos$ is comparable to or less than the mass density of the fluid $\rhof$.
%~ The immersed aspect of our ILE approach uses a non-overlapping discretization of the fluid-structure interface for easy handling of large structural motions and deformations.
%~ These capabilities are made available through the hybrid immersed/partitioned nature of the algorithm.
% The surface mesh representation that interacts with the fluid moves according to the local fluid velocity and exerts 
% an approximate Lagrange multiplier force distribution back to the fluid. 
% The question remains as to how the Lagrange multiplier force distribution that enforces the interfacial constraints needs to be calculated.
% In our approach, the approximate Lagrange multiplier force is generated from penalty forces that link the surface mesh to the boundary of its volumetric solid counterpart.
%~ At least formally, in the limit of infinite spring stiffness, the two interface representations become exactly conformal in their motion.
Finally, the partitioned formulation also provides the flexibility of using multi-rate time-stepping, which 
can relax time-step restrictions and boost the computational performance by allowing separate choices of time-step sizes in the fluid and
 solid mechanics solvers. 
 We demonstrate this capability of the present approach for the simulations of the modified Turek-Hron benchmark and damped structural instability 
 of an enclosed fluid reservoir.

\section{Continuous equations of motion}
\label{sec:formulation}

This section presents the
partitioned and ILE forms of the equations describing fluid-structure interaction with an immersed elastic body. 
A conventional partitioned formulation for FSI is outlined in Sec.~\ref{sec:partitioned} followed by
the general ILE formulation in Sec.~\ref{sec:ILE_formulation}.
Secs.~\ref{sec:penalty_ILE_incomp} and~\ref{sec:penalty_ILE_comp} introduce penalty methods for the ILE formulation that are used in practical computations 
of incompressible and compressible solid structures, respectively.
Sec.~\ref{sec:structural_mechanics} details the weak formulation of the structural mechanics along with specific material constitutive models used in 
our numerical tests.

\subsection{The partitioned formulation}
\label{sec:partitioned}

\begin{figure}[t!!]
		\centering
			\includegraphics[width=0.5\textwidth]{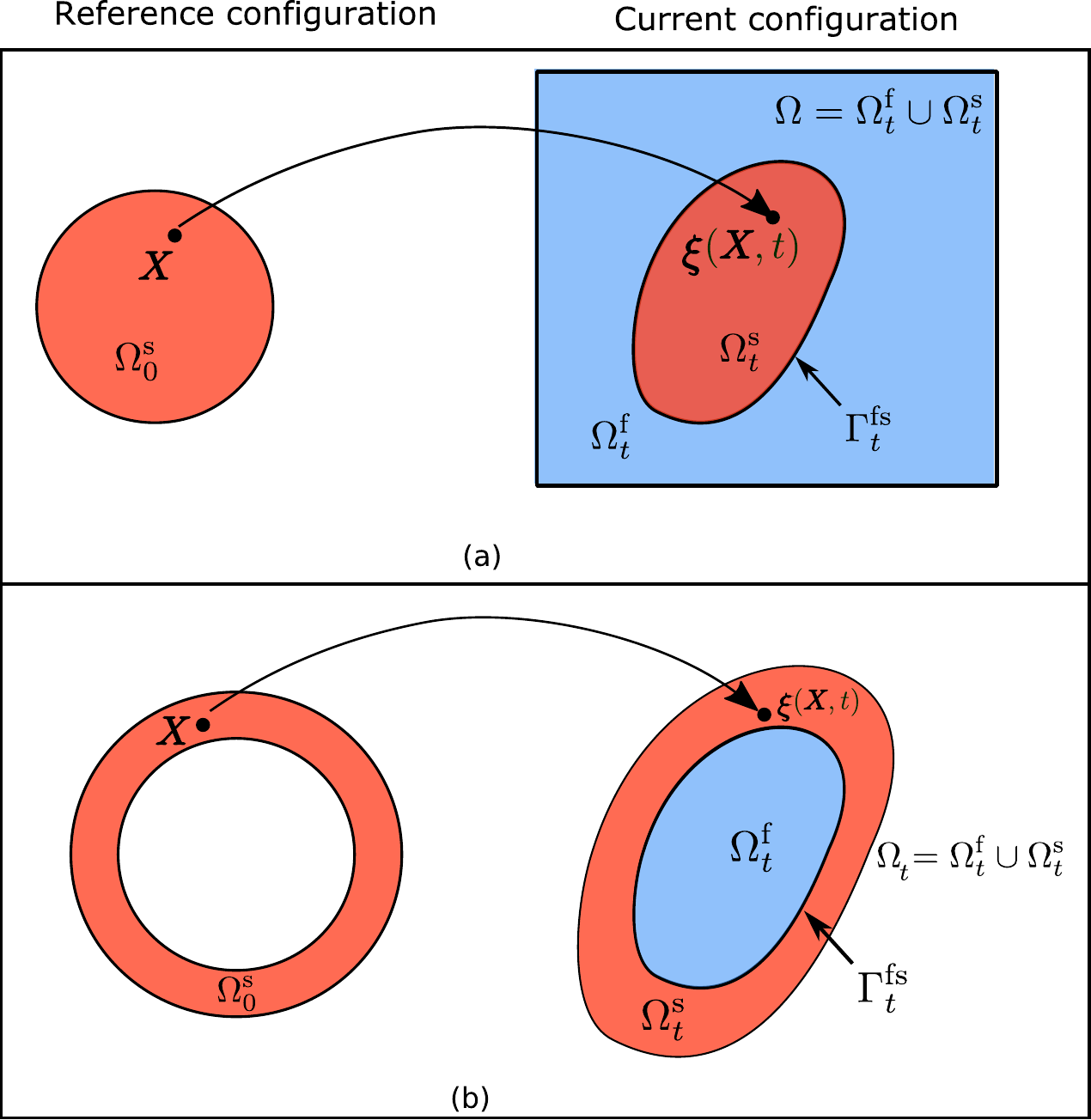}
		\caption{(a) Schematic of the domain $\Omega$ including time-dependent fluid and immersed solid subdomains, 
		$\Omega^{\text{f}}_t$ and $\Omega^{\text{s}}_t$, in a partitioned formulation.
	The solid is described using reference coordinates $\s \in \Omega^{\text{s}}_0$. Reference and current coordinates 
	are connected by the mapping $\Y : (\Omega^{\text{s}}_0 , t)\mapsto \Omega^{\text{s}}_t $.
	(b) Schematic of the partitioned formulation for a model in which the solid structure is
	 not fully immersed in the fluid.
	In the specific case shown here, the solid structure is a thick shell that encloses the fluid region, but its outer boundary
	does not come in contact with the fluid. Notice that in this case, the physical
	domain itself will generally be time dependent ($\Omega \equiv \Omega_t$).}
		\label{fig:Lag_Eul_schematic} 
\end{figure}

We consider a physical
domain $\Omega$ subdivided into time-dependent fluid and structure subdomains, $\Omega^{\text{f}}_t$ and $\Omega^{\text{s}}_t$,
indexed by time $t$, so that $\Omega=\overline{\Omega}^{\text{f}}_t \, \cup \, \overline{\Omega}^{\text{s}}_t$,
with $\x \in \Omega$ indicating fixed physical coordinates.
We first outline the partitioned formulation for an incompressible structure that is  
%~ defined by a simply connected domain, and
 fully immersed in an incompressible fluid 
within a fixed physical domain; see Fig.~\ref{fig:Lag_Eul_schematic}(a). 
%~ Note that in this case $\Omega \equiv \Omega_0 \equiv \Omega_t$. 
%~ but will demonstrate later in Sec.~\ref{sec:ILE_formulation} how our ILE formulation is able to
%~ handle the time-dependent $\Omega_t$ as in Fig.~\ref{fig:Lag_Eul_schematic}(b), thanks to the immersed aspect of our formulation.
 %~ In Sec.~\ref{subsubsec:pressure_wave} will demonstrate an example of a time-dependent computational domain in a fixed extended domain.
 We describe the fluid dynamics in Eulerian form and assume an incompressible Newtonian fluid with uniform density $\rhof$ and viscosity $\muf$. 
 The fluid stress tensor $\sigmaf$ 
is
\begin{equation}
 \sigmaf(\x,t) = -p(\x,t) \, \II + \muf \left( \grad\u(\x,t) + \grad\u^T(\x,t) \right),  \,  \x \in \Omega^{\text{f}}_t,
\end{equation}
in which $p(\x,t)$ is the fluid pressure and $\u(\x,t)$ is the fluid velocity.
The structural kinematics are described in Lagrangian form using reference coordinates $\s \in \Omega^{\text{s}}_0$ attached to the structure. We use the 
motion map $\Y : (\Omega^{\text{s}}_0 , t)\mapsto \Omega^{\text{s}}_t $ to determine the current position of material point $\s$ at time $t$.
The regions meet along the fluid-structure interface, $\Gamma^{\fs}_t = \overline{\Omega}^{\text{f}}_t \cap \overline{\Omega}^{\text{s}}_t$.
The equations of motion for the coupled fluid-structure system are
   \begin{align}
       \label{eq:iim_momentum} \rhof \, \frac{{\mathrm D} \u}{\mathrm {Dt}}(\x,t) &= \grad \cdot \sigmaf(\x,t), &  \x \in \Omega^{\text{f}}_t,\\
       \label{eq:iim_continuity}  \grad \cdot \u(\x,t) &= 0, &  \x \in \Omega^{\text{f}}_t,\\
       	  %~     \label{eq:traction_jump_codim1} \llbracket \sigmaf(\Y(\s,t),t) \cdot \vec{n}(\Y(\s,t),t) \rrbracket &= -\mathcal{I}^{-1}(\s,t) \, \F(\s,t),  &  \s \in \Gamma^{\fs}_0,\\
       	 \label{eq:solid_mom_immersed}  \rhos_0 \,\frac{\partial^2\Y}{\partial t^2}(\s,t) &= \grad_{\s} \cdot \PP(\s,t),  &  \s \in \Omegas_0,\\  %\int_{\Gamma^{\fs}_t} \tauf^{+}(\x,t) \, {\mathrm d} a,\\
     %~ \label{eq:solid_incomp} J(\s,t) &\approx 1,  &  \s \in \Omegas_0,\\ 
          \label{eq:solid_incomp} J(\s,t) &= 1,  &  \s \in \Omegas_0,\\ 
          		\label{eq:immersed_kinematic_cond} \frac{\partial\Y}{\partial t}(\s,t) &= \vec{u}(\Y(\s,t),t),  &  \s \in \Gamma^{\fs}_0,\\
     	\label{eq:force_balance}   j^{-1}(\s,t) \, {\PP(\s,t) \,  \vec{N}(\s)} &=  \bbsigma^{\text{f}}(\Y(\s,t),t) \, \vec{n}(\Y(\s,t),t), &  \s \in \Gamma^{\fs}_0,
 %~ J^{-1}(\s,t) \, \PP(\s,t) \, \FF^{T}(\s,t) \, \vec{n}(\Y(\s,t),t)
 %~ \label{eq:solid_ang_mom_immersed}  \frac{\mathrm d}{\mathrm{dt}}\int_{\Omega^{\text{s}}_0} \bigl(\s - \s_{\com}\bigl)
  %~ \times \bigl(\rhos \, \frac{\partial\Y}{\partial t}(\s,t)\bigl)\, {\mathrm d} \s &= \int_{\Gamma^{\fs}_t} \bigl(\x - \Y_{\com}(t)\bigl) \times  \tauf^{+}(\x,t)  \, {\mathrm d} a,
    \end{align} 
in which $\rhos_0$ is the structural mass density in the reference configuration 
(constant in time for an incompressible structure) and 
$J=\text{det}(\FF)$ is the (volumetric) Jacobian determinant associated with the structural deformation, 
in which $\FF(\s,t)={\partial \Y}/{\partial \s}$ 
is the deformation gradient tensor. In Eq.~(\ref{eq:force_balance}), 
$\vec{n}(\Y(\s,t),t)$ is the 
outward unit normal vector pointing into $\Omega^{\text{f}}_t$ along $\Gamma^{\fs}_t$, 
$\vec{N}(\s)$ is the outward normal vector in the reference configuration along $\Gamma^{\fs}_0$, and 
 $j(\s,t)$ is the surface Jacobian determinant
 that converts the surface force density from force per 
unit area in the current configuration to force per unit area in the reference configuration.
 $j(\s,t)$ is related to the volumetric Jacobian determinant  $J(\s,t)$ through Nanson's formula, 
 $j(\s,t)/J(\s,t)=\|\FF^{-T}(\s,t)\N(\s)\|$.  
$\PP(\s,t)$ is the first Piola-Kirchhoff stress tensor, which is defined for an incompressible a hyperelastic material by
\begin{equation}
 \PP(\s,t) = \frac{\partial \psi}{\partial \FF}(\s,t) - \Phi(\s,t) \, \FF^{-T}(\s,t),  \,  \s  \in \Omegas_0.\\
 \label{eq:PP_def}
\end{equation}
Here, $\psi(\s,t)$ is the strain-energy functional and $\Phi(\s,t)$
  is the hydrostatic pressure that enforces the incompressibility of the structure.
  %~ $\n(\x,t)$ is the unit normal vector pointing into $\Omega^{\text{f}}_t$ along $\Gamma^{\fs}_t$.
 %~ and $\grad_{\s} \cdot \mbox{}$  is the divergence operator in the reference coordinate system.
%~ $\F(\s,t)$ is an interfacial surface force density that 
%~ is the Lagrange multiplier for the kinematic condition, Eq.~(\ref{eq:immersed_kinematic_cond}), along the fluid-structure interface $\Gamma^{\fs}_t$.
%~ Although not considered further here, it is straightforward to relax the incompressibility constraint so that the structure is allowed to undergo local volume changes 
%~ so long as the total volume of the structure is constant.
In the aforementioned governing equations, 
Eqs.~(\ref{eq:iim_momentum}) and (\ref{eq:iim_continuity}) describe the momentum and incompressibility of the fluid, 
whereas Eqs.~(\ref{eq:solid_mom_immersed}) and (\ref{eq:solid_incomp}) describe the momentum and incompressibility of the structure.
Eqs.~(\ref{eq:immersed_kinematic_cond}) and (\ref{eq:force_balance}) respectively account for the kinematic and dynamic conditions
 at the fluid-structure interface. Eq.~(\ref{eq:immersed_kinematic_cond}) describes the no-slip and no-penetration condition
 at the FSI interface, whereas Eq.~(\ref{eq:force_balance}) describes the 
 force balance along the fluid-solid interface.

It is also possible to consider a partitioned fluid-structure formulation for a system in which the solid domain 
is only partially in contact with the fluid, as depicted for the 
non-simply connected ring-shaped solid domain in Fig.~\ref{fig:Lag_Eul_schematic}(b).
 Note that in this case, the physical domain $\Omega \equiv \Omega_t \equiv \overline{\Omega}^{\text{f}}_t \, \cup \, \overline{\Omega}^{\text{s}}_t$ may be time dependent.
In the following section, we will demonstrate how in this case, our ILE coupling 
allows us to use a fully compressible structural model.
 For compressible structures, the constraint in Eq.~(\ref{eq:solid_incomp}) is no longer enforced, and the density 
 of the solid in the current configuration is $\rhos = \rhos_0/J$.
 Although not considered here, it is also straightforward to relax the incompressibility constraint in the fully immersed
  case and allow the structure to experience local volume changes 
 so long as the total volume of the structure is constant.

\subsection{The ILE formulation}
\label{sec:ILE_formulation}

\begin{figure}[t!!!]
	\centering
	\includegraphics[width=0.85\textwidth]{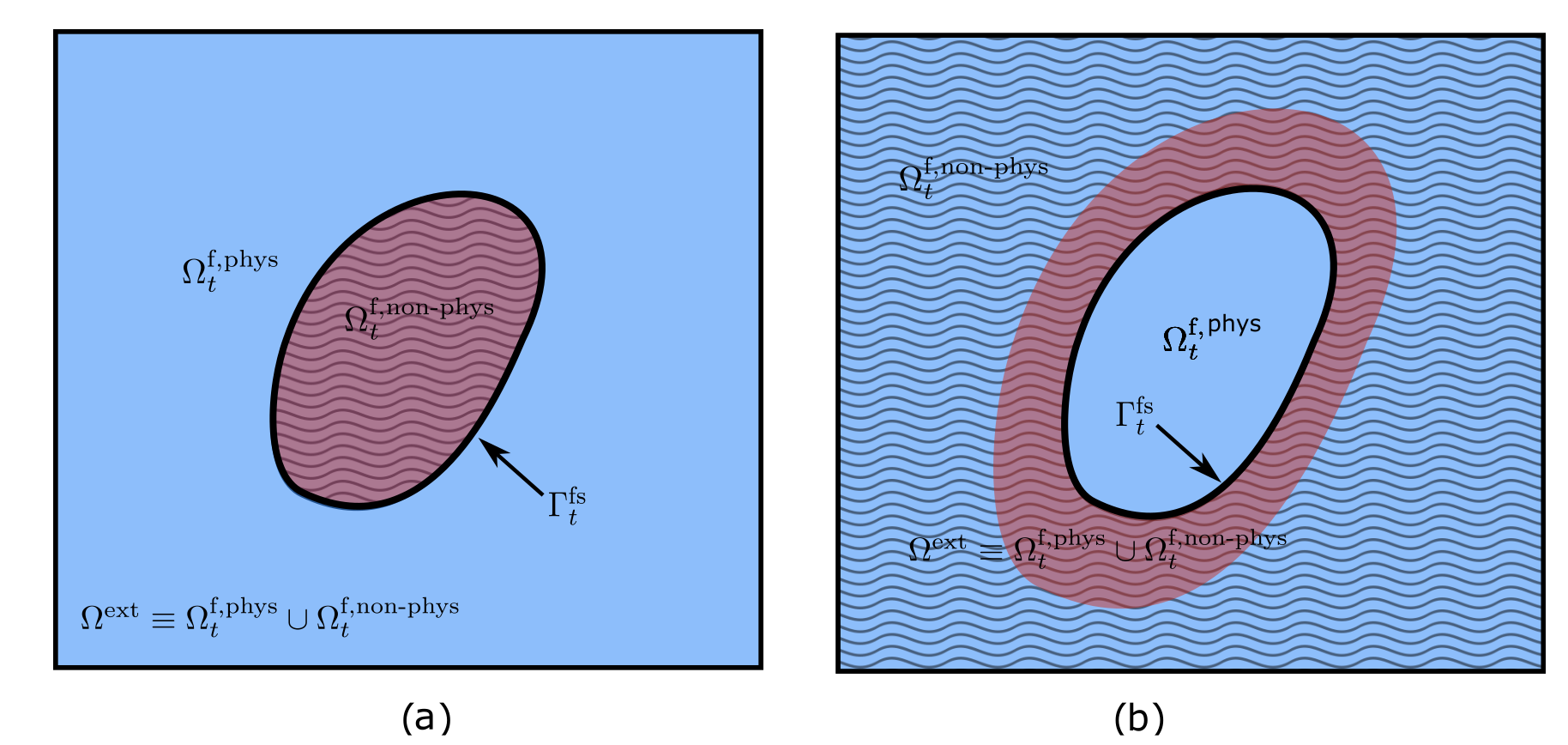}
	\caption{In the immersed Lagrangian-Eulerian method,
	the incompressible Navier-Stokes 
equations are solved on an extended computational domain $\Omega^{\text{ext}}$ that 
incorporates both fluid and solid subregions, and
is split
into a physical $\Omega^{\text{f},\text{phys}}_t$ and a non-physical $\Omega^{\text{f},\text{non-phys}}_t$ fluid regions. 
(a) In this configuration a simply connected domain represents the volumetric structure that is
 immersed in the fluid. The non-physical fluid region in this case is the region
  occupied by the solid and shown by wavy stripes. (b) In this case, the physical fluid region is 
   enclosed by the thick shell-like structure of the solid. The non-physical domain 
  (shown by wavy stripes) includes the solid domain in addition to an extended region to obtain a fixed  
  Cartesian computational domain.}
\label{fig:IIM_schematic} 
\end{figure} 
%~ \subsubsection{The ILE formulation for a nearly-incompressible solid structure}
The fundamental difference between the partitioned formulation detailed in Sec.~\ref{sec:partitioned} and the ILE 
formulation is that the ILE formulation solves the incompressible Navier-Stokes 
equations on an extended computational domain $\Omega^{\text{ext}}$ that
incorporates both fluid and solid subregions. To be more specific,
we split the computational domain $\Omega^{\text{ext}}$  into a \emph{physical} 
fluid region $\Omega^{\text{f},\text{phys}}_t$ and a \emph{non-physical} fluid region $\Omega^{\text{f},\text{non-phys}}_t$,
each parameterized by time $t$, with superscripts ``$\text{phys}$'' (``$\text{non-phys}$'') indicating values obtained from 
the physical (non-physical) side of the fluid at the fluid-structure interface; 
see Fig.~\ref{fig:IIM_schematic}.
Using this notation, we have $\Gamma^{\fs}_t = {\overline{\Omega}^{\text{f},\text{phys}}_t}  \cap  \overline{\Omega}^{\text{f},\text{non-phys}}_t$.
 In the case of fully immersed solid structure, as in Fig.~\ref{fig:IIM_schematic}(a), we have $ \Omegas_t  \equiv \Omega^{\text{f},\text{non-phys}}_t$.
For the case shown in Fig.~\ref{fig:IIM_schematic}(b), the structure is only in contact with the physical fluid  region on the inner side of its boundary, 
corresponding to the same physical problem as in Fig.~\ref{fig:Lag_Eul_schematic}(b).
In this case, $\Omegas_t \subset {\Omega^{\text{f},\text{non-phys}}_t}$, but the non-physical domain is further extended beyond the solid domain,  
allowing us to use a fixed computational domain. 
Note that in both Figs.~\ref{fig:IIM_schematic}(a) and (b), the extended computational domain is fixed, i.e. $\Omega^{\text{ext}} \equiv \Omega^{\text{ext}}_0 \equiv \Omega^{\text{ext}}_t$. 
We now extend the fluid velocity $\u$, pressure $p$, viscosity $\muf$, and stress $\sigmaf$ to be defined in the entire computational domain $\Omega^{\text{ext}}$, 
and define the \emph{extended fluid stress tensor} $\sigmaf$ as
\begin{equation}
 \sigmaf(\x,t) = -p(\x,t) \, \II + \muf \left( \grad\u(\x,t) + \grad\u^T(\x,t) \right), \, \, \, \, \, \, \, \,\, \x \in \Omega^{\text{ext}} \equiv {\overline{\Omega}^{\text{f},\text{phys}}_t}  \cup  \overline{\Omega}^{\text{f},\text{non-phys}}_t.
\end{equation}
The kinematic constraint, which requires that the fluid and solid move together along the fluid-structure 
interface, is imposed through an interfacial force density applied 
along $\Gamma^\text{fs}_t$. This singular force density implies 
a discontinuity or jump in the traction associated with the extended fluid stress along $\Gamma^{\fs}_t$.
We denote a jump in a scalar field $\psi(\x,t)$ at position $\x = \Y(\s,t)$ along the interface as
$\llbracket \psi(\vec{x},t) \rrbracket = \psi^{\text{phys}}(\vec{x},t) - \psi^{\text{non-phys}}(\vec{x},t),$
in which $\psi^{\text{phys}}(\vec{x},t)$ and $\psi^{\text{non-phys}}(\vec{x},t)$ are respectively the limiting values approaching the interface position $\vec{x}$ from the 
physical 
fluid region $\Omega^{\text{f},\text{phys}}_t$ and the non-physical fluid region $\Omega^{\text{f},\text{non-phys}}_t$ in the normal direction.
%~ We account 
%~ for the physical jump conditions 
%~ in the extended fluid stress
%~ along the fluid-solid interface
 %~ by leveraging a coupling scheme based on the
%~ immersed interface method for discrete surfaces \cite{kolahdouz2020immersed}.
The resulting ILE formulation is
    \begin{align}
       \label{eq:ILE_momentum} \rhof \, \frac{{\mathrm D} \u}{\mathrm {Dt}}(\x,t) &= \grad \cdot \bbsigma^{\text{f}}(\x,t), &  \x \in \Omega^{\text{ext}},\\
       \label{eq:ILE_continuity}  \grad \cdot \u(\x,t) &= 0, &  \x \in \Omega^{\text{ext}},\\
             \label{eq:ILE_traction_jump_codim1} \llbracket \sigmaf(\Y(\s,t),t) \cdot \vec{n}(\Y(\s,t),t) \rrbracket &= -j^{-1}(\s,t) \, \F(\s,t),  &  \s \in \Gamma^{\fs}_0,\\ %\bigg( \, \kappa\left(\Y-\X\right) + \eta \, \big(\frac{\partial \Y}{\partial t} - \frac{\partial \X}{\partial t}\big)\bigg),   &  \s \in \Gamma^{\fs}_0,\\
	       	 \label{eq:ILE_solid_mom_immersed}  \rhos_0 \,\frac{\partial^2 \Y}{\partial t^2}(\s,t) &= \grad_{\s} \cdot \PP(\s,t) ,  &  \s \in \Omegas_0,\\  %\int_{\Gamma^{\fs}_t} \tauf^{+}(\x,t) \, {\mathrm d} a,\\
		          \label{eq:ILE_solid_incomp} J(\s,t) &= 1,  &  \s \in \Omegas_0,\\ 
		\label{eq:ILE_immersed_kinematic_cond} \frac{\partial\Y}{\partial t}(\s,t) &= \vec{u}(\Y(\s,t),t),  &  \s \in \Gamma^{\fs}_0,\\
          %~ \label{eq:ILE_solid_incomp} J(\s,t) & \approx 1,  &  \s \in \Omegas_0,\\ 
     	\label{eq:ILE_force_balance} j^{-1}(\s,t) \, {\PP(\s,t) \,  \vec{N}(\s)} &=  \bbsigma^{\text{f},\text{phys}}(\Y(\s,t),t) \, \vec{n}(\Y(\s,t),t), &  \s \in \Gamma^{\fs}_0,
    \end{align}
 in which 
$\F(\s,t)$ is an interfacial surface force density that 
is the Lagrange multiplier for the kinematic condition in Eq.~(\ref{eq:ILE_immersed_kinematic_cond}).
 Eq.~(\ref{eq:ILE_traction_jump_codim1}) describes discontinuities in the extended 
fluid pressure and viscous stress that are related to the normal and tangential components of the interfacial Lagrange multiplier force, respectively.
Additionally, the right hand side of Eq.~(\ref{eq:ILE_force_balance}) is the \emph{exterior} 
fluid traction in current coordinates
\begin{align}
\vec{\tau}^{\text{f},\text{phys}}(\x,t) = \bbsigma^{\text{f},\text{phys}}(\x,t) \, \vec{n}(\x,t), \, \vec{x} \in \Gamma^{\fs}_t.
\end{align}
This equation accounts for the dynamic interface condition, which is
 the ``Neumann part'' of our Dirichlet-Neumann coupling approach, in which
the exterior fluid traction forces are imposed as boundary conditions to the solid domain.
Note that this equation also implies that \emph{only} the \emph{physical} fluid stresses from within $\Omega^{\text{f},\text{phys}}_t$ are able to 
drive the dynamics of the structure.
%~ Eq.~(\ref{eq:ILE_traction_jump_codim1}) implies discontinuities in the extended 
%~ fluid pressure and viscous stress that are related, respectively, to the normal and tangential components of the interfacial Lagrange multiplier force.
%~ Further, the jump discontinuity in Eq.~(\ref{eq:traction_jump_codim1}) can be decomposed into discontinuities in the pressure and viscous stress, which in current coordinates are expressed as
   %~ \begin{align}
		%~ \label{eq:pj_codim1} \llbracket p(\vec{x},t) \rrbracket &=  J^{-1}(\Y^{-1}(\x,t),t) \, \F(\Y^{-1}(\x,t),t) \cdot \vec{n}(\vec{x},t),  &  \x \in \Gamma^{\fs}_t,\\
		%~ \label{eq:du_dn_jump_codim1}     \left\llbracket \muf\frac{\partial\vec{u}}{\partial \vec{n}}(\x,t) \right\rrbracket  &= - (\II - \vec{n}(\x,t)\otimes\vec{n}(\x,t)) \, J^{-1}(\Y^{-1}(\x,t),t) \, \F(\Y^{-1}(\x,t),t),  &  \x \in \Gamma^{\fs}_t.
    %~ \end{align}
\begin{figure}[b!!!]
	\centering
	\includegraphics[width=\textwidth]{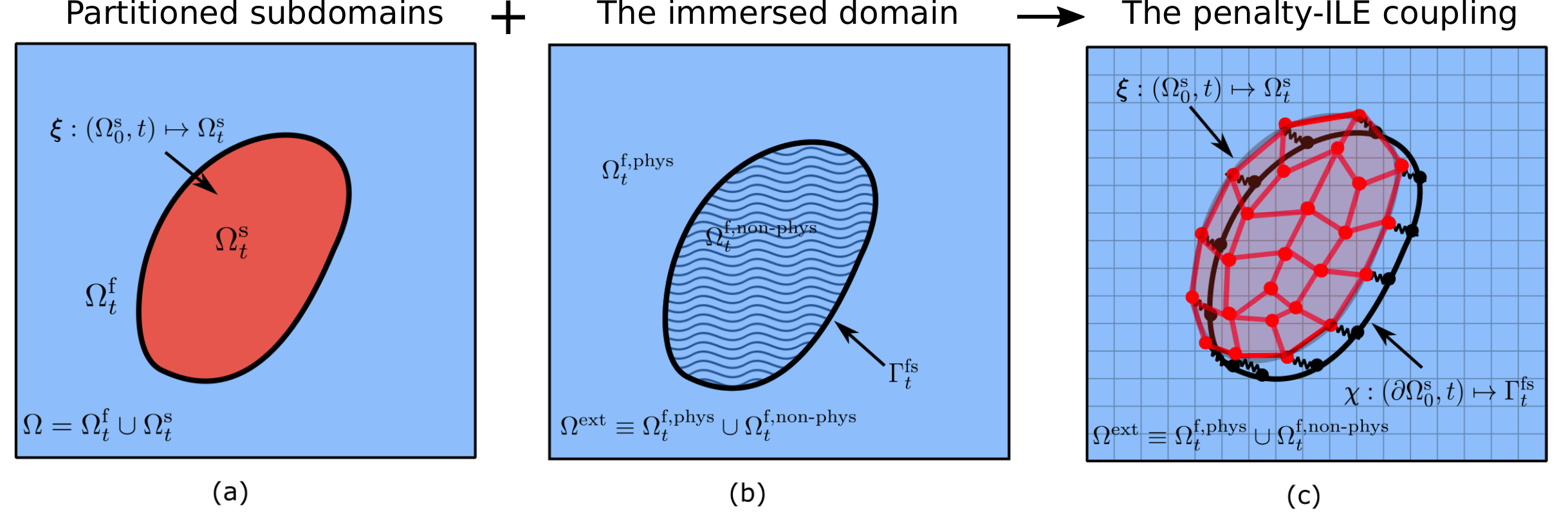}
	\caption{
Aspects of (a) partitioned and (b) immersed FSI formulations are integrated in the  (c) penalty-ILE coupling.
The fluid-structure interface representation (determined by $\vec{\chi}$) conforms to
 the boundary of the structure (determined by $\vec{\xi}$) in an approximate sense.
 The discrepancy between the two representations is exaggerated here for illustration purposes. 
In practice, we ensure that this maximum displacement is no greater than $0.1$ of the Cartesian grid spacing. 
%~ Otherwise, the code will abort and the penalty parameters have to be re-determined.
The local fluid velocity determines the motion of the surface representation $\Gamma^{\fs}_t$,
whereas the no-slip condition is satisfied in an approximate sense by spring-like forces that penalize displacements between the two representations 
of the fluid-structure interface. }
\label{fig:ILE_FSI_schematic} 
\end{figure}

%%%%%%%%%%%%%%%%%%%%%%%%%%%%%%%%%%%%%%%%%%%%%%%%%%%%%%%%%%%%%%%%%%%%%%%%%%%%%%%%%%%%%%%%%%%%%%%%%%%%%%%%%%%%%%%%%%%%%%%%%%%
\subsubsection{The penalty ILE formulation}
\label{sec:penalty_ILE_incomp}

We first describe our penalty ILE formulation for a solid structure fully immersed in a fluid domain.
Exactly imposing the kinematic constraint, Eq.~(\ref{eq:immersed_kinematic_cond}), 
in the formulation detailed in Sec.~\ref{sec:ILE_formulation} requires solving a saddle-point system 
that couples the Eulerian and Lagrangian variables.
As in our earlier work on rigid-body FSI \cite{kolahdouz2021sharp}, we develop a numerical method that avoids the need to solve such systems by relaxing
 the kinematic constraint. We do so by introducing two representations of 
the fluid-structure interface 
and applying forces that act to penalize displacements between the two representations. 
This penalty method thereby determines an approximate Lagrange multiplier force instead of solving for the exact Lagrange multiplier to impose the kinematic matching condition. 
Specifically, along with the mapping $\Y(\s,t)$ that 
determines the kinematics of the
structure, we introduce an explicit representation of the fluid-structure interface that is parameterized by the deformation mapping $\X(\s,t)$ and that moves 
with the fluid, so that $\partial\X(\s,t)/\partial t = \vec{u}(\X(\s,t),t)$; see Fig.~\ref{fig:ILE_FSI_schematic}(c).
%~ This surface representation corresponds to the restriction of the volumetric solid to its surface boundaries.

The approximate Lagrange multiplier force is determined by linear spring and dampers,
\begin{align}
 \F(\s,t) = \kappa\left(\Y(\s,t)-\X(\s,t)\right) + \eta \, \left(\frac{\partial \Y}{\partial t}(\s,t) - \frac{\partial\X(\s,t)}{\partial t}\right), \, \s \in \Gamma^{\fs}_0,
 \label{eq:penalty_force}
\end{align}
in which $\kappa$ is a tethering penalty parameter
 and $\eta$ is a damping penalty parameter. 
 This force penalizes deviations from
the constraint $ \Y(\s,t) =  \X(\s,t)$ and, in the discretized equations, 
acts to ensure that the two interface representations are at least approximately conformal in their motion.
It is possible to control the discrepancy between the two configurations, and 
 as $\kappa \rightarrow \infty$, the formulation exactly imposes the constraint that the two interface representations move together. 
 Although kinematic constraints can be imposed through spring forces alone, 
 we have found 
 that including a small amount of damping can reduce spurious numerical oscillations. We demonstrate in Sec.~\ref{subsubsec:soft_disc} that including damping can 
 reduce the value of $\kappa$ required to achieve a given displacement 
 tolerance, which, in turn, can allow for larger time step sizes.
 
 %~ Although in principle the constraint formulation can be imposed through spring forces alone, we have found out that
 %~ a slight amount of additional damping force 
 %~ could broaden the stable choices for $\kappa$ and help with a faster tuning process.
 %~ helps reducing spurious numerical oscillations.

Additionally, note that in the penalty ILE formulation the incompressibility of the fluid is discretely imposed through the solution of the Navier-Stokes equations.
 Therefore, for consistency,
 a structure that is fully immersed in fluid must also be incompressible in the present formulation.
 %~ and obtaining a sensible solution the incompressibility condition needs 
 %~ to be also satisfied within the solid structure surrounded by the fluid region. 
Exactly imposing this constraint in a Lagrangian framework is challenging.
%~ requires solving for a hydrodynamic pressure in another 
%~ saddle-point system as described in the previous section.
%~ Such solution needs to overcome
%~ the so-called \emph{inf-sup} condition restriction between the Lagrangian velocity and pressure.
To simplify the implementation, we relax this condition and use a nearly incompressible structural model that
is characterized by
a Poisson's ratio $\nu$ approaching 0.5 and a corresponding bulk modulus $\Ks$ approaching infinity. This implies 
that $J(\s,t) \rightarrow 1$. Our numerical tests demonstrate that the resulting volume conservation is superior to an immersed finite element-finite 
difference scheme that 
treats the structure as exactly incompressible \cite{BEGriffith17-ibfe}.

\subsubsection{The penalty ILE formulation for fully compressible solid structures}
\label{sec:penalty_ILE_comp}

\begin{figure}[b!!!]
	\centering
	\includegraphics[width=\textwidth]{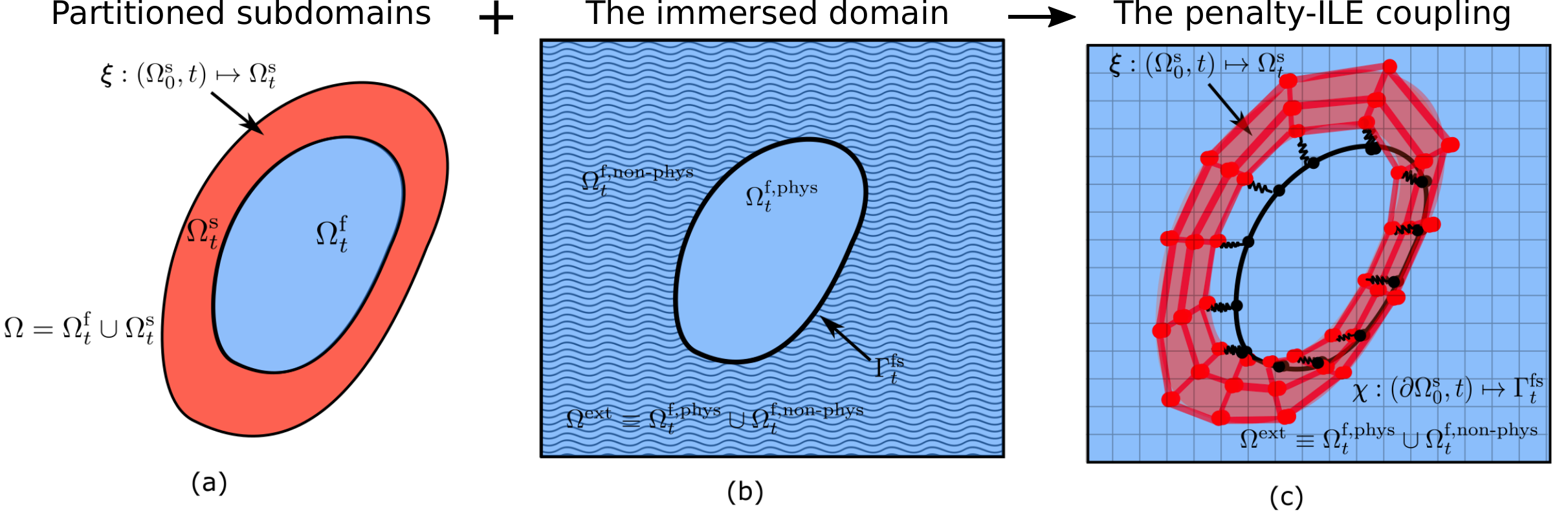}
	\caption{
	The coupling for the compressible structure combines 
 (a) the partitioned formulation and (b) immersed FSI formulations to create (c) the penalty-ILE coupling.
The fluid-structure interface representation in this case (determined by $\vec{\chi}$) conforms to
 the inner boundary of the structure (determined by $\vec{\xi}$) in an approximate sense.
 The displacement between the two representations is exaggerated here for illustration purposes.}
\label{fig:ILE_FSI_schematic-comp} 
\end{figure} 
It is straightforward to use
our penalty ILE formulation to 
 model fully compressible solid structures so long as
 the structure is not fully immersed in fluid (i.e., part of the solid boundary is not in contact with the fluid).
%~ Although not considered herein, our formulation also can readily handle fully immersed compressible structures that have a fixed total volume.
We enforce 
the kinematic constraint through the same penalty strategy as in the incompressible model, except in this case,
the additional representation of 
the fluid-structure interface is only a subset of the total solid boundary; see Fig.~\ref{fig:ILE_FSI_schematic-comp}.
Because part of the solid boundary 
does not contact the physical fluid subregion, the incompressibility of the extended fluid region does not constrain the volume of the structure.
%~ is allowed to freely move with no attachment or interaction with the physical fluid region, 
%~ compressible structures can be considered despite enforcing the fluid incompressibility in the extended computational domain.
An example of such a model is flow through an enclosed tube or a reservoir for which
  the FSI coupling is enabled only at the inner surface of the tube while the outer side is handled through 
  suitable boundary conditions for the structural model.
  %~ a traction-free condition in which
  %~ the volumetric structure freely displaces according to the deformations from the solid mechanics solution.
%~ in the formulation presented in Sec.~\ref{sec:partitioned} is related to

%~ Considering the approximations for imposing the kinematic constraint and the incompressibility of the solid structure, 

\subsection{The structural mechanics formulation}
\label{sec:structural_mechanics}

%~ In this section we introduce a weak formulation of Lagrangian equations of motion for the solid structure. 
We use a standard Galerkin finite element approach to model the structural deformation. 
Such methods rely on a weak form of the equations of elastodynamics. Briefly, by assuming that the immersed
solid structure has uniform density and integrating Eq.~(\ref{eq:solid_mom_immersed}) for all smooth test
 functions $\vec{\uppsi}(\s) \in V$ and 
$V \subseteq (H^1(\Omega^{\text{s}}_0))^3$, we obtain
\begin{equation}
 \rhos_0 \int_{\Omega^{\text{s}}_0} \left(\frac{\partial^2 \Y}{\partial t^2}(\s,t)\right) \, \cdot \, \vec{\uppsi}(\s) \, {\mathrm d} \s = \int_{\Omega^{\text{s}}_0} \, \left(\grad_{\s} \cdot \PP(\s,t)\right) \, \cdot \, \vec{\uppsi}(\s) \,{\mathrm d} \s.
\label{eq:weak_form1}
\end{equation}
(We remark that it is straightforward to consider nonuniform density structural models.) 
Integrating by parts yields
\begin{equation}
\rhos_0 \int_{\Omega^{\text{s}}_0} \left(\frac{\partial^2 \Y}{\partial t^2}(\s,t)\right) \, \cdot \, \vec{\uppsi}(\s) \, {\mathrm d} \s = 
-\int_{\Omega^{\text{s}}_0} \, \PP(\s,t) \, \colon \,  \grad_{\s} \, \vec{\uppsi}(\s) \,{\mathrm d} \s
+ \int_{\Gamma^{\fs}_0} \, \left(\PP(\s,t) \, \N(\s) \right) \, \cdot \, \vec{\uppsi}(\s) \,\DA.
\label{eq:weak_form2}
\end{equation}
%~ in which $\N(\s)$ is the outward unit normal on ${\partial \Omega^{\text{s}}_0}$.
Note that the second term on the right-hand-side of Eq.~(\ref{eq:weak_form2}) also appears in the force balance equation,  Eq.~(\ref{eq:force_balance}),
 and can be replaced by exterior fluid traction forces as
 %~ Using Nanson's formula, which relates surface normals and area elements reference and current coordinates, Eq.~(\ref{eq:weak_form2}) can be rewritten as
\begin{equation} 
\begin{split}
\rhos_0 \int_{\Omega^{\text{s}}_0} \left(\frac{\partial^2 \Y}{\partial t^2}(\s,t)\right) \, \cdot \, \vec{\uppsi}(\s) \, {\mathrm d} \s &= 
 -\int_{\Omega^{\text{s}}_0} \, \PP(\s,t) \, \colon \,  \grad_{\s} \, \vec{\uppsi}(\s) \,{\mathrm d} \s  \\
 & + \int_{\Gamma^{\fs}_0} \,  j(\s,t) \, \vec{\tau}^{\text{f},\text{phys}}(\X(\s,t),t) \, \cdot \, \vec{\uppsi}(\s) \,\DA.
 \end{split}
\label{eq:weak_form3}
\end{equation} 
Notice that in the surface integral in Eq.~(\ref{eq:weak_form3}), we use $j(\s,t)$ to obtain the fluid traction forces per unit \emph{reference} area.
We directly discretize this equation in our numerical scheme. 
%~ Note that for modeling compressible materials in a time-dependent computational domain, described in Fig.~\ref{fig:Lag_Eul_schematic}(b),
%~ part of the solid structure boundary may not have an interface with the fluid. In such scenario Eq.~(\ref{eq:weak_form3}) needs to be re-written as,
%~ \begin{multline}
%~ \rhos \int_{\Omega^{\text{s}}_0} \left(\frac{\partial^2 \Y}{\partial t^2}(\s,t)\right) \, \cdot \, \vec{\uppsi}(\s) \, {\mathrm d} \s = 
%~ \int_{\Gamma^{\fs}} \, \vec{\tau}^{\text{f},+}(\Y(\s,t),t) \, \cdot \, \vec{\uppsi}(\s) \,\Da \\
%~ + \int_{\partial \Omega^{\text{s}}_0\backslash\Gamma^{\fs}} \, \left(\PP(\s,t) \, \N(\s) \right) \, \cdot \, \vec{\uppsi}(\s) \,\DA
%~ -\int_{\Omega^{\text{s}}_0} \, \PP(\s,t) \, \colon \,  \grad_{\s} \, \vec{\uppsi}(\s) \,{\mathrm d} \s.
%~ \label{eq:weak_form4}
%~ \end{multline}

In this paper, we use hyperelastic structural models, for which the first Piola-Kirchoff 
stress tensor $\PP$
is related to the strain energy functional as in Eq.~(\ref{eq:PP_def}), and 
we use a nearly incompressible structural formulation as a penalty method for the exactly incompressible case. 
For nearly incompressible large-deformation elasticity,
we relax the incompressibility constraint by re-expressing $\PP$ in Eq.~(\ref{eq:PP_def}) via $\PP={\partial \psi}/{\partial \FF}$
and assuming an additive splitting of the elastic energies associated with
volume changing and volume preserving motions.
This splitting is motivated by the assumption that a uniform 
pressure should produce changes in size but not in shape \cite{sansour2008physical,vadala2020stabilization}. 
We first write $\psi$ for an isotropic material as a function of the scalar invariants $\psi=\psi(I_1,I_2,I_3)$ in which
\begin{align}
I_1 &= \text{tr}(\CC) = \sum_{i=1}^{3} \CC_{ii},\\
I_2 &=\frac{1}{2}\left(\text{tr}(\CC)^2 - \text{tr}(\CC^2) \right),\\
I_3 &= \text{det}(\CC),
\end{align}
with $\CC=\FF^T\FF$ being the right Cauchy-Green deformation tensor.
To achieve the desired split satisfying 
\emph{material frame invariance}, we use a standard approach based on the Flory decomposition \cite{flory1961thermodynamic} 
to reformulate the strain energy with the modified strain invariants, through
multiplicatively decomposing $\FF$ into dilatational and deviatoric parts, $\FF = J^{1/3} \bar{\FF}$.
This also results in a 
modified right Cauchy-Green strain tensor, $\bar{\CC}=J^{-2/3}\CC$.
Similarly, the first modified invariant of $\bar{\CC}$ is $\bar{I}_1=J^{-2/3}I_1$.
%--------------
%~ The strain energy functional $\psi$ can be decomposed and decoupled by isochoric $\mathcal{E}^{\text{iso}}$ and volumetric $\mathcal{E}^{\text{vol}}$ contributions,
%~ This allows for decoupling the energy  split associated with volume changing and volume preserving motions. 
The specific model used in this work is of the form
%~ The strain energy functional is then formulated as \cite{vadala2020stabilization},
\begin{equation}
\psi  = W(\bar{I}_1) + U(J), %= \frac{\Gs}{2}(\bar{I_1}-3) + \Ks\,(J\,\log J - J + 1)
\label{eq:energy_split}
\end{equation}
in which $W(\bar{I}_1) \equiv  \bar{W} \equiv W(\bar{\FF})$ characterizes the response of the material to shearing deformation
 in terms of the first modified invariant, and $U(J)$ is a volumetric energy that 
penalizes compression or expansion. It is preferred for $U$ to satisfy $U(1)=U^{\prime}(1)=0$.
Unless otherwise noted, we use a neo-Hookean material model 
for the elastic structure,
\begin{equation}
\psi (\FF) =  \frac{\Gs}{2}(\bar{I_1}-3) + \Ks\,(J\,\log J - J + 1).
\label{eq:neo_hookean}
\end{equation}
%~ It can also be beneficial to formulate the shearing energy $W$ 
%~ by multiplicatively decomposing $\FF$ into dilatational and deviatoric parts, $\FF = J^{1/3} \bar{\FF}$.
 %~ If we redefine the shearing energy to be dependent on the first modified invariant, $\bar{W}=W(\bar{\FF}) \equiv W(\bar{I}_1)$
Taking the derivative with respect to the deformation gradient tensor $\FF$ yields
%~ Finally, the first Piola-Kirchoff stress is calculated as,
\begin{equation}
\PP = \frac{\partial \psi}{\partial \FF} = \PP_{\text{dev}} + \PP_{\text{dil}},
\end{equation}
in which $\PP_{\text{dev}}=\frac{\partial W}{\partial \bar{I}_1} \, \frac{\partial \bar{I}_1}{\partial \FF}$ and 
$\PP_{\text{dil}}=\frac{\partial U}{\partial J} \, \frac{\partial J}{\partial \FF}$ are the
deviatoric and dilatational stresses, respectively. Using Eq.~(\ref{eq:neo_hookean}), these stresses are,
%~ is the deviatoric projection of the first Piola-Kirchoff stress given as,
\begin{align}
\PP_{\text{dev}} &= \Gs \, J^{-2/3}\, (\FF - \frac{I_1}{3}\FF^{-T}),\\
\PP_{\text{dil}} &=  \Ks \, J \, \log(J) \, \FF^{-T},
\end{align}
%~ and $\PP^{\text{dil}}$ is the dilatational stress given as,
%~ \begin{equation}
%~ \PP^{\text{dil}} =  J \, {\mathcal{E}^{\prime}}^{\text{vol}}(J) \, \FF^{-T} = \Ks \, J \, \log(J) \, \FF^{-T}
%~ \end{equation}
in which $\Gs$ and $\Ks$ are the shear and bulk moduli, which are related by the 
Poisson's ratio via $\Ks=\frac{2\Gs(1+\nu)}{3(1-2\nu)}$.
The above stress formulation can be easily extended to other isotropic and anisotropic constitutive laws.
In fact, the present ILE formulation is not restricted to neo-Hookean material models, and alternative 
hyperelastic or even inelastic material models can be readily used.

%~ is the deviatoric projection of the first Piola-Kirchoff stress given as,
%~ \begin{equation}
%~ \PP^{\text{dev}} = \text{dev}\bigl(\frac{\partial \mathcal{E}}{\partial \FF}\bigl) = \Gs \, J^{-2/3}\, (\FF - \frac{I_1}{3}\FF^{-T})
%~ \end{equation}
%~ and $\PP^{\text{dil}}$ is the dilatational stress given as,
%~ \begin{equation}
%~ \PP^{\text{dil}} =  J \, {\mathcal{E}^{\prime}}^{\text{vol}}(J) \, \FF^{-T} = \Ks \, J \, \log(J) \, \FF^{-T}
%~ \end{equation}

\section{Numerical discretization}
\label{sec:numerical}

We use a standard $H^1$-conforming Lagrangian finite element method to discretize the elastodynamics equations and
a finite difference method to discretize the Navier-Stokes equations.
 Further details of the discretization for each sub-problem are given below.

\subsection{Structural discretization}
\label{subsec:methodology-FE}

%~ The numerical approximation to the Lagrangian domain largely follows the approach in our previous work \cite{kolahdouz2021sharp}.
Let $\cT_h$ be a triangulation of $\Omegas_0$, the reference configuration of a three-dimensional solid structure, with $m$ nodes.
Let
\begin{equation}
  \label{eq:volumetric_fe_basis}
  S_h = \left\{\P_l(\s)\right\}_{l=1}^{m}
\end{equation}
be the standard interpolatory nodal basis functions in which each $\P_l(\s)$ is nonzero at exactly one node of the triangulation.
%~ In this work, these functions belong to the spaces of
%~ $Q^1$ bilinear (two dimensions) or trilinear (three dimensions) basis functions.
We let $\left\{\s_l\right\}_{l=1}^{m}$ and $\left\{\Y_l(t)\right\}_{l=1}^{m}$ be the coordinates of the triangulation's nodes in the reference and current configurations, respectively.
A continuous description of the configuration of the structure is defined by
\begin{equation}
  \label{eq:discrete_position}
  \Y_{h}(\s, t) = \sum_{l=1}^{m}\Y_{l}(t) \, \P_{l}(\s), \, \s \in \Omegas_0.
\end{equation}
%The reason is that the interpolatory~(nodal) basis functions associated with interior nodes 
%~ vanish on the surface. Similarly, the surface force density is determined by
 %~ \begin{equation}
%~ \F_{h} (\s,t) ={\displaystyle \sum_{l=1}^{M}\F_{l}(t)\phi_{l}(\s)},  \, \s \in \Gamma_0^{\fs}.
%~ \end{equation}
%~ Again, this sum only needs to be evaluated using the $M^\text{fs}$ surface nodes. 
Using the same FE basis functions, the deformation gradient is approximated as
\begin{equation}
 \label{eq:discrete_deformation} \FF_{h}(\s,t) = \frac{\partial \Y_{h}}{\partial \s}(\s,t) = \displaystyle \sum_{l=1}^{m}\Y_{l}(t) \, \frac{\partial \P_{l}}{\partial \s}(\s), \, \s \in \Omegas_0.
\end{equation}
This FE approximation to the deformation gradient tensor is used to compute 
the approximation to the first Piola-Kirchhoff stress tensor $\PP_{h}(\s,t)$ along with critical geometrical quantities.

Similarly, we construct the fluid-structure interface configuration described by the 
mapping $\X:(\Gamma_0^{\fs},t)\mapsto\Gamma_t^{\fs}$ by extracting the portion of the surface of the solid mesh that interacts with the fluid.
Let $m^\Gamma$ be the number of nodes on the surface of the triangulation and let
\begin{equation}
  \label{eq:surface_fe_basis}
  S_h^\Gamma = \{\P_l^\Gamma(\s)\}_{l=1}^{m^\Gamma}
\end{equation}
be the subset of $S_h$, the volumetric finite element space, which is nonzero on $\Gamma_0^{\fs}$.
Hence we define the configuration of the surface as
\begin{equation}
  \label{eq:discrete_deformation_surface}
  \X_{h}(\s, t) = \sum_{l=1}^{m^\Gamma}\X_{l}(t)\P_{l}^\Gamma(\s), \, \s \in \Gamma_0^{\fs}.
\end{equation}
For each $k=1, ..., m$ shape function in $S_h$, the standard Galerkin approach for the structural subproblem 
based on Eq.~(\ref{eq:weak_form3}) leads to the following evolution equation for the displacement of a single node:
\begin{equation}
\begin{split}
\rhos_0 \sum_{l=1}^{m} \left(\int_{\Omegas_0} \P_l(\s) \P_k(\s) \, {\mathrm d} \s \right)\frac{\partial^2 \Y_{l}}{\partial t^2}(\s,t) &=
- \int_{\Omegas_0} \, \PP_{h}(\s,t)  \, \cdot \, \grad_{\s} \, \P_k(\s) \,{\mathrm d} \s \\
& + \int_{\Gamma^{\fs}_0} \, j_{h}(\s,t) \, \vec{\tau}^{\text{f},\text{phys}}_{h}(\X_{h}(\s,t),t) \, \P_k(\s) \,\DA, 
\end{split}
\label{eq:weak_disc_form3}
\end{equation}
in which $\PP_{h}$ is a nonlinear function of $\FF_{h}(\s,t)$ through the constitutive relation as a function of the deformation gradient.
In our numerical scheme, we also compute the projection of the normal and tangential 
components of the surface force in Eq.~(\ref{eq:penalty_force}) per unit current area, $j_h^{-1} \F_h(\s,t)$, along 
the Lagrangian surface mesh. This is needed to specify the conditions for the pressure and the viscous stress and calculate 
the resulting exterior fluid traction $\vec{\tau}^{\text{f},\text{phys}}_{h}$.
Eq.~(\ref{eq:weak_disc_form3}) 
corresponds to a system of linear equations for the nodal accelerations, which we write as
\begin{equation}
 [ \vec{\CMcal{M}} ]\,[\vec{\ddot{\Y}}] = [ \vec{\CMcal{B}} ]
 \label{eq:algebraic_form}
 \end{equation}
in which $\vec{\CMcal{M}}$ is the mass matrix with entries $\vec{\CMcal{M}}_{i,j} = \rhos_0 \int_{\Omegas_0} \P_i(\s) \P_j(\s) \, {\mathrm d} \s$,
and $\vec{\CMcal{B}}$ is the applied load
vector, which here includes  
the exterior fluid traction supported on fluid-structure interface and
 the internal 
force density supported on the interior of the structure. 
Gaussian quadrature rules are used to approximate the integrals.
We use selective reduced integration \cite{malkus1978mixed} to avoid volumetric locking for nearly incompressible structural models.
%~ a third-order Gaussian quadrature is used to integrate the deviatoric stress $\PP_{\text{dev}}$ and 
 %~ first-order Gaussian quadrature is used for the dilatational stress $\PP_{\text{dil}}$.
%~ This is to avoid volumetric locking in nearly incompressible structures.
%~ to avoid volumetric locking in (nearly) incompressible structures, we use selective reduced integration \cite{malkus1978mixed}, in which 
%~ third-order Gaussian quadrature is used to integrate the deviatoric stress $\PP_{\text{dev}}$ and 
 %~ first-order Gaussian quadrature is used for the dilatational stress $\PP_{\text{dil}}$.
%~ We solve the variational problem using consistent mass matrix 
%~ on both sides, (i.e., no mass lumping on the left-hand side). See Wells et al. \cite{wells2021nodal} 
%~ for a discussion on the difference between using consistent and lumped mass matrices for projecting
%~ the divergence of the stress onto the finite element field.

\subsection{Fluid discretization}
\label{subsec:eulerian_discret}

We discretize the incompressible Navier-Stokes equations using a finite difference approximation on an adaptively refined marker-and-cell (MAC) 
staggered Cartesian grid \cite{griffith2009,griffith2012immersed}.
An unsplit, linearly-implicit discretization of the incompressible Navier-Stokes equations that 
solves time-dependent incompressible Stokes equations with 
the projection method used only as a preconditioner for the 
 flexible generalized minimal residual method (FGMRES) solver applied to the Stokes operator.
This allows for imposing physical boundary conditions along the boundaries of the Eulerian computational domain 
using methods detailed previously \cite{griffith2009,griffith2012immersed}. 
The divergence, gradient, and Laplace operators are discretized using compact second-order accurate differencing schemes.
The nonlinear advection terms in the Eulerian momentum equation are treated by
a staggered-grid variant of the piecewise parabolic method (PPM) \cite{griffith2012volume,griffith2012immersed, griffith2009}.

\subsection{Fluid-structure coupling}
\label{subsec:coupling_discretization}

We leverage a coupling scheme based on the immersed interface method for discrete surfaces  \cite{kolahdouz2020immersed}
that accounts for the physical jump conditions in the extended fluid stress along the fluid-solid interface 
and allows for the accurate evaluation of the exterior fluid traction.
To impose the jump in the stress, we first construct a continuous representation to the jump conditions in Eq.~(\ref{eq:ILE_traction_jump_codim1}) 
obtained by the $L^2$ projection along the surface mesh using the $S_h^\Gamma$ space defined by Eq.~(\ref{eq:surface_fe_basis}).
Geometrical quantities, including the surface normals and surface Jacobian determinant, that are needed by 
the IIM discretization are obtained by directly differentiating the discrete representation \cite{kolahdouz2020immersed,kolahdouz2021sharp}.
Briefly, given a scalar function $\psi \in L^{2}(\Gamma_0^{\fs})$, its $L^{2}$ projection $P_{h} \psi$ onto the subspace
$S_h^\Gamma$ is defined by requiring $P_{h} \psi$ to satisfy
\begin{equation}
  \int_{\Gamma_0^{\fs}} \big(\psi(\s) - P_{h}\psi(\s)\big) \, \P_l(\s) \,  {\mathrm d} A = 0, \quad \forall l=1,\ldots, m^\Gamma.
  \label{eq:L2_proj_def}
\end{equation}
Because the $L^{2}$ projection is defined via integration, the function $\psi$ does not need to be continuous or even to have well-defined nodal values.
By construction, however, $P_h \psi$ will inherit any smoothness provided by the subspace $S_h^\Gamma$.
In particular, for $H^1$-conforming Lagrangian basis functions, $P_h \psi$ will be at least continuous.
To solve for the projected jump conditions, linear systems of equations involving the boundary mass matrix ${\vec{\CMcal M}}^\Gamma$ need to be solved,
in which $\vec{{\CMcal M}}^\Gamma$ has
components $\vec{{\CMcal M}}_{i,j}^\Gamma = \int_{\Gamma_0^{\fs}} \P_i(\s) \P_j(\s) \, \DA$. The integral in
Eq.~(\ref{eq:L2_proj_def}) is evaluated using seventh-order Gaussian quadrature. 
Notice that these projections are computed only along the fluid-solid interface and involve only surface degrees of freedom. 
Consequently, the computational cost of evaluating these projections is 
asymptotically smaller than the solution of either the fluid or structure sub-problem.
Furthermore, this process can be repeated for each component of a vector-valued function independently.

To interpolate the discretized Eulerian velocity field $\u(\x,t)$  to the Lagrangian 
surface mesh, we use  a corrected  trilinear (or, in two spatial dimensions, bilinear) interpolation scheme that accounts for the known 
discontinuities of $\partial\vec{u}/\partial \vec{n}$ at intersection points with Cartesian finite difference stencils. 
The interpolated velocity field at intersection points is then projected
onto the space spanned by the nodal basis functions using Gaussian quadrature. 
A surface-restricted $L^2$ projection is also used to obtain nodal values  of the exterior fluid traction in Eq.~(\ref{eq:weak_disc_form3}), 
including the exterior fluid pressure and exterior viscous shear stress, see Kolahdouz et al.~\cite{kolahdouz2020immersed,kolahdouz2021sharp} for more details of these calculations.

In the Eulerian discretization, we account for the jump conditions in the pressure and the first normal derivatives of the velocity 
along the fluid-structure interface that occur in the ILE formulation by modifying 
the definitions of the pressure gradient and the viscous terms for those finite-difference stencils that cross the immersed 
interface using generalized Taylor series expansions; implementation details are discussed in prior work \cite{kolahdouz2020immersed, kolahdouz2021sharp}.

\subsection{Time integration}
\label{subsec:time_integration}

We define a vector $\vec{\Upsilon}$ that includes all of the Eulerian and Lagrangian quantities that need to be updated
as $\vec{\Upsilon}=[\u,  p, \X, \Y]$.
 We denote by $\dtf$ and $\dts$ the time steps associated with the fluid and solid mechanics solvers respectively, 
 for which we assume $\dtf = N_{\text{cycle}} \dts$ with an integer value of $N_{\text{cycle}} \geq 1$.
Starting from time step $n$ at time $t$ to time step $n+1$ at time $t + \dtf$,
we use the second-order Strang splitting scheme \cite{strang1968construction}, in which within three steps we: 
1) update the structure by solving the elastodynamic equations advancing over $N_{\text{cycle}}$ steps, each 
time solving with the time step $\dts/2$ and treating the fluid traction 
that loads the structure along $\Gamma^\text{fs}_t$ as fixed in time;
2) solve the IIM/FSI subproblem over a full time step $\dtf$, treating the structural configuration as constant in time; and
3) solve the elastodynamic equations over another $N_{\text{cycle}}$ steps each with the time step $\dts/2$, treating the fluid traction as fixed in time.
As described above, in this multirate (MR) approach, for 
 each half-time advancement of the structural solution in steps 1 and 3,  the time integration of the elastodynamics solver 
 can be subdivided into $N_{\text{cycle}}$ sub-steps.
This capability is particularly appealing for problems in which the stiffness of the 
structural model
requires structural  time step sizes that are substantially smaller than the time steps required to resolve the fluid dynamics.
  %~ In many FSI problems, the fluid solver is significantly more costly than the solid mechanics 
  %~ solution; therefore, this approach allows for taking large steps in the IIM/FSI solver 
  %~ while enabling finer time steps required to resolve the elastic time scales.
  %~ while controlling the stability through increasing the
  %~ the number of sub-steps in the solid mechanics solver. 
%~ Details of each time stepping for the solid mechanics solver and the IIM/FSI solver are given below.
We remark that unless otherwise mentioned, we set $N_{\text{cycle}}=1$ for results presented in Sec.~\ref{sec:results}.
   %~ We plan to explore the full capabilities of the multi-rate time-stepping in future work.

\subsubsection{Time integration of the elastodynamics equations}
\label{subsubsec:time_integration_solid}
%~ Although we ultimately employ a time step splitting approach that advances the solid configuration in two 
%~ half-steps $\dtf/2 = N_{\text{cycle}} \dts/2$, 
%~ to simplify the discussion, the 

The solution approach for a structural time step $\dts$ is detailed here.
%~ Additionally note that if the multi-rate time-stepping capability is turned on, 
%~ the elastodynamic equation is solved in a loop with $\dts =\dtf/N_{\text{cycles}}$ in which $N_{\text{cycles}}$ is the number of iterations.
Denote by $\cL_{\dts}$ the action of the elastodynamics solver over the time increment $\dts$ that
acts on a solution vector $\vec{\Upsilon}$. This solution vector includes all of the Eulerian and Lagrangian variables, but only advances 
the volumetric structural variables $\Y$ and $\PP$ while keeping the remaining variables 
fixed.
The discretized equation, Eq.~(\ref{eq:weak_disc_form3}), is solved over the time increment $\dts$ to obtain
${\Y^{n+1}}$ using a modified trapezoidal scheme. Defining $\vec{U}={\partial \Y}/{\partial t }$, the time-integration 
is performed in two steps. In the first step, we use a modified Euler method and obtain predicted values
\begin{align}
[\tilde{\vec{U}}^{n+1}] &= [\vec{U}^{n}] + \dts \, [\vec{{\CMcal{M}}^{-1}}^{}] \, [ \vec{\CMcal{B}}^{n} ],\\
{[\tilde{\Y}^{n+1}]} &= [\Y^{n}] + \dts \, [\tilde{\vec{U}}^{n+1}].
        %~ //    U^{n+1,*} := U^{n} + (dt/rho) F^{n}
        %~ //    X^{n+1,*} := X^{n} + dt       U^{n+1,*}
\end{align}
Using $[\tilde{\vec{U}}^{n+1}]$ and ${[\tilde{\Y}^{n+1}]}$, we calculate $\tilde{\vec{\CMcal{B}}}^{n+1}$ (recall
$\vec{\CMcal{B}}$ includes the terms on the right-hand-side of Eq.~(\ref{eq:discrete_deformation})). In the second (corrector) step,  
 we use a modified (explicit) trapezoidal rule, so that
\begin{align}
[\vec{U}^{n+1}] &=  [\vec{U}^{n}] + \frac{\dts}{2} \, [{\vec{\CMcal{M}}^{-1}}^{}] \, ([\vec{\CMcal{B}}^{n}] + [\tilde{\vec{\CMcal{B}}}^{n+1}]),\\
[{\Y^{n+1}}] &= [\Y^{n}] + \frac{\dts}{2} \, ([{\vec{U}}^{n}] + [{\vec{U}}^{n+1}]).
        %//    U^{n+1} := U^{n} + (dt/(2 rho)) (F^{n} + F^{n+1,*})
        %//    X^{n+1} := X^{n} + (dt/2)       (U^{n} + U^{n+1})
\end{align}
%~ Other explicit time integration schemes such as a second-order midpoint rule or a second-order
%~ strong-stability preserving Runge-Kutta time-steppers (SSP-RK) can be also implemented.
Note that in practice, we never directly evaluate $[\vec{{\CMcal{M}}}^{-1}]$, but instead approximately
evaluate its action via a Krylov solver. In particular, 
we use the symmetric iterative solver MINRES,
given the positive definiteness of the mass matrix.

\subsubsection{IIM time integration scheme}
\label{subsec:time-integration}

Denote by $\cE_{\dtf}$ the action of the IIM solver over a full time step that
acts on the solution vector $\vec{\Upsilon}$ that includes all of the Eulerian 
and Lagrangian variables but only advances $\u$,  $p$, and $\X$.
Starting from $\X^{n}$ and $\vec{u}^{n}$ at time $t^n$ and $p^{n-\frac{1}{2}}$ at 
time $t^{n-\frac{1}{2}}$, we compute $\X^{n+1}$, $\u^{n+1}$, and $p^{n+\frac{1}{2}}$. 
Briefly, using the discrete velocity restriction operator $\cJ^{n}=\cJ[\X^{n},\F^{n}]$, we first obtain 
initial approximations to the interface position at time $t^{n+\frac{1}{2}}$  via
\begin{equation}
\label{eq:trap1}	\X^{n+\frac{1}{2}} = \X^{n} + \frac{\dtf}{2} \cJ^{n} \vec{u}^{n}. \\
%~ \label{eq:trap2}	\X^{n+\frac{1}{2}}=\frac{\hat{\X}^{n+1}+\X^{n}}{2}.
\end{equation}
Next, we solve for $\X^{n+1}$, $\u^{n+1}$, and $p^{n+\frac{1}{2}}$ in Eqs.~(\ref{eq:ILE_momentum})--(\ref{eq:ILE_immersed_kinematic_cond}) via
\begin{align}
\label{eq:time-integ-momentum}\rho\left(\frac{\u^{n+1}-\u^{n}}{\dtf}+\vec{A}^{n+\frac{1}{2}}\right) &=-\vec{G} \, p^{n+\frac{1}{2}}+\muf \, \vec{L}\left(\frac{\u^{n+1} + \u^{n}}{2}\right)+ \vec{\cS}^{n+\frac{1}{2}}\, \F^{n+\frac{1}{2}}, \\
\vec{D} \cdot \vec{u}^{n+1} &=0, \\
\label{eq:update-X-FSI} \frac{\X^{n+1}-\X^{n}}{\dtf} &= \U^{n+\half} = \cJ^{n+\half} \left(\frac{\u^{n+1}+\u^{n}}{2}\right),
\end{align}
in which $\vec{A}^{n+\frac{1}{2}}=\frac{3}{2}\vec{A}^{n}-\frac{1}{2}\vec{A}^{n-1}$ is obtained from a high-order upwind spatial discretization
of the nonlinear convective term $\u \cdot \grad \u$ \cite{griffith2009}. In Eq.~(\ref{eq:time-integ-momentum}),
 $\vec{\cS}^{n+\frac{1}{2}}$ is the force-spreading operator at the half time step configuration and $\vec{\cS}^{n+\frac{1}{2}}\, \F^{n+\frac{1}{2}}$ is
the discrete Eulerian body force on the Cartesian grid that corresponds to the sum of all the correction terms due to physical jump conditions computed
using the Lagrangian force $\vec{F}^{n+\frac{1}{2}} = \vec{F}[\X^{n+\frac{1}{2}}, \U^{n}, t^{n+\half}]$.
 %~ $\vec{\cS}$ is the force-spreading operator that takes the form correction terms due to physical jump conditions, and 
 In Eq.~(\ref{eq:update-X-FSI}),
$\cJ^{n+\half}=\cJ[\X^{n+\half},\F^{n+\half}]$ is the velocity restriction operator at the half time step configuration.
$\vec{G}$, $\vec{L}$, and $\vec{D}\cdot\mbox{}$ respectively represent discrete gradient, Laplacian, and divergence operators. 
This time stepping scheme requires only linear solvers for the time-dependent incompressible Stokes equations.
We solve this system of equations by the FGMRES algorithm with a preconditioner based on the projection method 
that uses inexact subsolvers for the velocity and pressure \cite{griffith2009}.
In the initial time step, a two-step predictor-corrector method is used to determine the velocity, deformation, and pressure; 
see Griffith and Luo \cite{BEGriffith17-ibfe} for further details.

\subsubsection{Fluid-structure interaction time stepping scheme}

A loosely coupled scheme is considered for the time integration of the equations.
For the solution vector $\vec{\Upsilon}$ over a fluid time step $\dtf$ using the Strang splitting scheme we have,
\begin{equation}
\vec{\Upsilon}^{n+1} = \cL_{\dtf/2} \cE_{\dtf}\cL_{\dtf/2}\vec{\Upsilon}^{n}.
\end{equation}
With this time-staggered approach, the overall algorithm to solve the fluid-structure problem 
is

%~ \begin{algorithm}
    %~ \SetKwInOut{Input}{Input}
    %~ \SetKwInOut{Output}{Output}

    %~ \underline{function Euclid} $(a,b)$\;
    %~ \Input{Two nonnegative integers $a$ and $b$}
    %~ \Output{$\gcd(a,b)$}
    %~ \eIf{$b=0$}
      %~ {
        %~ return $a$\;
      %~ }
      %~ {
        %~ return Euclid$(b,a\mod b)$\;
      %~ }
    %~ \caption{Overall time-Stepping Scheme for advancing the solution vector from $\vec{\Upsilon}^{n}$ to $\vec{\Upsilon}^{n+1} $}
%~ \end{algorithm}

\begin{steps}
	\item{ 
		Solve the elastodynamic equation Eq.~(\ref{eq:algebraic_form})
		over a half time step $\dtf/2=N_{\text{cycle}} \dts/2$ using the modified trapezoidal scheme in Eqs.~(\ref{eq:trap1})--(\ref{eq:update-X-FSI}) 
		to advance from step $n$ to $n+\frac{1}{2}$ and obtain 
		${\Y^{n+\frac{1}{2}}}$. 
		If multi-rate time-stepping capability is used, the elastodynamic equation is solved with $\dts/2$ over a loop of size $N_{\text{cycle}} > 1$ with
		$N_{\text{cycle}}=\dtf/\dts$ number of steps.
		%and ${\PP^{n+\frac{1}{2}}}$ 
		%~ ${\disp_{\com}}^{n+\frac{1}{2}}$, ${\dispdot_{\com}}^{n+\frac{1}{2}}$, ${\vec{\omega}}^{n+\frac{1}{2}}$ and ${\QQ}^{n+\frac{1}{2}}$,
		%~ and the new position of the volumetric Lagrangian mesh.
		}
	\item{ 
		Calculate the penalty force using the most recent position of the volumetric structural mesh and the surface mesh
		that moves with the fluid. 
	}
	\item{ 
		Solve for the IIM subsystem in Eqs.~(\ref{eq:ILE_momentum})--(\ref{eq:ILE_immersed_kinematic_cond}) over a full time step 
		and obtain the updated Eulerian velocity field $\vec{u}^{n+1}$ and pressure $p^{n+\frac{1}{2}}$ as well as
		the Lagrangian velocity $\U^{n+\frac{1}{2}}$ and positions $\X^{n+1}$ of the surface mesh and the exterior 
		fluid traction forces $\vec{\tau}^{\text{f},\text{phys}}$.
	}
	\item{ 
		Using the exterior fluid traction force $\vec{\tau}^{\text{f},\text{phys}}$ from Step III, solve the elastodynamic equation 
		Eq.~(\ref{eq:algebraic_form})
		over a half time step $\dtf/2=N_{\text{cycle}} \dts/2$ using the modified trapezoidal scheme in Eqs.~(\ref{eq:trap1})--(\ref{eq:update-X-FSI}) 
		to advance from step $n+\frac{1}{2}$ to $n+1$ and obtain 
		$\Y^{n+1}$.
		If multi-rate time-stepping capability is used, the elastodynamic equation is solved with $\dts/2$ over a loop of size $N_{\text{cycle}}>1$ with
		$N_{\text{cycle}}=\dtf/\dts$ number of steps.
	}
	%~ \item{ 
			%~ Move the Lagrangian mesh of the bulk solid and obtain the new positions $\Y^{n+1}$.	}
\end{steps}

Our numerical experiments 
indicate that we can find values of 
$\dtf$ 
that are small enough so that the 
scheme is stable, 
as was also suggested in a 
recent analysis by Hua and Peskin \cite{hua2022analysis} for the standard immersed boundary method that uses the same form of 
		 		approximate Lagrange multiplier force. 
Additionally, notice that the ILE formulation easily allows for taking smaller solid time-step
				in the explicit finite element solid mechanics solver.
				This will be an important feature for problems that pose larger stiffness in the solid domain than the fluid and 
				is achieved through increasing the $N_{\textrm{cycle}}$ parameter without changing the 
				 overall fluid time-step or penalty parameters.

\subsection{Software implementation}
\label{sec:implementation}

The ILE approach has been implemented in the open-source IBAMR software \cite{ibamr}, a C++ framework for FSI modeling using immersed formulations.
IBAMR provides support for large-scale simulations through the use of distributed-memory parallelism and adaptive mesh refinement (AMR).
IBAMR relies on other open-source software libraries, 
including SAMRAI \cite{samrai, hornung2002managing}, PETSc \cite{petsc-web-page}, \textit{hypre} \cite{hypre,falgout2002hypre}, 
and \texttt{libMesh} \cite{Kirk:2006:LCL:1230680.1230693}.
For example,  mesh handling tools in \textit{libMesh} allows for the use of separate data structures for the volumetric and the corresponding
 surface structural meshes throughout the entire computation by storing the so-called \emph{interior parent} elements.
 All the capabilities of the present approach and benchmark examples are included in the IBAMR software framework, version 0.11.0 or newer.

\section{Numerical results}
\label{sec:results}

We use \textbf{Q1} structural elements in our computational tests (bilinear elements in two spatial dimensions 
and trilinear elements in three spatial dimensions). 
The grid spacing on the finest level of the hierarchical, locally refined Cartesian grid is $h_{\text{finest}}=r^{-(N-1)}h_{\text{coarsest}}$,
in which $h_{\text{coarsest}}$ is the coarsest level grid spacing, $r$ is the refinement ratio, and $N$ is the number of refinement levels.
Unless otherwise noted, the initial ratio of the Lagrangian element size of the structural meshes to the Eulerian grid spacing, denoted by $\Mfac$,
is approximately $\Mfac \approx 2$ along the fluid-structure interface.
In all computations, we require that the maximum absolute displacement between the 
two representations of the fluid-structure interface is no greater than $0.1$ of the Eulerian meshwidth.
Note that to achieve $\|\vec{\xi}(\s,t)-\vec{\chi}(\s,t)\| = O(h^2)$, the spring penalty parameter needs to satisfy $\kappa = O(1/h^2)$.
Additionally, we take the damping parameter to satisfy $\eta = O(1/h)$.
Under grid refinement, we expect that the penalty force
 $\F = \kappa\left(\Y(\s,t)-\X(\s,t)\right) + \eta \,(\frac{\partial \Y}{\partial t}(\s,t) - \frac{\partial \X}{\partial t}(\s,t))$ 
will converge to the physical loading forces 
and satisfy $\|\F\| = O(1)$.
In particular, for our convergence tests, we maintain $\dt = O(h)$  to keep the advective Courant-Friedrichs-Lewy (CFL) number fixed under grid refinement. 
It is then convenient to choose $\kappa = \kappa_0 / \dt^2$ and also optionally choose $\eta = \eta_0 / \dt$.
 $\kappa_0$ is determined by, first, choosing the finest grid spacing to be considered, say $h_\text{min}$, 
along with the corresponding time step size $\dt_\text{min}$.
We then empirically determine approximately the largest value of $\kappa_0$ that allows the scheme to remain stable with $h = h_\text{min}$ and $\dt = \dt_\text{min}$.
We then use the prescribed relationship between $\kappa$ and $\dt$ to determine $\kappa$ (and $\eta$) for all coarser cases.
This ensures that the numerical parameters are stable for all grid spacings considered in the convergence test while using scalings 
that ensure that the method achieves its targeted convergence rate.
Note that in our explicit time stepping approach to coupling the fluid and solid degrees of freedom, we expect
a largest stable value of spring penalty parameter $\kappa$ for fixed spatial and temporal discretization parameters.
In a given problem, if there is already a physical estimate for the total fluid forces acting on the fluid-structure interface, 
this estimate along with a targeted discrepancy tolerance $\epsilon_\text{tol}=\norm{\Y(\s,t)-\X(\s,t)}_{\infty}$, 
	can help with an initial guess for the spring stiffness, e.g., via $\kappa \approx \| \vec{F} \| / \epsilon_\text{tol}$. 
		 		Regardless of whether a reasonable estimate is available however, in practice 
		 		 we find the approximation for
the $\kappa$ close to the stability limit by tuning the parameter using the method of bisection.
Damping appears to be optional, and a slight amount of damping can be added while finalizing the tuning of $\kappa$.
All computations use a tight relative convergence threshold of $\epsilon_\text{rel}=\text{1e-10}$ for all iterative linear solvers.
%~ For consistency with the scaling description of penalty parameters in Sec.~\ref{sec:kappa_scaling}, 
%~ For time-stepping choice of $\dt \propto h_{\text{finest}}$, we express the spring penalty parameter in 
%~ form of $\kappa=C_0/\dt$, consistent with the force per unit volume unit of $\kappa$. 
%~ Similarly, we express the damping parameter as $\eta=C_1$.
%~ In these relations, $C_0$ and $C_1$
 %~ are constant values that need to be tuned using the approach described in Sec.~\ref{sec:kappa_scaling}.
To avoid volumetric locking in incompressible cases, we use selective reduced integration \cite{malkus1978mixed}, in which 
third-order Gaussian quadrature is used to integrate the deviatoric stress $\PP_{\text{dev}}$ and 
 first-order Gaussian quadrature is used for the dilatational stress $\PP_{\text{dil}}$.
For simplicity in the rest of the paper we drop the superscript ``f'' from the fluid time step and 
use the general notation $\Delta t \equiv \dtf$ and $N_{\text{cycle}}=1$, unless otherwise mentioned.
Tests include two dimensional cases in Secs.~\ref{subsubsec:soft_disc}, \ref{subsubsec:Turek_Hron} and \ref{subsubsec:fluid_reservoir}
and three dimensional cases in Secs.~\ref{subsubsec:Hessenthaler_model} and \ref{subsec:clot_ivc}.

\subsection{Soft disk in a lid driven cavity \protect\footnote{This benchmark test is provided in IBAMR version 0.11.0 or newer, within the directory \texttt{examples/IIM/ex7}.}}
\label{subsubsec:soft_disc}

This benchmark test involves a soft structure in a lid-driven cavity flow. Slightly different versions of this example have been adapted in previous studies 
\cite{wang2010interpolation, roy2015benchmarking, BEGriffith17-ibfe}.
The computational domain is $\Omega =[0.0 \, \text{cm}, 1.0 \, \text{cm}]\times[0.0 \, \text{cm}, 1.0 \, \text{cm}]$,
a square of size $L_x = L_y = 1.0 \, \text{cm}$. The immersed structure is initially a disk of radius $R= 0.2 \, \text{cm}$, centered at $(0.5 \, \text{cm}, 0.6 \, \text{cm})$.
A uniform velocity $\vec{u}=(u_{\infty} = 1 \, \text{cm}\cdot\text{s}^{-1} ,0 \, \text{cm}\cdot\text{s}^{-1})$ is imposed along the top boundary and 
zero-velocity conditions are imposed at all other boundaries.
The fluid has a uniform density $\rhof =1.0$~$\text{g}\cdot\text{cm}^{-3}$ and dynamic viscosity $\muf = 0.01$~$\text{g}\cdot(\text{cm}\cdot \text{s})^{-1}$.
The soft structure has a shear modulus $\Gs=0.1~\text{dyn}\cdot\text{cm}^{-2}$ 
and Poisson's ratio $\nu = 0.499$.
%~ we consider three different Poisson's ratios of $\nu = 0.49, \, 0.499$ and $0.4999$.
The density ratio is set to $\rhos_0/\rhof=1$. 
We consider the time span $0.0\,\text{s}\, \leq t \leq 10.0 \,\text{s}$, during which the disk goes through slightly more than one rotation inside the cavity.
The Eulerian domain is
 discretized using a refinement ratio of $r=2$ with $N=2$ Cartesian grid levels with a grid spacing of $h_{\text{coarsest}}=\frac{L_x}{48}$ on the coarsest 
  level and $h_{\text{finest}}=\frac{L_x}{48 \times 2}$  on the finest level.
  The time step size is $\dt=(0.05 \, \text{s}\cdot\text{cm}^{-1})\,h_\text{finest}$, and the force penalty parameters are set 
  to $\kappa=(1\times 10^{-2} \, \text{g}\cdot\text{cm}^{-2})/{\dt}^2$ and $\eta=(5 \times 10^{-2} \, \text{g}\cdot\text{cm}^{-2})/\dt$.
  Unless otherwise noted, these forms of penalty parameters are used for all the following tests related to the soft disk in a lid driven cavity.
%~ to $\kappa=(1.92 \times 10 \, \text{g}\cdot\text{cm}^{-2}\cdot s^{-1})/{\dt}$ and $\eta= 1.0 \, \text{g}\cdot\text{cm}^{-2} \cdot {s}^{-1}$.
 %~ \begin{figure}[hb!!!!]
		%~ \centering
			%~ \includegraphics[width=0.8\textwidth]{fig/soft_disk_contours.pdf}
		%~ \caption{The computed velocity field and motion of the soft disk in a lid driven cavity 
		%~ at times (a) $t = 0.0$~s, (b) $t=2.0$~s, 
		%~ (c) $t=3.5$~s, (d) $t= 5.25$~s, (e) $t=6.5$~s, and (f) $t=9.0$~s. }	
		%~ \label{fig:soft_disk_contours} 
%~ \end{figure}

\begin{figure}[t!!!!]
		\centering
			\includegraphics[width=0.9\textwidth]{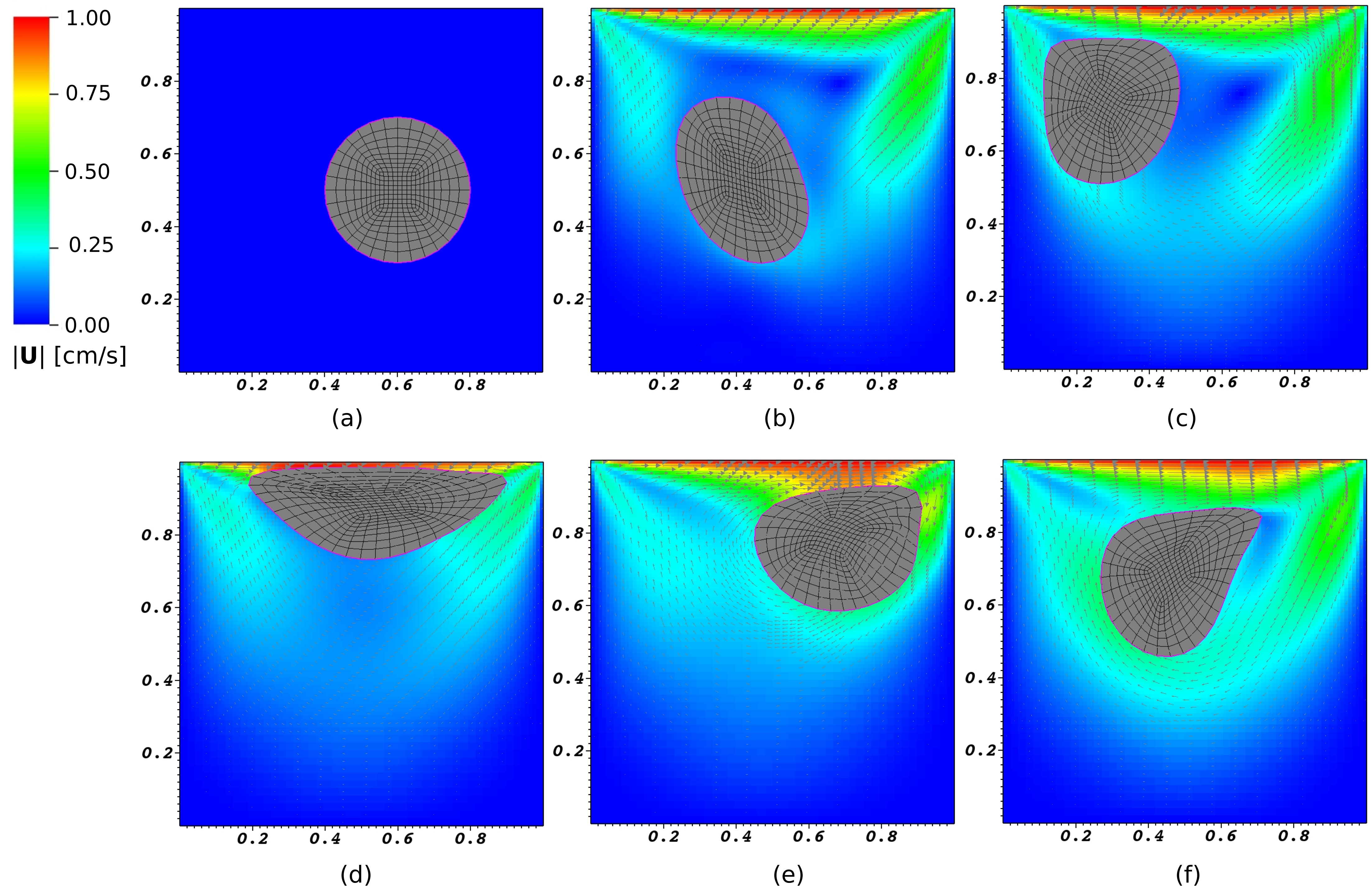}
		\caption{The computed velocity field and motion of the soft disk in a lid driven cavity 
		at times (a) $t = 0.0$~s, (b) $t=2.0$~s, 
		(c) $t=3.5$~s, (d) $t= 5.25$~s, (e) $t=6.5$~s, and (f) $t=9.0$~s. }	
		\label{fig:soft_disk_contours} 
\end{figure}

Fig.~\ref{fig:soft_disk_contours} shows the snapshots and the corresponding velocity fields 
at six different times. The mesh deformation of the volumetric solid structure is plotted along 
with the surface representation (colored in pink), which clearly shows that the two are conformal in their motion.
Additionally, notice that the induced flow by the driven lid brings the structure into near contact with the upper boundary of the domain.
 Starting from approximately time $t= 3.5 \, \text{s}$ in Fig.~\ref{fig:soft_disk_contours}(c) to Fig.~\ref{fig:soft_disk_contours}(e), 
 the entrained structure undergoes very large deformation in a near-contact condition with the domain's boundary. 
 Fig.~\ref{fig:soft_disk_contours}(d) shows how the near-contact dynamics are captured using the ILE algorithm.
\begin{figure}[b!!]
		\centering
			\includegraphics[width=0.75\textwidth]{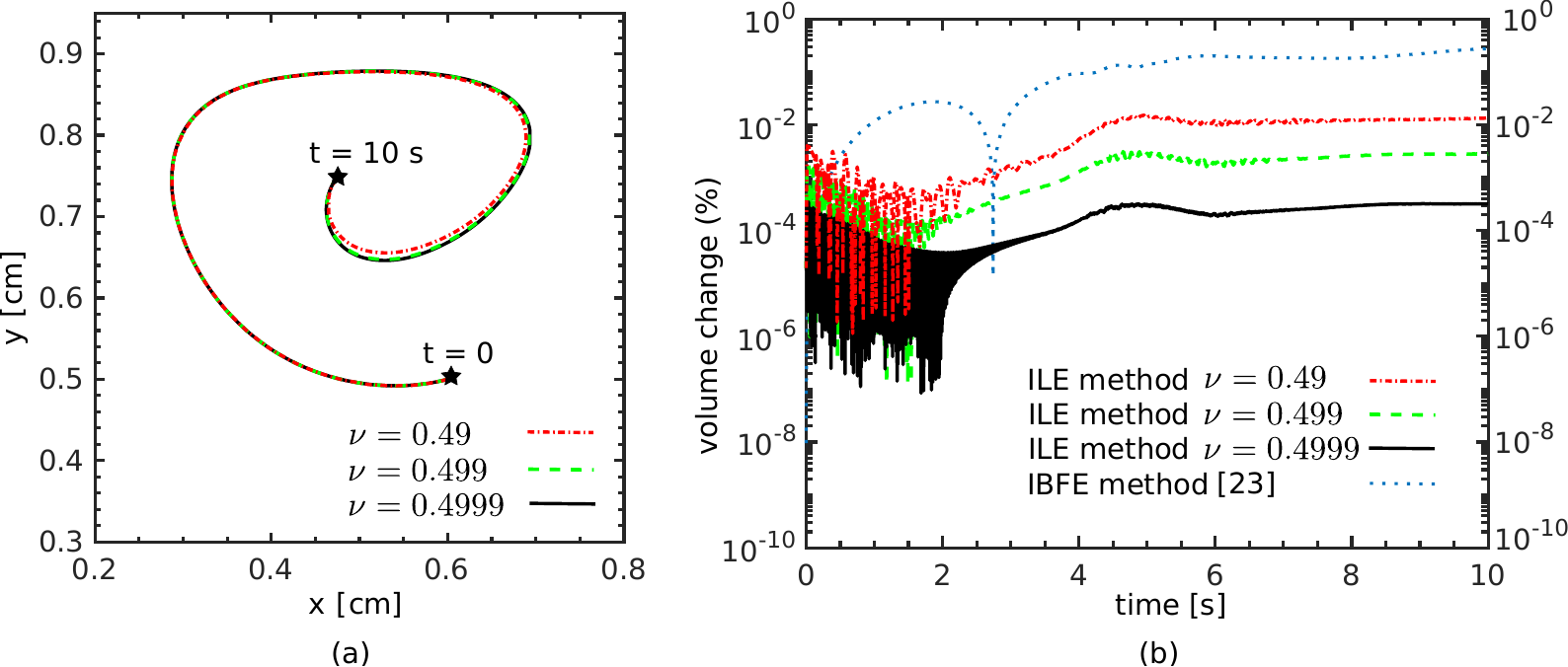}
		\caption{(a) Trajectory of the centroid of the soft disk over time for three different values of the Poisson's ratio $\nu$. 
				(b) The percentage of volume change over time for three values of the Poisson's ratio $\nu$ along with 
				the percentage of volume change using the immersed finite element/difference method \cite{BEGriffith17-ibfe}. }	
		\label{fig:soft_disk_trag_vol} 
\end{figure}
 \begin{figure}[ht!!]
		\centering
			\includegraphics[width=\textwidth]{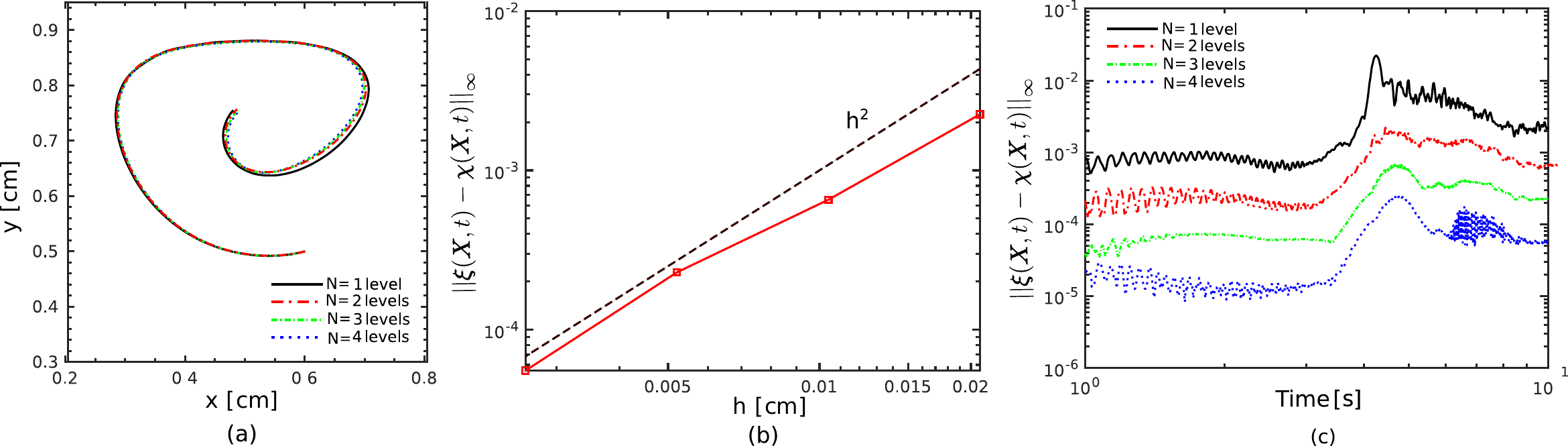}
		
		\caption{
		(a) The trajectory of the soft disk's centroid for a series of grid resolutions.
		(b) Norm of the discrepancy in the fluid-structure interface configurations at the final time $t=10$~s.
		(c) Time history of the ${\Linf}$ difference 
		between the positions of the 
		two representations of the fluid-structure interface over time.}	
		\label{fig:Lmax_flexible_2D} 
\end{figure}

To investigate 
how our use of a nearly incompressible structural formulation as a penalty method for the exactly incompressible case
 affects volume conservation, we consider
two additional Poisson's ratios $\nu = 0.49$ and $0.4999$. 
Fig.~\ref{fig:soft_disk_trag_vol} shows the trajectory 
of the centroid of the soft disk along with the percentage of volume change over time.
 For comparison of the volume conservation, 
 we compare against a solution based on the immersed finite element/difference method reported by Griffith and Luo \cite{BEGriffith17-ibfe}
 using a three-point regularized delta function kernel and the same time step size and grid spacing as the ILE cases.
  Notice in Fig.~\ref{fig:soft_disk_trag_vol}(a) that there is an excellent agreement between the trajectories 
  of the case with $\nu=0.499$ and $\nu=0.4999$, whereas slight differences are observed for the case with $\nu=0.49$.
   However, as depicted in Fig.~\ref{fig:soft_disk_trag_vol}(b), even for $\nu=0.49$, 
   the ILE formulation yields substantially improved volume conservation 
  compared to the immersed finite element/difference method, generating volume conversation errors that are at least two orders of magnitude smaller.
   Volume conservation further improves as $\nu \to 0.5$. 
   Given these observations, it appears that a choice of $\nu=0.499$ is a good compromise 
   for our FSI simulations considering both the convergence of the centroid's trajectory, overall volume conservation,
   and time step size restrictions associated with increased bulk moduli.
    (Because we use a fully explicit time stepping scheme for the structural dynamics, as $\nu \rightarrow 0$, $\Delta t^\text{s} \downarrow 0$.)

As a further verification, a mesh refinement study using a series of three grid resolutions is conducted.
The grid refinement is performed by changing the number of AMR levels.
The case considered thus far, with $N=2$, is used as an intermediate resolution, along with locally refined grids with $N=1$, $3$, and $4$ levels. 
The Lagrangian meshes are consistently refined (or, for $N=1$ case, coarsened) with the Eulerian grid to maintain $\Mfac \approx 2$.  
The same grid-dependent time step size and penalty parameters are used as detailed above. 
Fig.~\ref{fig:Lmax_flexible_2D}(a) shows the trajectories of the centroid of the soft disk for all the four grid sizes.
The trajectory clearly converges under grid refinement.
Further, for the choice of penalty parameters used here, we expect 
pointwise second-order convergence for the discrepancy between the positions of
the two interface representations. 
Fig.~\ref{fig:Lmax_flexible_2D}(b) shows $\norm{\Y(\s,t)-\X(\s,t)}_{\infty}$,
%~ which is the $\Linf$ norm of the discrepancy between the Lagrangian quadrature points at the surface 
%~ mesh and the corresponding points on the boundary of the volume mesh, computed
at the final time $t=10$~s. Second-order convergence in the maximum norm is apparent.
Fig.~\ref{fig:Lmax_flexible_2D}(c) further details the change of the $\Linf$ norm of the discrepancy over time, 
%~ Note that to achieve $\|\vec{\xi}(\s,t)-\vec{\chi}(\s,t)\| = O(h^2)$, the penalty parameters need to satisfy $\kappa = O(1/h^2)$ and $\eta = O(1/h)$, 
%~ because the penalty force $\F = \kappa\left(\Y(\s,t)-\X(\s,t)\right) + \eta \,(\frac{\partial \Y}{\partial t}(\s,t) - \frac{\partial \X}{\partial t}(\s,t))$ 
%~ satisfies $\|\F\| = O(1)$ under grid refinement \cite{kolahdouz2020immersed}. 
%~ This justifies our aforementioned choices of the time step size and penalty parameters and their relation with the grid size $h_\text{finest}$ 
%~ for the purpose of performing the grid convergence.    
 \begin{figure}[hb!!]
		\centering
			\includegraphics[width=0.75\textwidth]{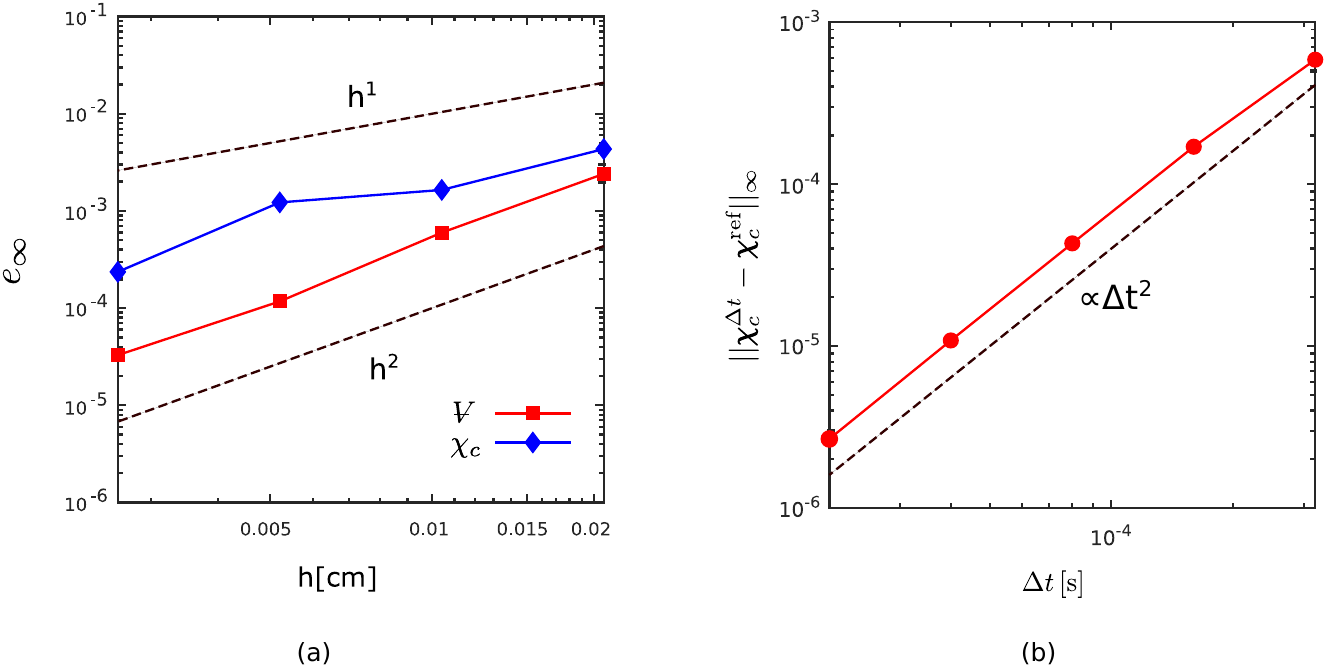}
		\caption{
		(a) The spatio-temporal convergence rates of the volume and centroid of the soft disk in the
		lid driven cavity test case over the 
		time interval $0 \,\text{s} \leq t \leq 10 \,\text{s}$.
		(b) The temporal convergence of the disk's centroid for a fixed grid with $N=3$ at $t=10$~s.
		Second-order convergence is apparent.}
		\label{fig:vol_x_dt_conv} 
\end{figure}
   \begin{figure}[t!!]
		\centering
			\includegraphics[width=0.8\textwidth]{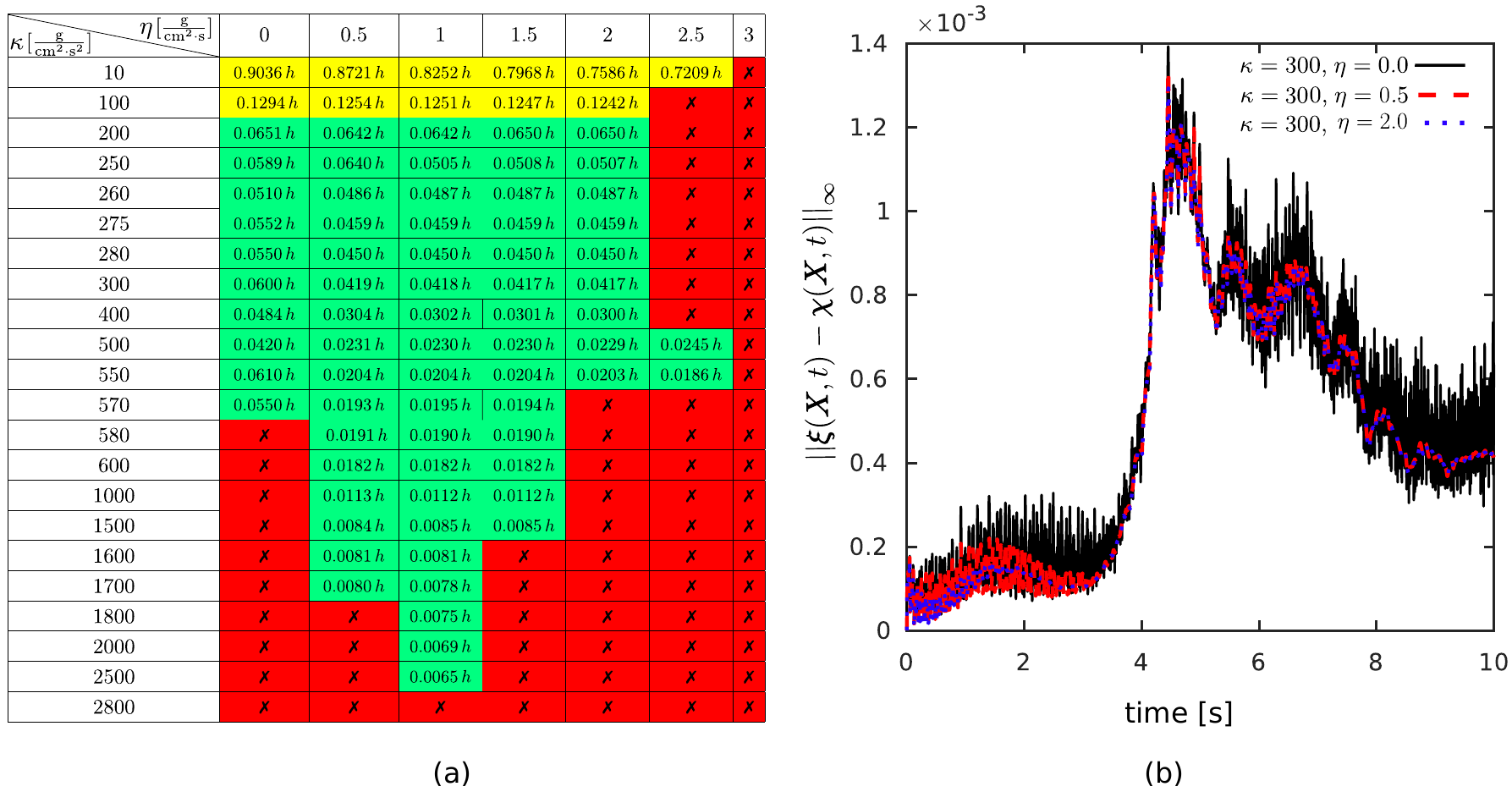}
		
		\caption{
		%~ (a) Limits of stable $\kappa$ values ($\eta = 0$) with the choice of the time step $\dt=(0.05 \, \text{s}\cdot\text{cm}^{-1})\,h_\text{finest}$
		 %~ and maximum discrepancy $\norm{\Y(\s,t)-\X(\s,t)}_{\infty} \leq 0.1 \, h_\text{finest}$.
		 %~ A clear proportionality of $\kappa \propto \dt$ is observed, consistent with the overall scaling of $\kappa \propto h/\dt^2$ in the continuous limit.
(a) A table of stable vs. unstable penalty pairs $(\kappa, \eta)$
 with red color cells marked by $\xmark$ symbol showing unstable cases. The numbers reported for stable pairs is the 
 maximum discrepancy between the two Lagrangian representations $\norm{\Y(\s,t)-\X(\s,t)}_{\infty}$ at $t=10 \, \text{s}$, shown as a factor of 
 $h \equiv h_\text{finest}$. Green cells are the acceptable cells with the additional condition of $\norm{\Y(\s,t)-\X(\s,t)}_{\infty} \leq 0.1 \, h_\text{finest}$ imposed 
 in our simulations.
   A non-zero damping parameter does not seem to be required to achieve a stable result with an acceptable maximum discrepancy.
   Additionally, note that there can always be a small enough time-step for which the currently unstable pairs would become stable. 
  (b) The maximum discrepancy between the two Lagrangian representations over time for 3 cases from the table  with fixed
  $\kappa=300$ and different choices of $\eta$.
  Although the oscillations are all smaller than $0.1 \, h$, 
 a small amount of non-zero damping is shown to further reduce spurious numerical oscillations.}	
		\label{fig:kappa_eta_disk} 
\end{figure}
   \begin{figure}[b!!]
		\centering
			\includegraphics[width=0.7\textwidth]{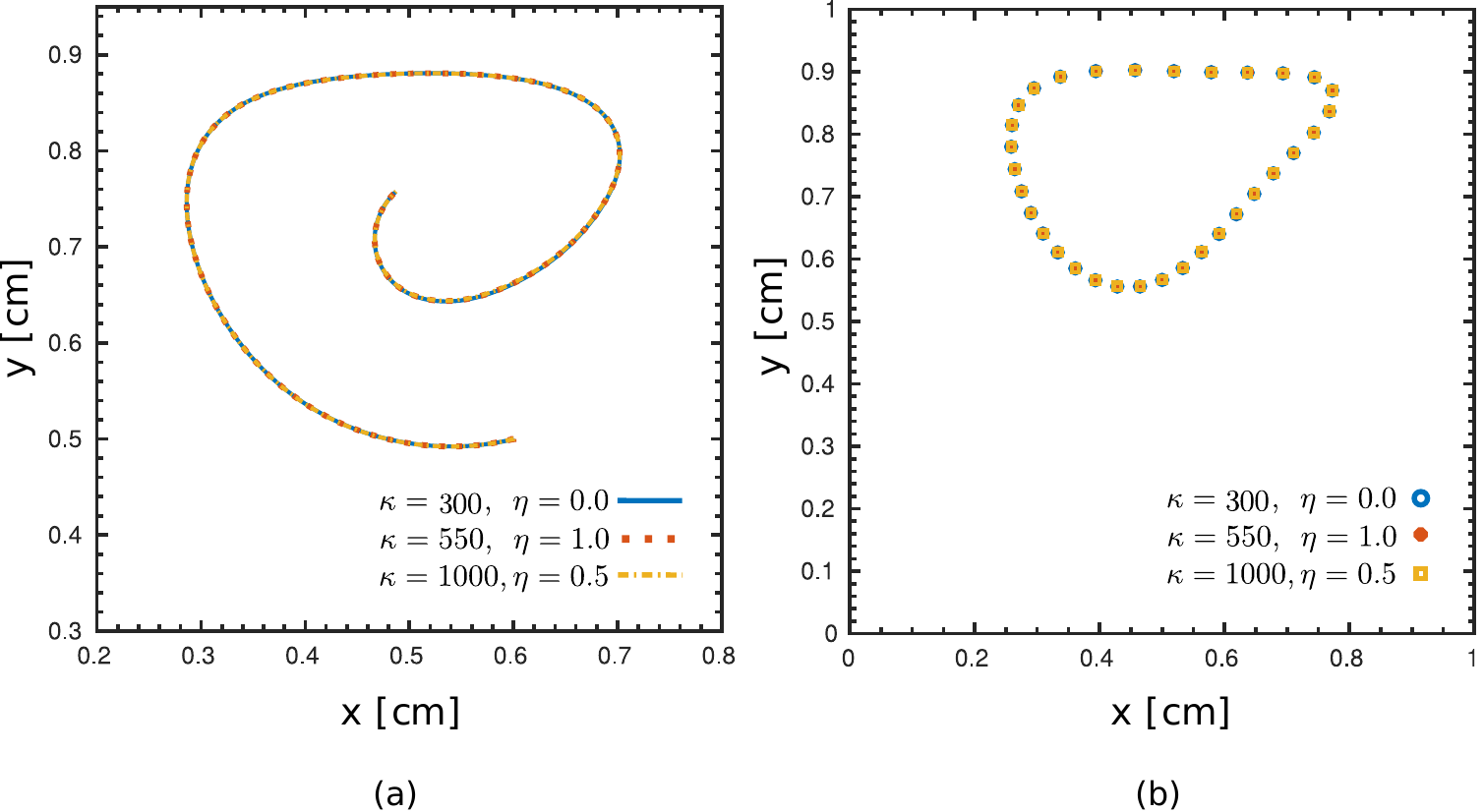}
		
		\caption{
		(a) The trajectory of the centroid of the soft disk for three selected pairs of $(\kappa, \eta)$ 
		within the stable cells in Fig.~\ref{fig:kappa_eta_disk}(b), including a case with $\eta =0$.
		(b) Nodal positions of the surface of the structure for the same pairs of $(\kappa, \eta)$ at $t=10$~s.
		 The results are largely insensitive to the values of penalty parameter values considered here.}	
		\label{fig:trajectory_softdisk} 
\end{figure}

%~ (a) Heat map of the penalty pair $(\kappa, \eta)$ showing the region with stable results. 
		%~ The green color shows the region with stable penalty pairs while the red color points to the unstable region.

In addition to the pointwise Lagrangian displacement, we have computed empirical ${\Linf}$ convergence rates for the structural volume of the soft disk over 
the entire course of the simulation.
%~ as well as 
%~ the Lagrangian pressure and wall shear stress at final time.
Because there is no exact solution available for the full dynamics of this problem,
we use a Richardson extrapolation procedure 
to estimate the convergence rate $q_{\infty}$ in the ${\Linf}$ norm of the error in the  soft disk's centroid $\X_c$, %volume $\volume$,
\begin{equation}
q_{\infty}[\X_c;N_{\text{coarse}}]=\log_2\left(\frac{e_{\infty}[\X_c;N_{\text{coarse}}]}{e_{\infty}[\X_c;N_{\text{intermediate}}]}\right)
\label{eq:conv_v}
\end{equation}
in which $e_{\infty}[\X_c;N_{\text{coarse}}]$ and $e_{\infty}[\X_c;N_{\text{intermediate}}]$ are defined as
\begin{align}
e_{\infty}[\X_c;N_{\text{coarse}}]=\max\limits_{0 \,\text{s} \leq t \leq 10 \,\text{s}}|\X_c^{N_{\text{coarse}}}-\X_c^{N_{\text{intermediate}}}|,\\
e_{\infty}[\X_c;N_{\text{intermediate}}]=\max\limits_{0 \,\text{s} \leq t \leq 10 \,\text{s}}|\X_c^{N_{\text{intermediate}}}-\X_c^{N_{\text{fine}}}|.
\end{align}
%~ Similarly, the empirical estimate for the convergence rate $q_{\infty}$ in the ${\Linf}$ norm of the soft disk's centroid $\X_c$ over the course of the simulation 
%~ time is defined as
%~ \begin{equation}
%~ q_{\infty}[\X_c;N_{\text{coarse}}]=\log_2\left(\frac{e_{\infty}[\X_c;N_{\text{coarse}}]}{e_{\infty}[\X_c;N_{\text{intermediate}}]}\right).
%~ \label{eq:conv_x}
%~ \end{equation}
The convergence rates of the centroid of the soft disk along with the volume of the disk
are reported in Fig.~\ref{fig:vol_x_dt_conv}(a), 
clearly showing second-order convergence in volume conservation and between first and second-order convergence
for the centroid's position. To further examine the contribution to the displacement error due to the time-integration algorithms in the 
FSI system,
we analyze the  temporal convergence of the disk's centroid at the final time $t=10$~s.
The grid-spacing is fixed in this case, corresponding to $N=3$ refinement levels. We use the solution from a very small time-step of $\dt=1\times{10}^{-5}$~s 
as the reference solution. As shown in Fig.~\ref{fig:vol_x_dt_conv}(b), full second-order convergence is evident.
This indicates that the reduction in the accuracy of the disk's centroid in Fig.~\ref{fig:vol_x_dt_conv}(a) can be attributed
to the spatial discretization error and not the temporal part.

As our last investigation for the soft disk moving in a lid driven cavity, 
we study the sensitivity and stability of our simulations for the choice of $\kappa$  and the impact of damping.
Once again, we consider the case with $N=2$ grid levels, $\nu=0.499$, and the same time step size as before 
to study the sensitivity of the numerical algorithm to penalty parameters based on changes in the spatio-temporal 
discretization.
Fig.~\ref{fig:kappa_eta_disk}(a) examines
 the sensitivity of the numerical algorithm to damping $\eta$ and its influence on the accuracy of the 
computed solution by varying $\kappa$ and $\eta$ to obtain a table demarcating cells with stable vs. unstable behavior.
The numbers reported for cells with stable pairs is the 
 maximum discrepancy between the two Lagrangian representations at time $t=10\, \text{s}$, shown as a factor of 
 $h \equiv h_\text{finest}$.
 Stable results are obtained for $\kappa \leq 570$,  without using any damping. The narrower range $200 \leq \kappa \leq 570$ includes 
 results with an acceptable accuracy based on the maximum tolerance, i.e.
 $\norm{\Y(\s,t)-\X(\s,t)}_{\infty} \leq 0.1 \, h_\text{finest}$.
 On the other hand, choosing $\eta > 0$ can reduce the displacement error and broaden the range of effective $\kappa$ value.
  Additionally, using a very large damping parameter $\eta$ leads to instability, regardless of the choice for $\kappa$.
     %~ Additionally, note that there can always be a small enough time-step for which the currently unstable pairs would become stable. 
Fig.~\ref{fig:kappa_eta_disk}(b) shows the maximum discrepancy between the two Lagrangian representations over time for three 
acceptable cases (maximum discrepancy equal or smaller than $0.1 \, h_\text{finest}$) from the table with a fixed
  $\kappa=300$ and different choices of $\eta$.
  Non-zero damping values are shown to further reduce spurious numerical oscillations.
 %~ Large values of stable $\kappa$ closer to the upper bound appear to need non-zero damping in order to maintain stability 
 %~ while stable $\kappa$ values 
 %~ closer to the lower bound do need seem to require additional damping.  
 Finally, Fig.~\ref{fig:trajectory_softdisk}(a) and Fig.~\ref{fig:trajectory_softdisk}(b) assess the sensitivity of the result to stable choices of $(\kappa,\eta)$, 
 the centroid's trajectory and the deformation at $t=10$~s for three distinctive and selected choices of $(\kappa,\eta)$.
  Results obtained for different choices of penalty parameters appear to be virtually identical and 
  	are largely insensitive to the values of $\kappa$ and $\eta$ within the acceptable region, i.e.,
  	the region with the maximum discrepancy equal or smaller than $0.1 \, h_\text{finest}$.

\subsection{Modified Turek-Hron benchmark \protect\footnote{This benchmark test is provided in IBAMR version 0.11.0 or newer, within the directory \texttt{examples/IIM/ex6}.}}
\label{subsubsec:Turek_Hron}

\begin{figure}[b!!!!]
		\centering
			\includegraphics[width=0.8\textwidth]{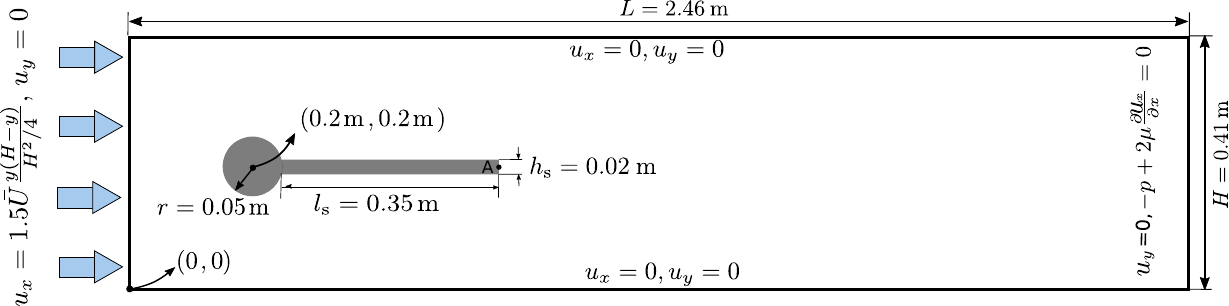}
		\caption{Schematic of the modified Turek-Hron benchmark including the geometry and boundary conditions. }	
		\label{fig:schematic-TH} 
\end{figure}

This test was originally proposed by Turek and Hron \cite{turek2006proposal} for benchmarking fluid-structure interaction algorithms.
We examine the performance of our methodology for models with
large added effects and large deformation resulting from the presence of a light structure.
The structure is composed of a stationary circular disk of radius $r=0.05 \, \text{m}$ centered at $(0.2 \, \text{m}, 0.2 \, \text{m})$
 and an elastic tail that has a length of $0.35 \, \text{m}$ and thickness $0.02 \, \text{m}$, see the schematic in Fig.~\ref{fig:schematic-TH}.
The rectangular domain surrounding the structure has a height $H=0.41$~m and length $L=6\,H=2.46$~m
 extending over $\Omega =[0.0 \, \text{m}, 2.46 \, \text{m}]\times[0.0 \, \text{m}, 0.41 \, \text{m}]$.
 A horizontal parabolic velocity profile of $u_x=1.5\,\bar{u}\,{y(y-H)}/(H^2/4)$ is imposed upstream of the structure in which $\bar{u}$ is the mean velocity.
Zero normal traction and zero tangential velocity are imposed at the outlet ($x=L$). Zero velocity conditions are imposed at $y=0.0$~m and $y=0.41$~m. 
The Eulerian domain is
 discretized using a refinement ratio of $r=4$ with $N=3$ Cartesian grid levels with a
  grid spacing of $h_{\text{coarsest}}=\frac{L}{24}$ on the coarsest 
  level.
Note that given the dimensions of the computational domain and the center position of the circular disk as in the original study,
the immersed body is positioned asymmetrically in the $y$-direction.
The circular disk is described by a nearly rigid material.
The left end of the tail is fixed and attached to the cylinder.
 We consider two cases with different density ratios, shear moduli, and mean entrance velocities.
The associated parameter sets including physical properties of both fluid and solid are reported in Table~\ref{table:TH_exp_info}.
In the original specification, a compressible St.~Venant-Kirchhoff model was used to describe the elastic tail with Poisson's ratio $\nu=0.4$,
but in this example we instead consider a nearly incompressible neo-Hookean material with Poisson's ratio $\nu=0.499$.
  The time step size is $\dt=(0.05\,\text{s}\cdot\text{m}^{-1})\,h_\text{finest}$, and the force penalty parameters are set 
to $\kappa=(2\times 10^{-2} \, \text{kg}\cdot\text{m}^{-2})/{\dt}^2$ and $\eta=(4\times 10^{-2} \, \text{kg}\cdot\text{m}^{-2})/\dt$ with $\Ncycle=4$ for case 1
 and $\Ncycle=8$ for case 2 described in Table~\ref{table:TH_exp_info}.
%~ to consistently enforce
%~ incompressibility on both solid and fluid in our formulation.
\begin{table}[t!!!!]
	\centering	
	\caption{Parameters for the modified Turek-Hron benchmark (Sec.~\ref{subsubsec:Turek_Hron}).}
	\label{table:TH_exp_info}	
\begin{tabular}{l*{6}{c}r}
\hline
Case  &   $ \bar{u} \, [\frac{\text{m}}{\text{s}}]$ & $\muf \, [\frac{\text{kg}}{\text{m}\cdot\text{s}}]$ &  $\rhof \, [10^3 \times \frac{\text{kg}}{\text{m}^{3}}]$ &  $\rhos_0 \, [10^3 \times \frac{\text{kg}}{\text{m}^{3}}]$  & $\Gs \, [10^6 \times \frac{\text{kg}}{\text{m}\cdot\text{s}^{2}}]$ & $\nu$ \\
\hline
1 (FSI2 of Turek \& Hron \cite{turek2006proposal})    & 1.0  & 1.0 & 1.0 & 10.0 & $0.5$  &  0.499  \\
2 (FSI3 of Turek \& Hron \cite{turek2006proposal})         & 2.0  & 1.0 & 1.0 & 1.0 & $2.0$ &  0.499
\end{tabular}
\end{table}
\begin{table}[t!!!]
	\centering	
	\caption{Comparison of computational values for the modified Turek-Hron benchmark
	for case 1 \& 2 described in Table~\ref{table:TH_exp_info}. Other simulation parameters include $h_\text{finest}=0.0293$, $\dt=(0.05\, \text{s}\cdot\text{m}^{-1})\,h_\text{finest}$, and $\mfac=2$.}
	\label{table:TH_comparison}
\begin{tabular}{l*{7}{c}r}
             & \hspace{1cm} Case 1 (FSI2 in \cite{turek2006proposal}) & & \hspace{1cm} Case 2 (FSI3 in \cite{turek2006proposal}) \\
\cmidrule(lr){2-3}\cmidrule(lr){4-5}
 & $y$-disp. of A $[10^{-3} \times m]$ & $\St$ & $y$-disp. of A $[10^{-3} \times m]$ & $\St$  \\
\hline
Turek \& Hron~\cite{turek2006proposal}         & $0.123 \pm 8.06$  &  0.19  &  $0.148 \pm 3.438$ & 0.265   \\
Zhang et al.~\cite{zhang2013three}      & $0.1 \pm 8.3$ & 0.192 & $0.1 \pm 3.6$   & 0.263  \\
Lee \& Griffith~\cite{lee2022lagrangian}      & - & - & $0.145 \pm 3.33$ & 0.25 \\
Present method    & $0.12 \pm 8.18$ & 0.193 &  $0.146 \pm 3.41$ & 0.265\\
\end{tabular}
\end{table}
\begin{figure}[t!!]
		\centering
			\includegraphics[width=0.95\textwidth]{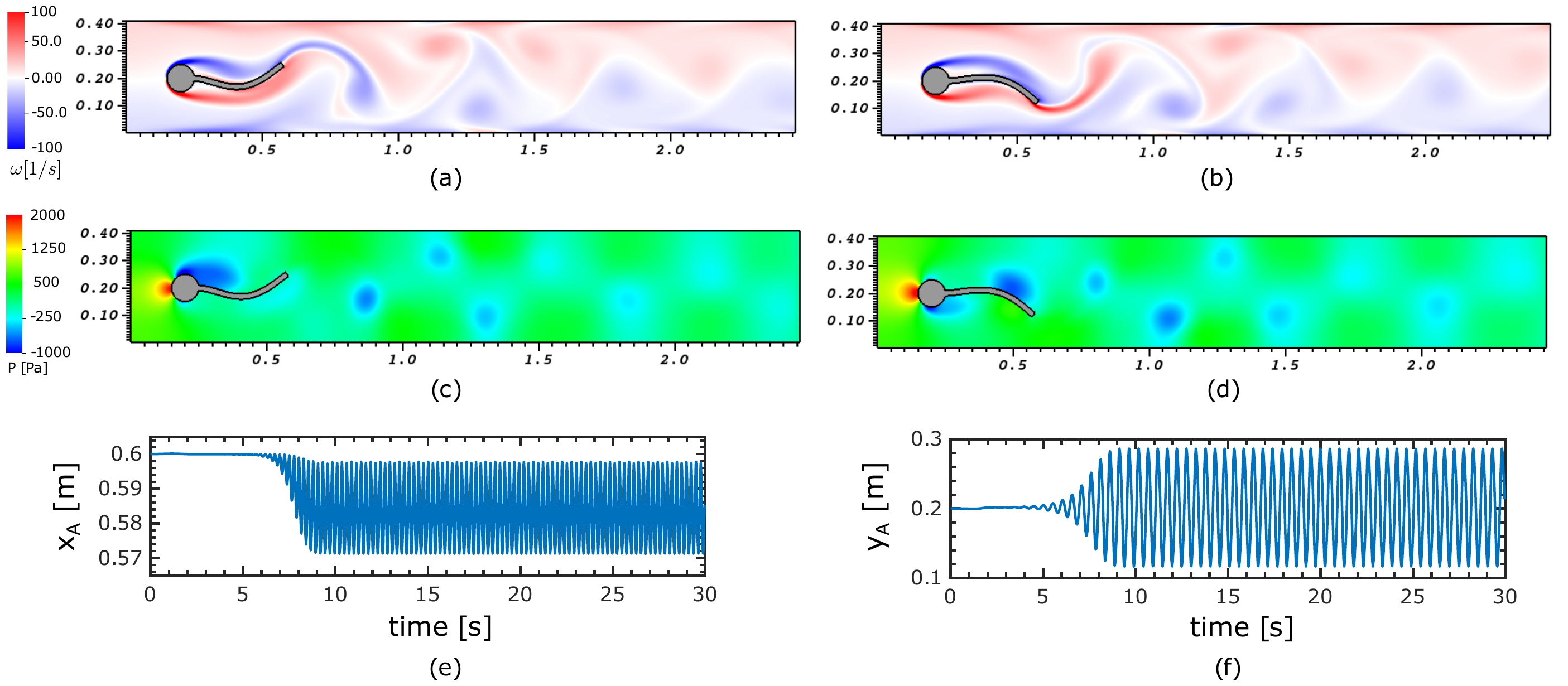}
		
		\caption{Dynamics of modified Turek-Hron benchmark corresponding to case 1 (FSI2 in \cite{turek2006proposal})
		(a) vorticity field at $t=28.4$~s, (b) vorticity field at $t=28.6$~s, (c) pressure field at $t=28.4$~s, (d) pressure field at $t=28.6$~s, 
		(e) time history of horizontal displacement of point A at the tail and (f) vertical displacement of point A at the tail.}	
		\label{fig:TH_RHO10} 
\end{figure}
%~ Tail_conv_Turek_Hron.pdf
Table~\ref{table:TH_comparison} compares the $y$-displacement of point A and the Strouhal number ($\St$) to previous studies 
\cite{turek2006proposal,bhardwaj2012benchmarking,zhang2013three,tian2014fluid,lee2022lagrangian} for the two cases 
reported in the original work by Turek \& Hron. 
The Strouhal number is $\St=f D/\bar{u}$, in which $f$ is the oscillation frequency for the $y$-displacement of point A 
and $D$ is the diameter of the circular disk.
We observe good agreement between values produced by the present method and previous
work, reflected also in the periodic vertical ($y_{\text{A}}$) and horizontal ($x_{\text{A}}$) positions of point A at the tip of the tail.
 Note that the study by Lee and Griffith \cite{lee2022lagrangian} uses an exactly incompressible material for the beam.
Fig.~\ref{fig:TH_RHO10} shows the dynamics of the beam and its deformation in the interaction with surrounding fluid
with density ratio $\rhos_0/\rhof=10$ and the Reynolds number $Re=100$. 
Vorticity magnitudes and 
pressure contours are respectively shown in panels (a)--(b)  and (c)--(d) 
for times after self induced periodic oscillations have developed in the structure and flow.
Both vorticity and pressure contours indicate the periodic behavior downstream of the flow.

\begin{figure}[b!!]
		\centering
			\includegraphics[width=0.9\textwidth]{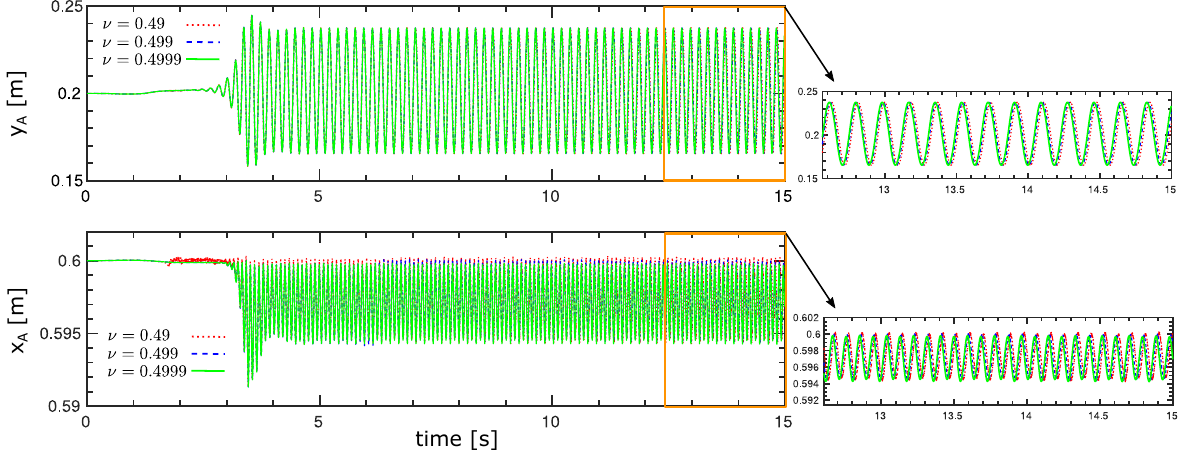}
		
		\caption{Time histories of the vertical (top) and horizontal (bottom) displacements of the flexible tail 
		in the modified Turek-Hron benchmark (case 2) for three different values of the Poisson's ratio $\nu$ approaching the fully incompressible limit $\nu \to 0.5$.}	
		\label{fig:TH_nu} 
\end{figure}

As an additional verification test, we have explored the sensitivity of our results in the limit of the Poisson's ratio approaching $\nu \to 5$. Similar to 
the soft disk example in Sec.~\ref{subsubsec:soft_disc}, we choose two additional Poisson's ratios, $\nu=0.49$ and $\nu=0.4999$ along 
with the default $\nu=0.499$ 
and study the tip of the tail's position (point A in Fig.~\ref{fig:schematic-TH}) over time. 
For the case of $\nu=0.4999$, we maintain stability for the overall scheme by using a smaller $\dts$, which is achieved by taking $\Ncycle=12$.
Notice that the fluid time step size $\dtf$ does not need to be reduced to accommodate this increase in the bulk modulus.
As shown in Fig.~\ref{fig:TH_nu}, there is very little difference 
between all the three cases and in particular, there is an 
excellent agreement 
between the tip's position for $\nu=0.499$  and $\nu=0.4999$. 
\begin{figure}[t!]
		\centering
			\includegraphics[width=0.9\textwidth]{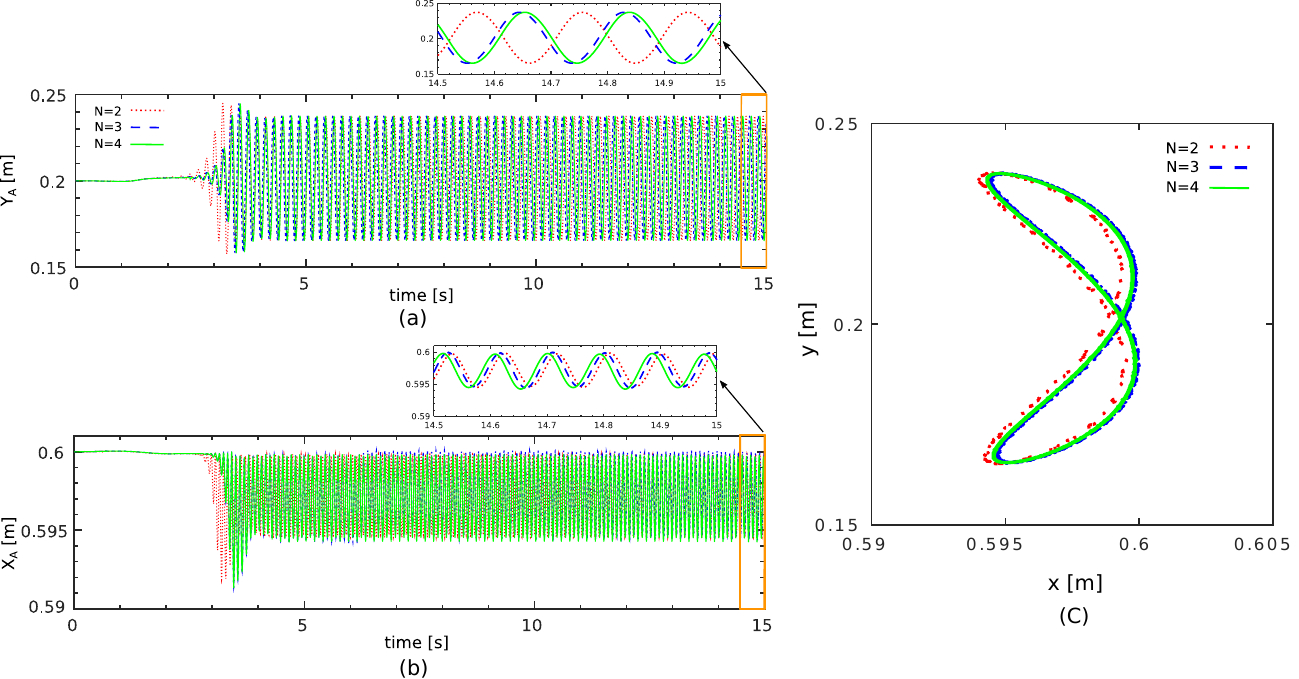}
		
		\caption{(a) Time histories of the vertical and (b) horizontal displacements of point A at the flexible tail 
		in the modified Turek-Hron benchmark (case 2) under grid refinement. (c) The trajectory of point A at the flexible tail 
		 under grid refinement for multiple cycles.}
		\label{fig:TH_conv} 
\end{figure}
\begin{figure}[b!!]
		\centering
			\includegraphics[width=0.75\textwidth]{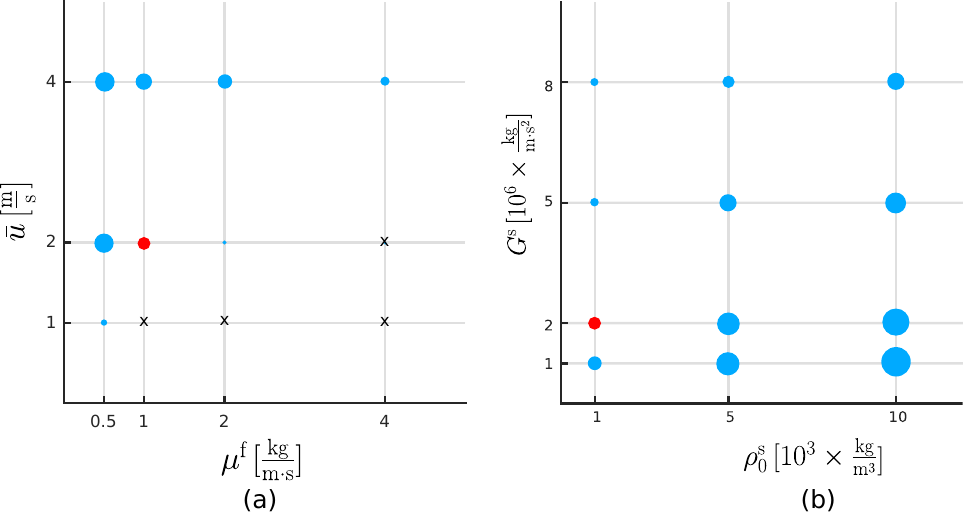}
		
		\caption{(a) The vertical oscillation amplitude of point A at the flexible tail shown for
		a range of  fluid viscosity $\muf$ and inlet mean velocity  $\bar{u}$ using case 2 of Table~\ref{table:TH_exp_info} as the default case 
		with fixed values of $\rhos_0 = 10^3 \, {\text{kg}}\cdot{\text{m}^{-3}}$ and  $\Gs = 2 \times 10^6 \,  {\text{kg}}\cdot{\text{m}^{-1}\cdot\text{s}^{-2}}$.
		The vertical amplitude sizes are depicted as circles of different diameters. The amplitude of $3.41 \times 10^{-3} \, \text{m}$
        given for case 2 of Table \ref{table:TH_comparison} is shown as a reference size for comparison. % corresponding to $ \bar{u}= 2 \, {\text{m}}\cdot{\text{s}}$ and $\muf = 1 \, {\text{kg}}\cdot{\text{m}^{-1}\cdot\text{s}^{-1}}$
		 (b) The vertical displacement of point A at the flexible tail shown for
		a range of  solid density $\rhos$ and shear modulus  $\Gs$ 
		using case 2 of Table~\ref{table:TH_exp_info} as the default case with fixed values of  $ \bar{u}= 2 \, {\text{m}}\cdot{\text{s}}$, $\rhof = 10^3 \, {\text{kg}}\cdot{\text{m}^{-3}}$
		and  $\muf = 1 \, {\text{kg}}\cdot{\text{m}^{-1}\cdot\text{s}^{-1}}$.
		As in panel (a), the amplitude of $3.41 \times 10^{-3} \, \text{m}$ corresponding to case 2 of Table \ref{table:TH_comparison} is shown as a reference size for comparison.}
		%~ The change in the solid density $\rhos$ 
		%~ can be interpreted as a change in the structure-to-fluid density ratio with the fixed fluid density considered here.
		\label{fig:TH_solid_fluid} 
\end{figure}

Next, we perform a grid refinement study for case 2 in Table \ref{table:TH_exp_info}.
 In our AMR framework, we vary the number of refinement levels by taking  
two additional grid spacings with levels $N=2$ and $N=4$ along 
with the default $N=3$. We use the same time step size and penalty parameters as before 
and examine the vertical and horizontal position of the tip of the tail (point A in Fig.~\ref{fig:schematic-TH}). 
Figs.~\ref{fig:TH_conv}(a) and (b) show the vertical and horizontal position of point A over time, respectively.
There is a good agreement across all grid spacings for the magnitudes of both vertical and horizontal variations. 
Additionally, there is a much smaller phase lag between the positions for $N=3$ and $N=4$ than the case with $N=2$.
 The same trend is also apparent in the convergence of the trajectory of point A under grid refinement, as shown
 in Fig.~\ref{fig:TH_conv}(c).

Finally, we study the effect of changes in the fluid dynamics parameters and solid structure properties 
on the vertical oscillation amplitude of point A at the flexible
tail.
To study the fluid dynamics parameters, we fix the solid properties based on case 2 in Table \ref{table:TH_exp_info} 
and vary the fluid viscosity $\muf$ and inlet mean velocity  $\bar{u}$.
Fig.~\ref{fig:TH_solid_fluid}(a) shows the vertical amplitude as circles with different diameters. The amplitude of $3.41 \times 10^{-3} \, \text{m}$
   associated with case 2 of Table \ref{table:TH_comparison} is shown as a reference size for comparison.
   Clearly, as we increase the viscosity the vertical amplitude decreases.
 Similarly, Fig.~\ref{fig:TH_solid_fluid}(b) indicates the vertical amplitudes of point A at the tail for a range 
 of shear modulus $\Gs$ and solid density $\rhos$  of the solid structure. 
 This time,  fluid dynamics properties are fixed corresponding to the reference point of case 2 in Table \ref{table:TH_comparison}.
 Increasing the shear modulus tends to damp the vertical amplitude of the periodic oscillation of the tail while larger density ratios result in larger 
 vertical amplitude.
  Regarding the influence of material properties on penalty parameters, our numerical experiments show that the penalty parameters are largely insensitive to changes in the shear modulus. Additionally, 
		 we noticed weak dependencies on the density ratio, fluid viscosity, and the mean inlet flow that are readily resolved by slightly adjusting the penalty parameters.
 %~ In particular, we aim to verify the adjustable tuning based on the functional relation of $\kappa$ previously given in Sec.~\ref{sec:kappa_scaling}. 
 %~ As previously explained in Sec.~\ref{sec:kappa_scaling}, we have verified that with our proper 
 %~ functional form of  $\kappa  \propto  \rhof h/(\sqrt{\Re} \, \Delta t^2)$ in terms of the Reynold number $\Re$,  no manual re-tuning is needed for 
 %~ the simulations reported in Fig.~\ref{fig:TH_solid_fluid}(a).
 %~ successfully achieved while 
 %~ Additionally, the changes in the solid properties in Fig.~\ref{fig:TH_solid_fluid}(b) are only reflected in adjusting the choice of $\Ncycle$ while
  %~ keeping $\kappa$ a fixed value. 

\subsection{Damped structural instability of a fully enclosed fluid reservoir \protect\footnote{This benchmark test is provided in IBAMR version 0.11.0 or newer, within the directory \texttt{examples/IIM/ex8}.}}
\label{subsubsec:fluid_reservoir}

  \begin{figure}[b!!]
		\centering
			\includegraphics[width=0.55\textwidth]{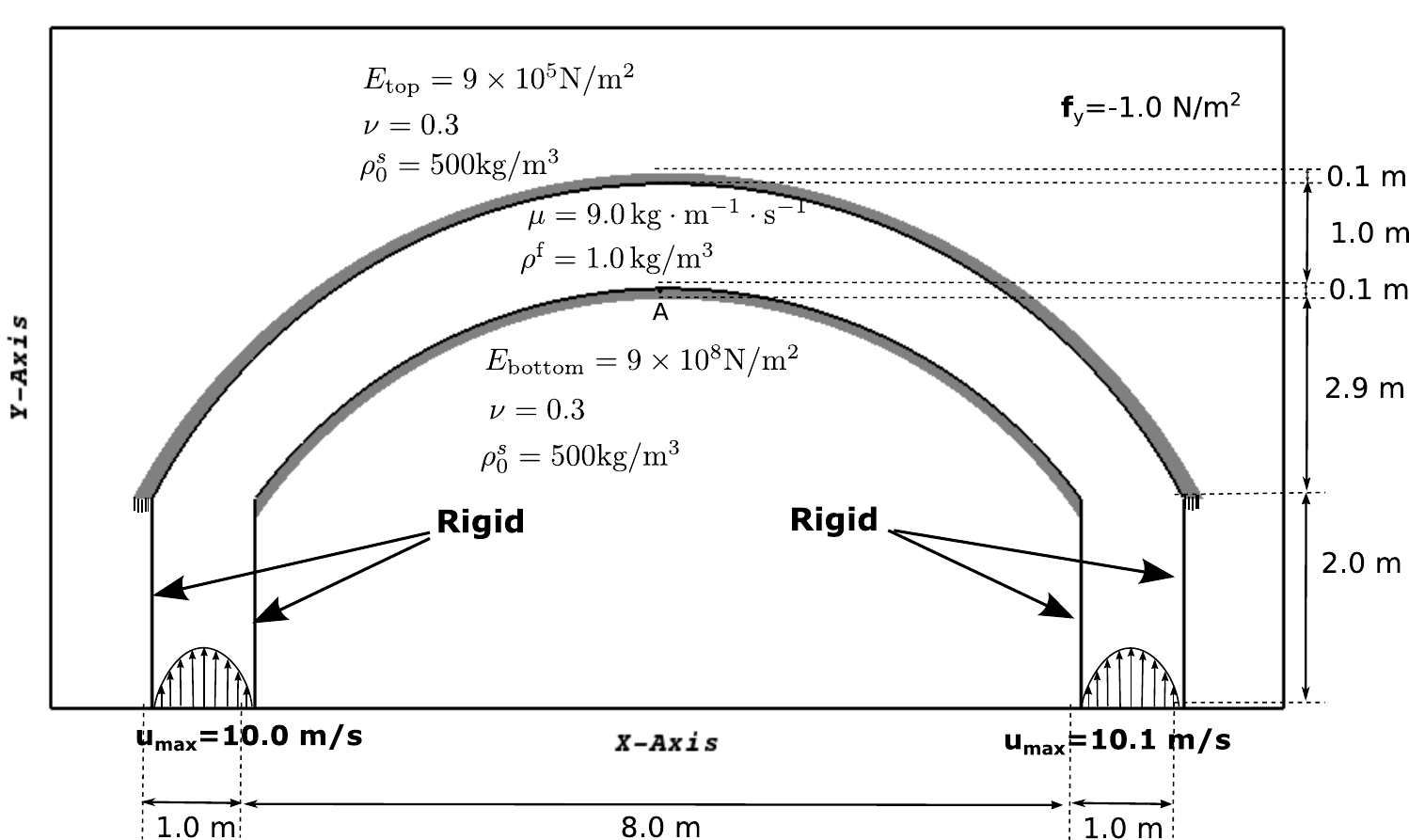}
		
		\caption{The problem description and boundary conditions of the fully enclosed fluid reservoir.} 
						\label{fig:damped_balloon_schematic} 
\end{figure}
  \begin{figure}[t!!]
		\centering
			\includegraphics[width=0.9\textwidth]{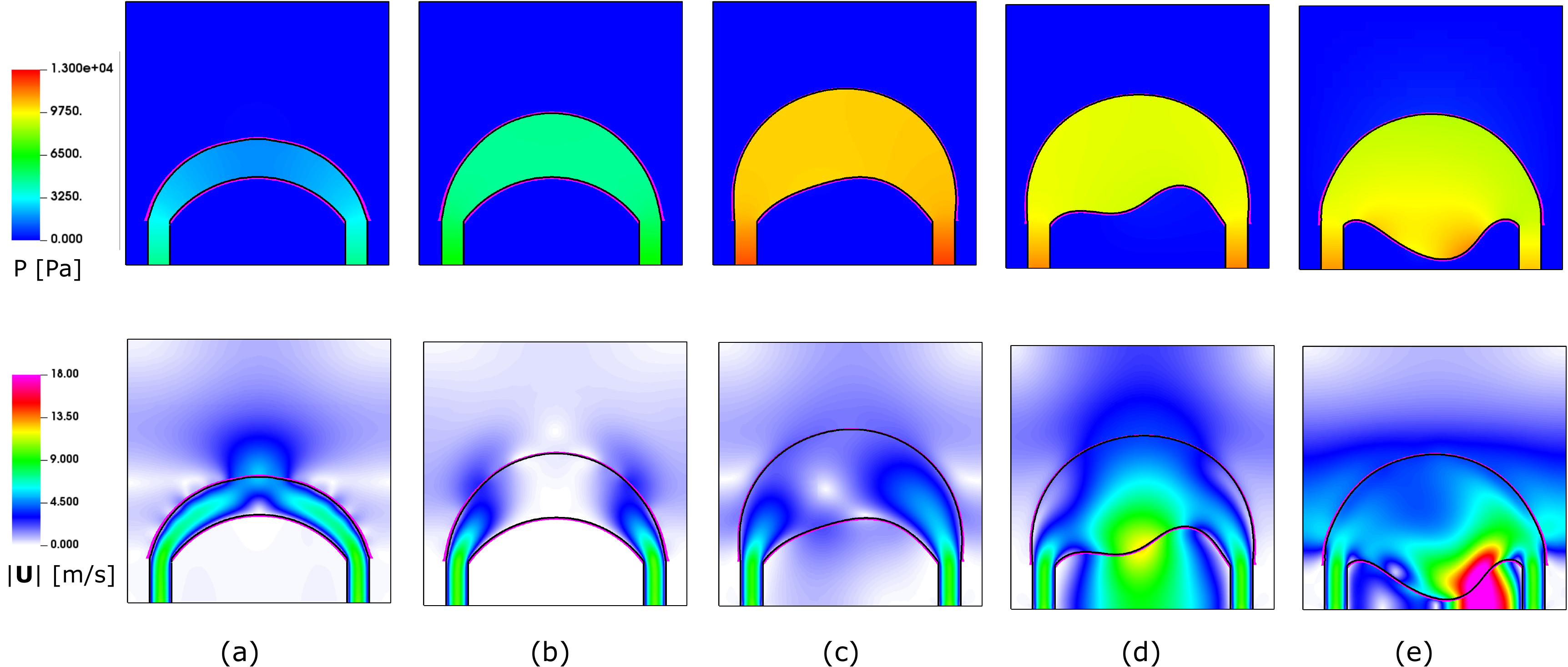}
		
		\caption{The contours of pressure (top) and velocity magnitudes (bottom) for the damped structural instability 
		of the enclosed fluid reservoir at times (a) $t=0.5$~s, (b) $t=1.2$~s, (c) $t=2.2$~s, (d) $t=2.4$~s, and (e) $t=2.6$~s. } 
				\label{fig:damped_balloon_contours} 
\end{figure}
\begin{figure}[b!!]
		\centering
			\includegraphics[width=0.4\textwidth]{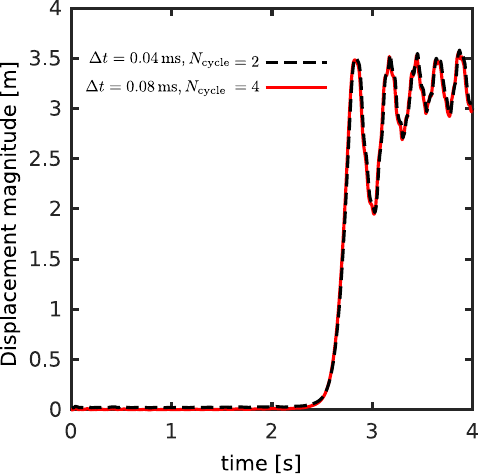}
		
		\caption{The displacement magnitude of the bottom band's mid-point over time. The sudden collapse of the lower band after some time is apparent.
		Results show two choices with two $\Ncycle$ values showing the multirate time-stepping capability of the algorithm in predicting essentially identical results.} 
						\label{fig:balloon_displacement} 
\end{figure}
To test our ILE formulation for compressible solid material interacting with an enclosed fluid domain, we consider a banded fluid domain that is surrounded by 
two rigid inlet channels and two thin curved structures with neo-Hookean material and different stiffness. 
 This benchmark was introduced by K{\"u}ttler et al.~\cite{kuttler2006solution}. Here we consider a version of the problem with connecting inlet channels 
 extended twice as long to ensure that the deformed configuration of the lower band will not fall outside of the computational domain.
The schematic of the system, boundary conditions and
 the dimensions are shown in Fig.~\ref{fig:damped_balloon_schematic}. 
The full computational domain is $\Omega =[-7.0 \, \text{m}, 7.0 \, \text{m}]\times[0.0 \, \text{m}, 6.0 \, \text{m}]$,
and is
 discretized using a refinement ratio of $r=4$ with $N=3$ Cartesian grid levels with a
  grid spacing of $h_{\text{coarsest}}=\frac{L_x}{16}$ on the coarsest 
  level. As usual, we use \textbf{Q1} elements for the solid structure with $\Mfac\approx 2$.
  The time step size is $\dt=0.04$~ms, and the force penalty parameters are set 
  to $\kappa=10^{6} \, \text{kg}\cdot\text{m}^{-2}\cdot \text{s}^{-2}$ and 
$\eta=2.0 \, \text{kg}\cdot\text{m}^{-2}\cdot \text{s}^{-1} $ with $\Ncycle=4$.
 At both inflow boundaries a fully developed velocity profile is prescribed 
 with the right inlet having a slightly larger maximum velocity than 
 the left inlet ($U_{\text{max}}=10.1\, \text{m}/\text{s}$ vs. $U_{\text{max}}=10 \, \text{m}/\text{s}$). 
 Zero normal traction and tangential velocity conditions are imposed at remaining sides of the embedding Cartesian domain.
 The fluid domain is loaded with the body force $\vec{f}=(0, -1)\, \text{N}/\text{m}^{2}$. 
 A Saint Venant-Kirchhoff constitutive model for the undamped shell is assumed 
 using Young moduli $E_{\text{top}}=9\times 10^5 \, \text{N}/\text{m}^2$ and $E_{\text{bottom}}=9\times 10^8\, \text{N}/\text{m}^2$ for the top and bottom
  bands, respectively. It is important to note that in this example, fluid-structure interaction only takes place between the inner fluid region within the reservoir and 
  the interior sides of flexible bands. Zero displacement conditions
   are imposed at the tails of both bands. Traction-free conditions are imposed on the exterior sides of both bands. 
   %~ ILE formulation allows for traction-free boundary condition by only coupling the interior side of the bands to the surrounding fluid. 
   In other words, although 
   the reservoir is embedded in a larger computational domain, in practice, the flexible structures do not ``feel'' the effect of the nonphysical 
   fluid region that is present in the computation but outside of the reservoir. 
   Because the FSI coupling only takes place on one side of the structures, this also makes it possible to solve for compressible
    structures, as detailed in Sec.~\ref{sec:penalty_ILE_comp}.

Because of the injected inflow at both inlets, the fluid pressure gradually increases inside the reservoir. It is expected that the softer top band first deforms 
to create space for the entering fluid. Once a critical pressure builds inside the reservoir, the stiffer bottom band collapses, and from that point forward 
the pressure rapidly increases in the domain. Fig.~\ref{fig:damped_balloon_contours} shows the pressure field and velocity 
magnitudes for multiple snapshots of the simulation. 
It is important to note that our ILE formulations naturally handles the solution for the fluid pressure within the reservoir 
as part of the solution in the larger computational domain. Typical partitioned formulations 
with Dirichlet-Neumann coupling require special treatments to handle an undetermined pressure that lies
in the null space of the pressure solver for a fluid region fully embedded within a solid structure \cite{kuttler2006solution,akbay2021boundary}.
Around time $t=2.2$~s, it is clear that the onset of the lower band collapsing is observed. This is concurrent with large pressure built-up in the reservoir prior to 
that shown in Fig.~\ref{fig:damped_balloon_contours}.
We observe that
despite minor differences in the problem setup, our results are qualitatively in good agreement with K{\"u}ttler et al.~\cite{kuttler2006solution},
 which uses an augmented Dirichlet-Neumann approach; 
specifically,  the collapse of the lower band and the structural instability have been captured in our results. 
We note that the time at which the bottom band starts to buckle was not explicitly reported in the original study by K{\"u}ttler et al.~\cite{kuttler2006solution}, 
but in a later work by
Fern{\'a}ndez et al. \cite{fernandez2015fully} using a fully decoupled time-marching scheme, the midpoint deformation of the lower band revealed a buckling time closer to $t=3$~s.  
A more recent study by Akbay et al. \cite{akbay2021boundary} using a 
partitioned boundary pressure projection method reported $t\approx 2$~s for the onset of the bottom band buckling. 

In addition, we have tested the multirate time-stepping feature of our ILE formulation for this problem. Fig.~\ref{fig:balloon_displacement} shows the 
magnitude of the displacement vector for a material point located at the mid-point of the lower flexible band over time. The first simulation uses 
a time step of $\dt=0.04$~ms and $\Ncycle=2$ as before, but in the second simulation, we have increased the fluid time step to $\dt=0.08$~ms 
while keeping the structural time step size fixed by using $\Ncycle=4$. The displacement magnitude of the midpoint at the lower flexible band is recorded over time 
and plotted for both simulations.
There is an excellent agreement between the cases. The sudden collapse of the lower band after inaction time of resistance against the pressure load is clearly observed.
 The maximum deformation of the mid-point around $3.5$~m appears to match the maximum deformation previously reported by Fern{\'a}ndez et al. \cite{fernandez2015fully}.

\subsection{Horizontal flexible plate inside a flow phantom \protect\footnote{This benchmark test is provided in IBAMR version 0.11.0 or newer, within the directory \texttt{examples/IIM/ex9}.}}
\label{subsubsec:Hessenthaler_model}

 \begin{figure}[t!!]
		\centering
			\includegraphics[width=0.6\textwidth]{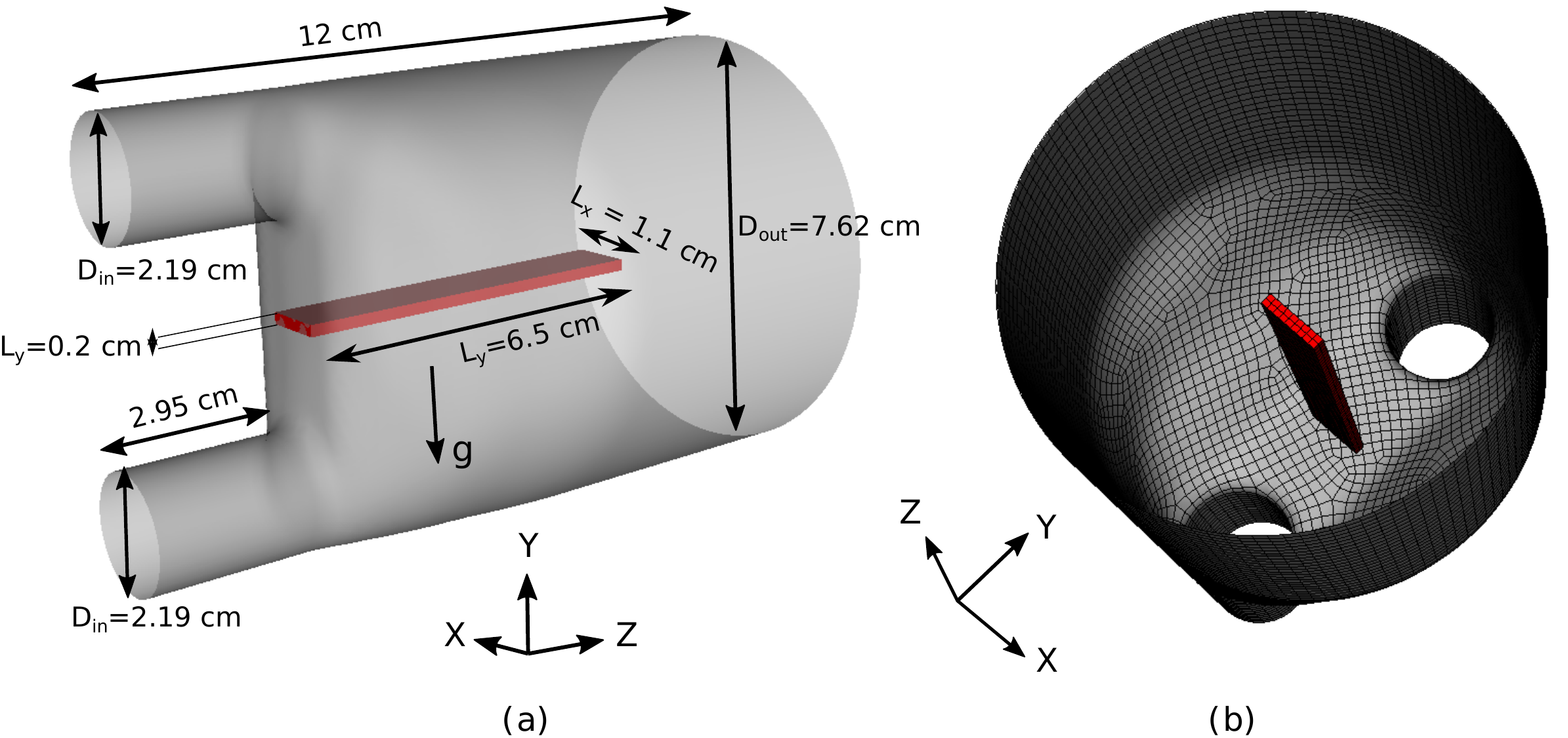}
		
		\caption{(a) Geometrical dimensions and (b) mesh for the flexible plate inside a flow phantom. 
		The gravitational force is directed towards the negative $y$-axis.}	
		\label{fig:geom-Hessenthaler} 
\end{figure}
 \begin{table}[b!!]
	\centering	
	\caption{Physical parameters for the horizontal flexible plate inside a flow phantom benchmark (Sec.~\ref{subsubsec:Hessenthaler_model}).  We use 
	a nearly incompressible material unlike the original work that uses $\nu=0.3151$.}
	\label{table:HT_exp_info}	
\begin{tabular}{l*{9}{c}r}
\hline
Case  & $ {\vec{u}}^{{\text{max}}} \, [\frac{\text{cm}}{\text{s}}]$ & $ {\vec{u}}^{{\text{bottom}}}_{\text{max}} \, [\frac{\text{cm}}{\text{s}}]$ &  $f_\text{r} \, [\frac{1}{\text{s}}]$ &  $\muf \, [\frac{\text{g}}{\text{cm}\cdot\text{s}}]$ & $\rhof \, [\frac{\text{g}}{\text{cm}^{3}}]$ & $\rhos_0 \, [\frac{\text{g}}{\text{cm}^{3}}]$ & $\Gs \, [\frac{\text{dyn}}{\text{cm}^{2}}]$ & $\nu$ \\
\hline
Hydrostatic &  (0, 0, 0) & (0, 0, 0) & - & 0.1250 &  1.1633  & 1.0583  & $6.1 \times 10^5$ & 0.499  \\
Phase I   &  (0, 0, 61.5) & (0, 0, 63.0) & - &  0.1250 & 1.1633  & 1.0583  & $6.1 \times 10^5$ & 0.499  \\
Phase II       & See Fig.~\ref{fig:Hessenthaler-phaseII}(a) & See Fig.~\ref{fig:Hessenthaler-phaseII}(a) & 1/6 & 0.1337 & 1.1640 & 1.0583 & $6.1 \times 10^5$ & 0.499
\end{tabular}
\end{table}
This section considers a three-dimensional FSI benchmark based on experimental results originally 
proposed by Hessenthaler et al.~\cite{hessenthaler2017experiment}.
This test involves an elastic beam clamped inside a flow phantom on one end and 
free on the other end moving inside the enclosure; see Fig.~\ref{fig:geom-Hessenthaler} for the geometrical dimensions of the setup.
 Two inlet pipes of diameter $D_{\text{in}} = 2.19 \, \text{cm}$ merge into one large outlet with diameter $D_{\text{out}}=7.62 \, \text{cm}$. 
 The housing is a stationary surface structure tethered in place by spring forces. 
 The filament (shown in red in Fig.~\ref{fig:geom-Hessenthaler}) has dimensions $L_x\times L_y \times L_z = 0.2 \, \text{cm} \times 1.1 \, \text{cm} \times 6.5$~cm.
  The computational domain is $\Omega =[-6.0 \, \text{cm}, 6.0 \, \text{cm}]\times[-6.0 \, \text{cm}, 6.0 \, \text{cm}]\times[0.0 \, \text{cm}, 12.0 \, \text{cm}]$,
 a cuboid of size $L_x \times L_y \times L_z = 12.0 \, \text{cm} \times 12.0 \, \text{cm} \times 12.0 \, \text{cm}$. 
 This Eulerian domain is
 discretized using a refinement ratio of $r=4$ with $N=3$ Cartesian grid levels 
 yielding a grid spacing of $h_{\text{coarsest}}=\frac{L_x}{12}$ on the coarsest 
  level. The time step size is $\dt=0.1$~ms with $\Ncycle=1$.
  The rigid housing is described using our underlying immersed interface formulation 
  for discrete surfaces \cite{kolahdouz2020immersed}  with a spring penalty parameter of 
  $\kappa = 3.50 \times 10^{5} \, \text{g}\cdot\text{cm}^{-2}\cdot \text{s}^{-2}$. 
  The penalty parameters for the filament are 
  $\kappa =3.75 \times 10^{4} \, \text{g}\cdot\text{cm}^{-2}\cdot \text{s}^{-2}$  and  $\eta = 5.0 \, \text{g}\cdot\text{cm}^{-2}\cdot \text{s}^{-1}$. 
   The physical parameters of the test are given in Table~\ref{table:HT_exp_info}. $f_\text{r}$ is the frequency of the applied velocity field.
  Note that in the original work the Poisson's ratio is $\nu=0.3151$, but here we use $\nu=0.499$ in 
 our nearly incompressible formulation.
 The numerical tests are performed in three stages: the hydrostatic; a steady phase I; and a transient phase II. 
 \begin{figure}[t!!]
		\centering
			\includegraphics[width=0.9\textwidth]{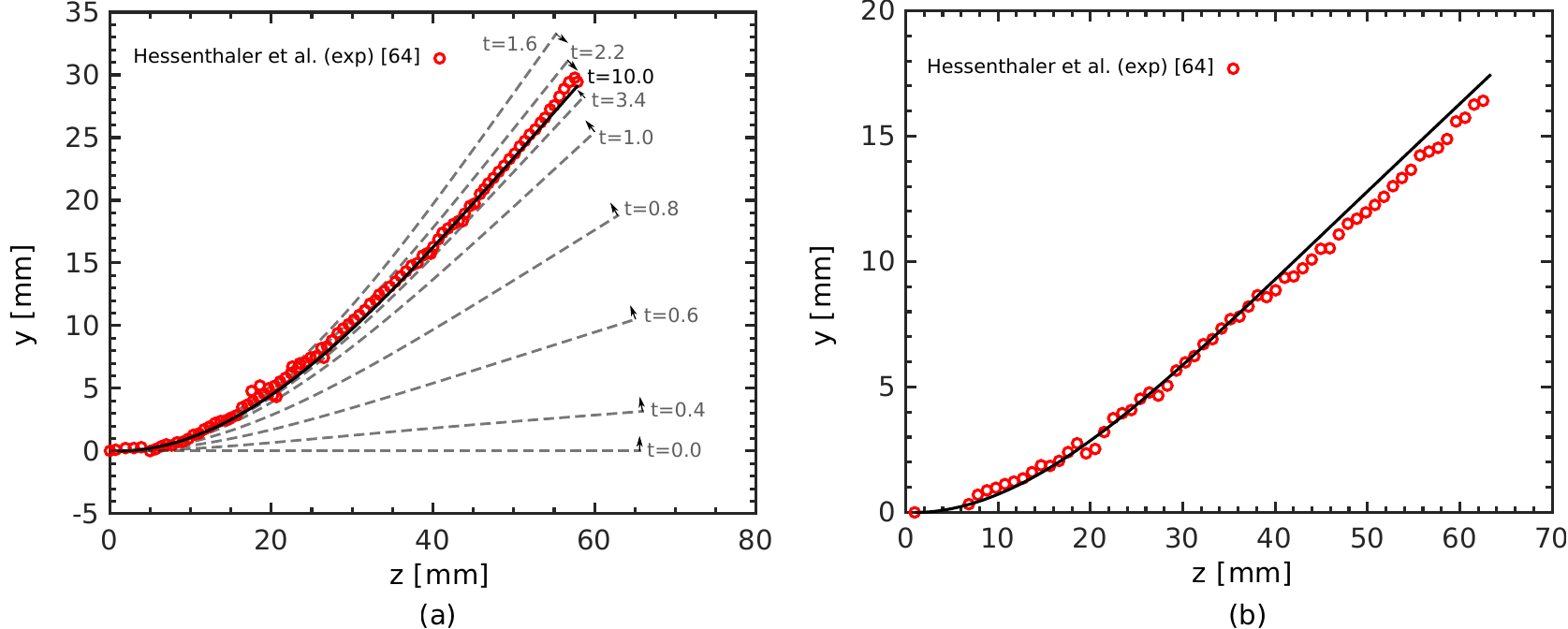}
		
		\caption{(a) The hydrostatic response for the $y$-position of the silicone filament centerline
		inside the flow phantom, and (b) the displacement of the silicone filament center-line in the Phase I experiment.
		 Red circles show measured experimental data by Hessenthaler et al.~\cite{hessenthaler2017experiment}.}	
		\label{fig:Hessenthaler_HS} 
\end{figure}
 \begin{figure}[b!!]
		\centering
			\includegraphics[width=0.7\textwidth]{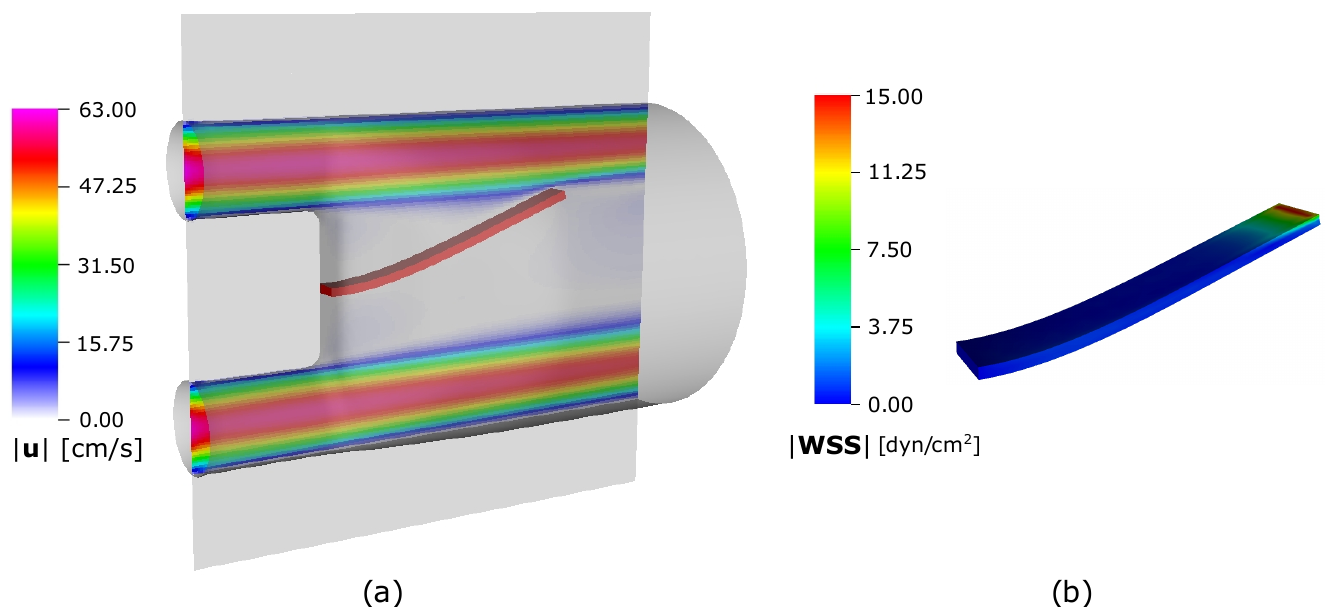}
		
		\caption{Simulation results of the horizontal flexible plate inside a flow phantom 
		including (a) steady-state velocity magnitudes at $x=0$  and (b) equilibrium position of the silicon beam including contour of steady-state 
		wall shear stress (WSS) magnitudes.}	
		\label{fig:Hessenthaler-phaseI} 
\end{figure}
 \begin{figure}[t!!]
		\centering
			\includegraphics[width=0.9\textwidth]{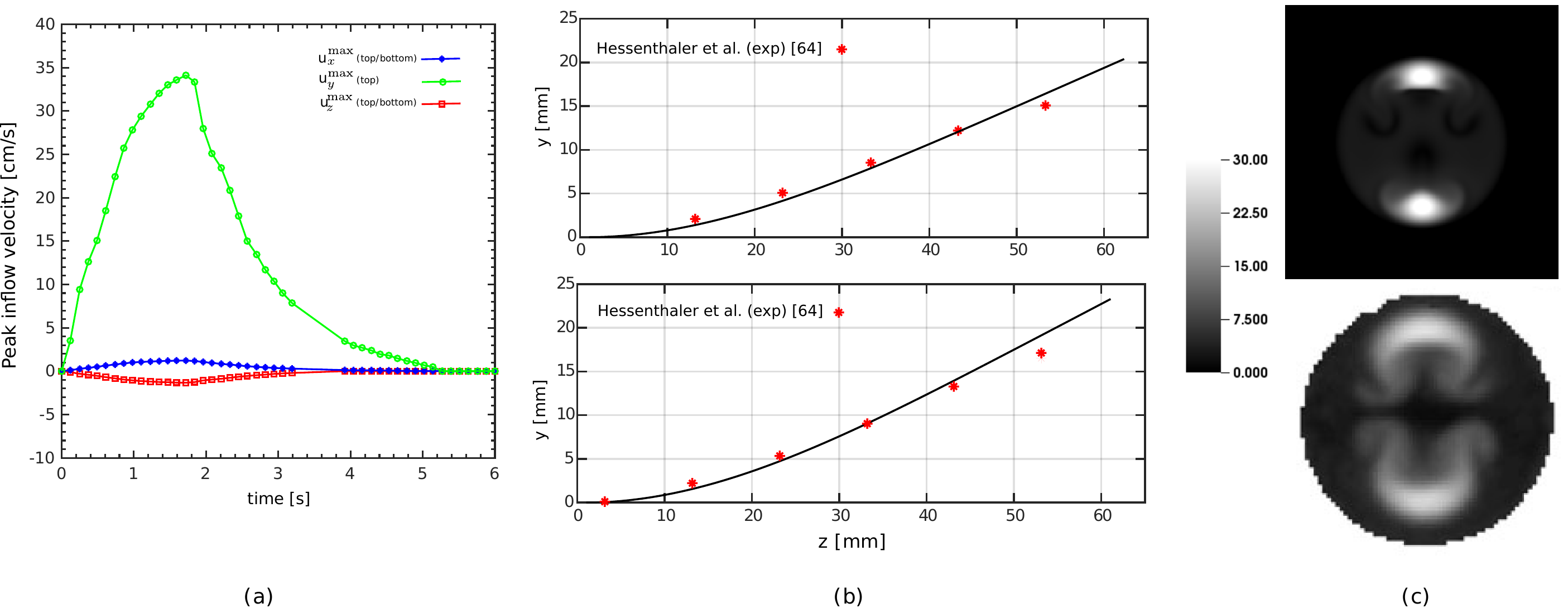}
		\caption{The Phase II benchmark of 
		the flexible plate inside a flow phantom:
		(a) Measured peak inflow velocity of the parabolic profiles 
		based on data from Hessenthaler et al.~\cite{hessenthaler2017experiment} used as boundary conditions of the top inlet. The bottom inlet uses the 
		same velocity profile as the top inlet except using $v_y^{\text{max}}=0$ for the y-component.
		(b) Sample results of the deformation of the flexible plate including 
		the centerline's $y$-position at $t=1.58$~s (top) and $t=2.88$~s (bottom). 
		Red stars show measured data by Hessenthaler et al.~\cite{hessenthaler2017experiment}.
		(c) $z$ velocity component magnitudes at $z=9.3$~cm at the time $t=2.02$~s for the present simulation (top) and 
		the MRI measurements \cite{hessenthaler2017experiment} (bottom).}	
		\label{fig:Hessenthaler-phaseII} 
\end{figure}

 Fig.~\ref{fig:Hessenthaler_HS}(a) shows the hydrostatic equilibrium solution under gravity with $g=980.665 \, \text{cm}/\text{s}^2$.
  The solution at $t=10$~s has reached the equilibrium state matching the measurements from the experiment. For the phase I experiment, fully developed 
  parabolic inflow boundary conditions are used with zero velocity in the $x$ and $y$ directions and peak velocity in the $z$ direction given as, 
  \begin{align}
   {{u}}_z^{{\text{max}}}(\text{top}) &=
      \begin{cases}
       61.5 \, (12 \, [\frac{\text{cm}}{\text{s}^3}] \, t^2 - 16 \, [\frac{\text{cm}}{\text{s}^4}] \, t^3)  &  t < 0.5\,[s], \\
        61.5 \,[\frac{\text{cm}}{\text{s}}] & t \geq 0.5\,[s],
     \end{cases} \\
   {{u}}_z^{{\text{max}}}(\text{bottom}) &=
      \begin{cases}
       63.0  \, (12 \, [\frac{\text{cm}}{\text{s}^3}] \, t^2 - 16 \, [\frac{\text{cm}}{\text{s}^4}] \, t^3)  &  t < 0.5\,[s], \\
        63.0 \, [\frac{\text{cm}}{\text{s}}] & t \geq 0.5\,[s].
     \end{cases} 
\end{align}
  Zero fluid traction is applied at the outflow and the boundary conditions at the remaining sides of the Eulerian domain's boundary is set to zero velocity.
  The simulation is run until a steady-state condition is reached. Fig.~\ref{fig:Hessenthaler_HS}(b) shows the steady-state solution for the
  displacement of the center-line of the silicon beam at $t=25$~s. Overall, there is a good agreement with the experimental data reported 
  by  Hessenthaler et al.~\cite{hessenthaler2017experiment}. Some discrepancies 
  could be attributed to the fact that the outflow condition is imposed at an outlet that is 
  not far enough from the filament. In Fig.~\ref{fig:Hessenthaler-phaseI}(a) the velocity magnitudes of the 
  steady-state solution at $x=0$ is shown along with the equilibrium position of the silicon beam 
  in Fig.~\ref{fig:Hessenthaler-phaseI}(b) colored by the magnitude of the steady-state wall shear stress. 
  
  To perform phase II of the benchmark, a pulsatile inflow, according to the peak velocity values in Fig.~\ref{fig:Hessenthaler-phaseII}(a), is 
  imposed  at the inlets
   except that the $y$ component of the velocity at the bottom inlet is set to $u_y=0$. The frequency of the pulsatile velocity is $f_r=1/6$ and we run the simulation for 
   two periods up to $t=12$~s. Fig.~\ref{fig:Hessenthaler-phaseII}(b) shows sample displacements of the filament's centerline for two separate times and their comparisons 
   with experimental data. Additional qualitative comparison of the magnitude of the $z$ component velocity with MRI measurement obtained by
   Hessenthaler et al. \cite{hessenthaler2017experiment} on $z=9.3$~cm 
   plane at time $t=2.02$~s is shown 
   Fig.~\ref{fig:Hessenthaler-phaseII}(c). Although some discrepancies can be observed in the details of the velocity distribution, 
   the overall magnitude and peak locations are generally in good agreement.

\subsection{Clot trapping of the FDA generic IVC filter}
\label{subsec:clot_ivc}

IVC filters are medical devices implanted in the IVC to capture blood clots from the lower extremities before they migrate to a patient's lungs and cause a pulmonary embolism.
Here we perform simulations to demonstrate 
the capability of using the present method to model the migration
 and trapping of a realistic flexible blood clot
by an IVC filter. 
The challenging aspects of this model result from the relatively large size of the clot that affects the local fluid dynamics, 
the large deformations 
that occur, and contact between the clot and the IVC filter.
In principle, the ILE formulation will naturally handle contact between immersed structures given sufficient spatial resolution. 
However, resolving contact for thin structures like the IVC filter is challenging because of the practical difficulties of resolving not just the 
thin struts, which have a width and thickness of approximately $300 \mu\text{m}$, but also the fluid boundary layers in the vicinity of the struts. 
For this case, 
we implemented a 
simple penalty-based contact model \cite{kamensky2015immersogeometric} to prevent interpenetration. This was done by locally 
adding opposing concentrated loads on nodes of 
the boundary of the volumetric solid structure where the fronts of two nearby structures come in close contact. In our computations, this force is 
added along with the interface forces of the exterior fluid in Eq.~(\ref{eq:weak_form3}). The weak formulation takes the form
\begin{equation}
\begin{split}
\rhos_0 \int_{\Omega^{\text{s}}_0} \left(\frac{\partial^2 \Y}{\partial t^2}(\s,t)\right) \, \cdot \, \vec{\uppsi}(\s) \, {\mathrm d} \s & = 
-\int_{\Omega^{\text{s}}_0} \, \PP(\s,t) \, \colon \,  \grad_{\s} \, \vec{\uppsi}(\s) \,{\mathrm d} \s \\
 &  + \int_{\Gamma^{\fs}_0} \,   j(\s,t) \, \vec{\tau}^{\text{f},\text{phys}}(\X(\s,t),t) \, \cdot \, \vec{\uppsi}(\s) \,\DA \\
 & + \int_{\Gamma^{\fs}_0}  \, \vec{F}^{\text{ct}}(\s,t) \, \cdot \, \vec{\uppsi}(\s) \,\DA.
\end{split}
\label{eq:weak_form_contact}
\end{equation}
in which $\vec{F}^{\text{ct}}(\s,t)$ is the surface force density obtained from the simple penalty-based contact model by Kamensky et al. 
\cite{kamensky2015immersogeometric}.
    \begin{figure}[t!!]
		\centering
			\includegraphics[width=0.6\textwidth]{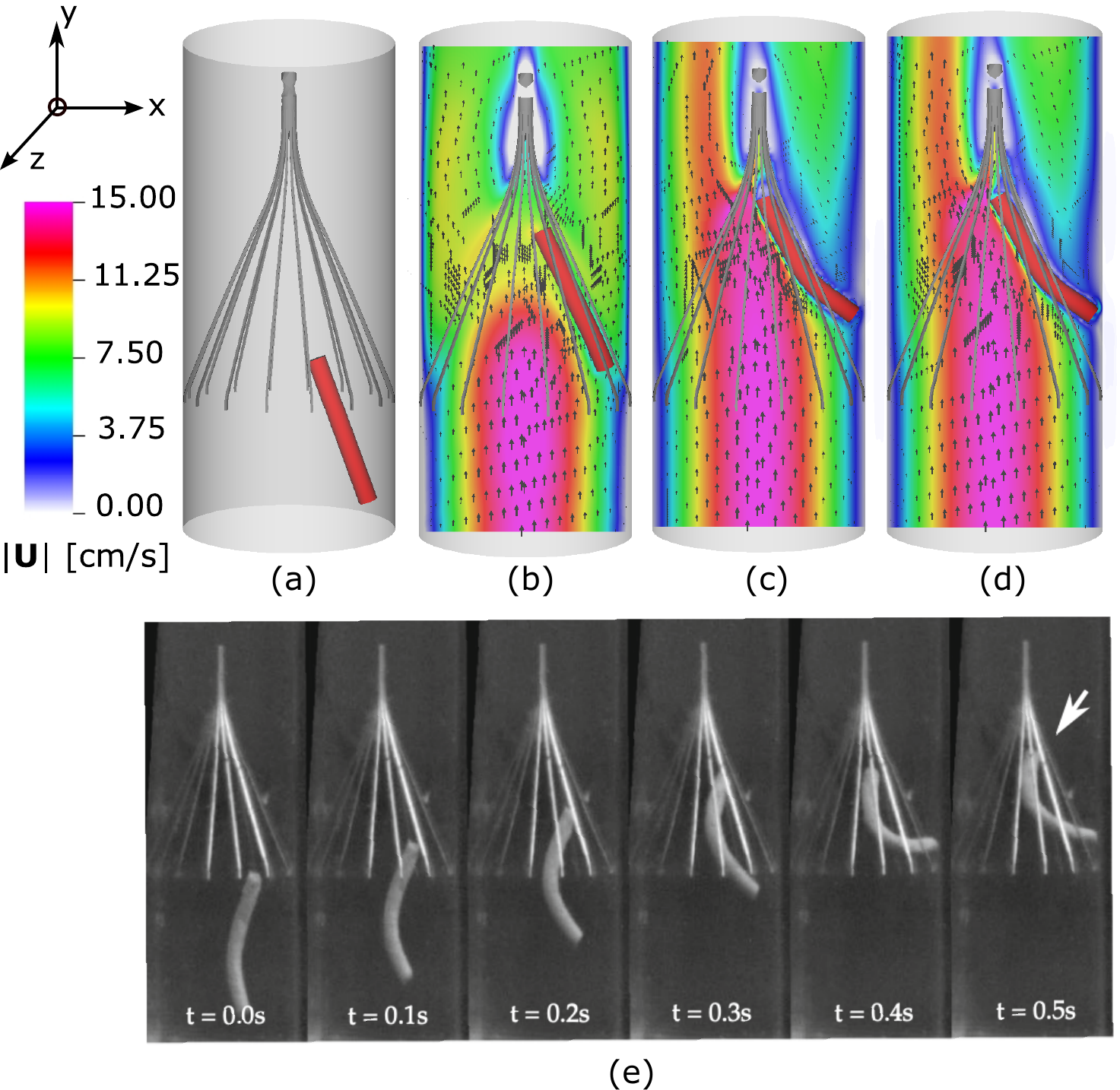}
		
		\caption{Simulations of the capturing of a geometrically realistic, deformable clot in the IVC filter 
		deployed in a rigid circular tube showing 
		(a) the initial clot-filter-tube configuration. 
		Dynamics of the clot deformation along with contours of velocity magnitude and vectors on the $x$-$z$ plane are shown for times:
		(b) $t=0.3$~s, (c) $t=1.45$~s, and (d) $t=1.8$~s. (e)
		Image sequence from high-speed videography of a similar 
		cylindrical blood clot trapped
		in the IVC filter \cite{riley2021vitro} showing similar clot trapping dynamics. }	
		\label{fig:filter_clot_cylinder} 
\end{figure}

The IVC filter used in this study is a 
new generic conical-type filter, referred to as the GENI (GEneric NItinol) filter, that is 
designed for research purposes by U.S. Food and Drug Administration (FDA) scientists and collaborators \cite{riley2021vitro}.
It resembles a cage with support structure of an umbrella consisting of a central hub with 16 evenly spaced radial struts.
We simulate the dynamics of clot transport and capture 
with the filter deployed in a rigid circular tube that is representative of the IVC.
Blood is modeled as a Newtonian fluid with the density 
$\rhof =1.0$~$\text{g}\cdot\text{cm}^{-3}$ and
dynamic viscosity $\muf = 0.03$~$\text{g}\cdot(\text{cm}\cdot \text{s})^{-1}$ .
A realistic cylindrical blood clot with length $2.0$~cm and diameter $2.5$~mm is considered.
The blood clot is modeled as a nearly incompressible hyperelastic material
 with Poisson's ratio $\nu=0.499$, shear modulus $\Gs=2000~\text{dyn}\cdot\text{cm}^{-2}$, and density $\rhos_0=1.2\,\rhof$.
 Here we ignore the influence of gravity.
 The clot and the filter are deployed in a circular tube of length $6.4$~cm 
  and diameter $2.8$~cm. 
  A fully developed parabolic velocity profile with the maximum velocity of 
  $U_{\text{max}}=15.0 \, \text{cm}\cdot\text{s}^{-1}$ is imposed at the 
inlet of the tube at the bottom boundary ($y=-3.2$~cm), corresponding to the total flow rate of $2.77\, \text{L}/\text{min}$.
Zero normal traction and tangential velocity conditions are imposed at the outlet of the tube at the top boundary. No-slip conditions 
are imposed along all the 
other parts of the Cartesian domain boundary.
 The computational domain is 
 a cuboid of size $L_x \times L_y \times L_z = 3.2 \, \text{cm} \times 6.4 \, \text{cm} \times 3.2 \, \text{cm}$. 
 This Eulerian domain is
 discretized using mixed refinement ratios, including two Cartesian grid levels with a refinement 
 ratio of four and additional three Cartesian grid levels with
  refinement ratios of two, resulting effective Cartesian grid spacings of $h_{\text{coarsest}}=\frac{3.2}{12} \, \text{cm}$ on the coarsest 
  level and $h_{\text{finest}}=\frac{3.2}{12\times4\times2^3} \, \text{cm}$ 
  on the finest level. 
  Note that while it is not necessary to use a mixed refinement ratio,
our numerical experiments show that in some large scale simulations a mixed refinement ratio provides a
reasonable compromise between the overall computational cost and additional refinement in regions where higher accuracy is desired.
  Both the IVC filter and the circular tube are assumed to be rigid, stationary structures 
that are modeled using penalty methods \cite{kolahdouz2020immersed}.
The time step size is $\dt=0.04$~ms and 
and the force penalty parameters for the clot are 
 set 
  to $\kappa=10^{5} \, \text{g}\cdot\text{cm}^{-2}\cdot \text{s}^{-2}$ and 
$\eta=1.0 \, \text{g}\cdot\text{cm}^{-2}\cdot \text{s}^{-1} $ with $\Ncycle=4$.

Fig.~\ref{fig:filter_clot_cylinder}(a) shows the initial setup of the problem including the geometries of the clot, filter, and the circular tube.
  The clot is positioned at an initial angle of $\theta_0=20^{\circ}$ with the $y$ axis, and its center of mass is initially located at $(0.7\,\text{cm}, \, 2.15 \, \text{cm}, \, 0.0 \, \text{cm})$.

In the present methodology, the dynamics of the clot deformation along with the velocity distribution in the x-z mid-plane are shown in Fig.~\ref{fig:filter_clot_cylinder}(b)--(d) at three 
different times. 
The clot is captured by the filter at $t \approx 0.35$~s.
  As the simulation proceeds, the clot continues to contact the filter while undergoing deformations
  caused as a result of the clot interaction with the local fluid field.
The dynamics of the clot's 
deformation follow similar patterns as experimental observations 
from recent in vitro experiments using the same model IVC device \cite{riley2021vitro}.

\section{Summary and conclusions}
\label{sec:conclusion}

This paper introduces 
a partitioned sharp interface immersed approach
for modeling flexible structures interacting with a Newtonian fluid.
It extends our recently developed rigid body immersed Lagrangian-Eulerian method for FSI \cite{kolahdouz2021sharp} to models involving immersed flexible structures.
The partitioned aspect of our formulation provides the flexibility of taking
a broad range of solid to fluid density ratios and 
using multi-rate time-stepping 
that allows for 
separate choices of time-step sizes in the fluid and
 solid mechanics solvers. 
 These features are in addition to general advantages of immersed formulations, including geometrical flexibility
 and fast domain solution capability for large deformation FSI.
The Dirichlet-Neumann coupling strategy only connects fluid and solid subproblems through interface conditions, unlike 
alternative fictitious domain approaches that use  
distributed Lagrange multiplier to hold the constraint in the entire space occupied by the structure.
The deformations of the structure are determined 
as a result of interplay between the forces of intrinsic stresses associated with the solid constitutive model 
and the exterior fluid traction forces that are computed from the fluid solver using physical jump conditions and imposed as boundary conditions to the solid domain.
We then use the deformation of the flexible structure to determine the motion
 of the fluid along the fluid-structure interface via a penalty method.
The dynamics of the volumetric structural mesh are determined using a standard finite element approach to
large-deformation nonlinear elasticity via a nearly incompressible solid mechanics formulation.
We also demonstrate the method's capability for modeling compressible solid structures
 if the fluid subdomain is completely enclosed within the structure.
 Pointwise second-order accuracy of the difference norm
between the positions of the two representations of the fluid-structure interface 
as well as a second-order accuracy
 for the structural volume of the soft disk in a lid driven cavity are achieved.
 For the structural position of the disk's centroid, between first and second-order spatio-temporal convergence is observed in the
 $\Linf$ norm over the entire simulation.
 To assess and validate the robustness and accuracy of the new algorithm,
comparisons are made with 
computational and experimental FSI benchmarks in two and three dimensions. 
Through our benchmark test cases, we demonstrate test cases with different density ratios including smaller, equal and larger than one. 
 We previously reported  the ability of the method in modeling problems with various density ratios for rigid-body FSI \cite{kolahdouz2021sharp}.
  Though rigorous stability analyses and analytical investigations 
 still remain to be done and are beyond the aim of the current extension of the method,
 our numerical experiments herein provided further evidence that the ILE formulation can handle flexible-body FSI problems with large added mass effects.
A rigorous mathematical analysis of the stability of the methodology could be the subject of future investigation.
Finally, we show representative results from applications of this methodology
 to the transport and capture of a geometrically realistic, deformable blood clot in an inferior vena cava filter.

Future work could address further applications of this method to more challenging problems in biomedical fluid-structure interaction. 
Numerical approaches that use fully incompressible structural formulations, including a more sophisticated exactly incompressible and implicit elasticity 
formulation, can be developed.
Additionally, the current formulation requires the structure to have a finite thickness. 
A surface representation with higher regularity can be considered in future work that allows for imposing higher order jump conditions, more accurate traction calculations on the solid boundary and 
higher order structural deformation. Further extensions on the fluid solver for more accurate simulation of large Reynolds number phenomena
could consider higher-order spatial discretizations and large-eddy simulation (LES) turbulence modeling.

\section*{Acknowledgements}

We acknowledge research support through 
NIH (Awards U01HL143336 and R01HL157631) and NSF (DMS 1664645, CBET 175193, OAC 1652541, and OAC 1931516) 
and the U.S.~FDA Center for Devices and Radiological Health (CDRH) Critical Path program.
This research was supported in part by an appointment to the Research Participation Program at the U.S.~FDA administered by 
the Oak Ridge Institute for Science and Education through an interagency agreement between the U.S.~Department of Energy and FDA.
The findings and conclusions in this article have not been formally disseminated by the FDA and should not be construed to represent
 any agency determination or policy.
The mention of commercial products, their sources, or their use in connection with material reported herein is not to be construed as either an actual or implied endorsement of such products by the Department of Health and Human Services.
Computations were performed using facilities provided by University of North Carolina at Chapel Hill through the 
Research Computing division of UNC Information Technology Services.

%~ \clearpage

\bibliography{paper-ILE-flexible}

\end{document}